\documentclass[a4j,12pt]{article}

\usepackage{amsmath,amsthm,amssymb}
\usepackage{amsmath,amssymb}
\usepackage{amscd}
\usepackage{cases}
\usepackage{comment}
\usepackage[dvipdfmx]{graphicx}
\usepackage{epic,eepic}
\usepackage[small]{caption}

\allowdisplaybreaks[1]

\theoremstyle{plain}
\newtheorem{thm}{Theorem}[section]
\newtheorem{lem}[thm]{Lemma}
\newtheorem{prop}[thm]{Proposition}
\newtheorem{cor}[thm]{Corollary}

\newtheorem{mthm}{Main Theorem}[section]

\theoremstyle{definition}
\newtheorem{df}[thm]{Definition}
\newtheorem{rem}[thm]{Remark}
\newtheorem*{rem*}{Remark}
\newtheorem{ex}[thm]{Example}

\makeatletter

\@addtoreset{equation}{section}
\makeatother
\newcommand{\R}{\mathbb{R}}
\newcommand{\C}{\mathbb{C}}
\newcommand{\Z}{\mathbb{Z}}

\newcommand{\Aut}{\mathop{\mathrm{Aut}}\nolimits}
\newcommand{\ord}{\mathop{\mathrm{ord}}\nolimits}
\newcommand{\Graph}{\mathop{\mathrm{Graph}}\nolimits}

\newcommand{\Tr}{\mathop{\mathrm{Tr}}\nolimits}
\newcommand{\divisor}{\mathop{\mathrm{div}}\nolimits}
\newcommand{\reduced}{\mathrm{red}}
\newcommand{\Br}{\mathbf{Br}}
\newcommand{\br}{\mathbf{br}}
\newcommand{\barkarrow}{\xrightarrow{\: \mathrm{bark} \:}}

\newcommand{\keywords}[1]{{\renewcommand{\thefootnote}{\fnsymbol{footnote}}\footnote[0]{\textit{Key words and phrases.} #1 }}}
\newcommand{\subclass}[1]{{\renewcommand{\thefootnote}{\fnsymbol{footnote}}\footnote[0]{\textit{MSC Classification.} #1 }}}
\newcommand{\email}[1]{{\renewcommand{\thefootnote}{\fnsymbol{footnote}}\footnote[0]{e-mail: {\text #1}}}}

\title{\textsf{Singular Fibers in Barking Families 
of Degenerations of Elliptic Curves}}

\author{Takayuki OKUDA}

\date{}

\begin{document}

\maketitle

\begin{abstract}
Takamura \cite{Ta3} established 
a theory of splitting families of degenerations of complex curves of genus 
$g \geq 1$. 
He introduced a powerful method for constructing a splitting family, 
called a \emph{barking family}, 
in which the resulting family of complex curves has 
a singular fiber over the origin (the \emph{main fiber}) 
together with other singular fibers (\emph{subordinate fibers}). 
He made a list of barking families for genera up to 
$5$ 
and determined the main fibers appearing in them. 
This paper determines \emph{most} of the subordinate fibers of 
the barking families in Takamura's list for the case 
$g = 1$. 
(There remain four undetermined cases.) 
Also, we show that 
some splittings never occur in a splitting family. 
\subclass{Primary 14D06; Secondary 14H15, 14D05, 32S50.}
\end{abstract}

\section{Introduction}  \label{sec:Intro}

Let 
$\pi : M \to \Delta$ 
be a proper surjective holomorphic map 
from a smooth complex surface 
$M$ 
to an open disk 
$\Delta 
:= \left\{s \in \C \, : \,  \left| s \right| < \delta\right\}$ 
in 
$\C$ 
with radius 
$\delta > 0$. 
We call 
$\pi : M \to \Delta$ 
a \emph{family of complex curves} of genus 
$g \geq 1$ 
over 
$\Delta$ 
if 
$\pi$ 
has at most finitely many singular fibers 
and the other fibers are smooth complex curves of genus 
$g$.
In particular, 
$\pi : M \to \Delta$ 
is called 
a \emph{degeneration} of complex curves of genus 
$g$ 
if the fiber 
$X_0 
:= \pi^{-1}(0)$ 
over the origin is singular 
and the other fibers 
$X_s 
:= \pi^{-1}(s) \, (s \neq 0)$ 
are all smooth.

In this paper, 
we consider the following problem: 
\emph{How does a singular fiber split in a deformation?} 
Let us recall the concept of a splitting family of degenerations. 
Let 
$\mathcal{M}$ 
be a smooth complex $3$-manifold 
and set 
$\Delta^{\dagger} 
:= \left\{t \in \C \, : \, \left| t \right| < \varepsilon\right\}$, 
an open disk with sufficiently small radius 
$\varepsilon > 0$. 
Consider a proper flat surjective holomorphic map 
$\Psi : \mathcal{M} \to \Delta \times \Delta^{\dagger}$. 
For 
$t \in \Delta^{\dagger}$, 
set 
$\Delta_t := \Delta \times \left\{t\right\}$, 
$M_t := \Psi^{-1}(\Delta_t)$ 
and 
$\pi_t := \Psi \big|_{M_t} : M_t \to \Delta_t$. 
Suppose that 
$\pi_0 : M_0 \to \Delta_0$ 
coincides with a given degeneration 
$\pi : M \to \Delta$.
Then we call 
$\Psi : \mathcal{M} \to \Delta \times \Delta^{\dagger}$ 
a \emph{deformation family} of the degeneration 
$\pi : M \to \Delta$ 
and each 
$\pi_t : M_t \to \Delta_t 
\, (t \in \Delta^{\dagger} \setminus \left\{0\right\})$ 
a \emph{deformation} of the degeneration 
$\pi : M \to \Delta$. 
In particular, 
$\Psi : \mathcal{M} \to \Delta \times \Delta^{\dagger}$ 
is called a \emph{splitting family} 
if every deformation 
$\pi_t : M_t \to \Delta_t$ 
of the degeneration 
$\pi : M \to \Delta$ 
is a family of complex curves with at least two singular fibers. 
Set 
$X_{s, t} 
:= \Psi^{-1}(s, t) \, 
(= \pi_t^{-1}(s))$. 
Clearly 
$X_{0, 0}$ 
is the original singular fiber 
$X_0$ 
of the degeneration 
$\pi : M \to \Delta$. 
For a fixed 
$t \in \Delta^{\dagger} \setminus \left\{0\right\}$, 
let 
$s_1, s_2, \dots, s_N \, (N \geq 2)$ 
be the singular values of 
$\pi_t$, 
that is, 
$X_{s_1, t}, X_{s_2, t}, \dots, X_{s_N, t}$ 
are the singular fibers of 
$\pi_t : M_t \to \Delta_t$. 
(note: 
The singular values 
$s_1, s_2, \dots, s_N$ 
depend on 
$t$, 
but the number of them and the types of the singular fibers do not.) 
In this case, 
we say that 
the singular fiber 
$X_0$ 
\emph{splits} into the singular fibers 
$X_{s_1, t}, X_{s_2, t}, \dots, X_{s_N, t}$.

To classify \emph{atomic degenerations} --- 
degenerations admitting no splitting family --- 
Takamura \cite{Ta3} introduced a powerful method 
for constructing splitting families. 
Splitting families obtained by this construction 
are called \emph{barking families}. 
In a barking family, 
the original singular fiber 
$X_0$ 
of the degeneration 
$\pi : M \to \Delta$ 
is deformed to a simpler singular fiber of its deformation 
$\pi_t : M_t \to \Delta_t$ 
in such a way that a part of 
$X_0$ 
looks ``barked'' off from 
$X_0$. 
See Figure \ref{ta:fig:BarkingExample} in Section \ref{sec:Ta}. 
The resulting singular fiber appears over the origin of 
$\Delta_t$ 
under Takamura's construction, 
so we denote it by 
$X_{0, t}$. 
In such a situation, 
we write%
\footnote{
In the same situation, Takamura \cite{Ta3} wrote 
$X_0 \longrightarrow\ X_{0, t}$. 
In this paper, 
we use 
``$\longrightarrow$'' 
only for splittings and distinguish 
``$\barkarrow$'' 
from it. 
} 
$$
X_0 \barkarrow X_{0, t}, 
$$
and call 
$X_{0, t}$ 
the \emph{main fiber}. 

In \cite{Ta3}, 
for genera up to $5$,
Takamura made a list of barking families which enabled him to show that 
a degeneration is \emph{absolutely atomic} --- that is, 
any degeneration topologically equivalent to it is atomic --- if and only if 
its singular fiber is either a Lefschetz fiber or a multiple of a smooth complex curve. 
For instance, 
he listed thirty five barking families 
for degenerations of complex curves of genus 
$g = 1$, 
that is, 
for degenerations of elliptic curves, 
and determined the type of the main fiber of each of them 
as follows, 
where 
we use Kodaira's notation%
\footnote{
See Table \ref{ta:table:KodairaNotation} in Section \ref{sec:Ta}. 
} 
for singular fibers 
(see also the list in Section \ref{sec:Talist}): 
\begin{align}  \label{eq:talist}
&\quad \textrm{\bf Takamura's list}\rm\\
&\boldsymbol{[II.1]} \quad II \barkarrow I_1 &\hspace{30pt}
 &\boldsymbol{[III^*.5]} \quad III^* \barkarrow I_6 \nonumber\\
&\boldsymbol{[II.2]} \quad II \barkarrow I_1 &\hspace{30pt}
 &\boldsymbol{[III^*.6]} \quad III^* \barkarrow I_2^* \nonumber\\
&\boldsymbol{[II^*.1]} \quad II^* \barkarrow III^* &\hspace{30pt}
 &\boldsymbol{[III^*.7]} \quad III^* \barkarrow I_7 \nonumber\\
&\boldsymbol{[II^*.2]} \quad II^* \barkarrow IV^* &\hspace{30pt}
 &\boldsymbol{[III^*.8]} \quad III^* \barkarrow I_6 \nonumber\\
&\boldsymbol{[II^*.3]} \quad II^* \barkarrow I_2^* &\hspace{30pt}
 &\boldsymbol{[III^*.9]} \quad III^* \barkarrow IV^* \nonumber\\
&\boldsymbol{[II^*.4]} \quad II^* \barkarrow I_5 &\hspace{30pt}
 &\boldsymbol{[IV.1]} \quad IV \barkarrow I_3 \nonumber\\
&\boldsymbol{[II^*.5]} \quad II^* \barkarrow I_3^* &\hspace{30pt}
 &\boldsymbol{[IV.2]} \quad IV \barkarrow I_2 \nonumber\\
&\boldsymbol{[II^*.6]} \quad II^* \barkarrow I_3^* &\hspace{30pt}
 &\boldsymbol{[IV.3]} \quad IV \barkarrow I_2 \nonumber\\
&\boldsymbol{[II^*.7]} \quad II^* \barkarrow I_8 &\hspace{30pt}
 &\boldsymbol{[IV.4]} \quad IV \barkarrow II \nonumber\\
&\boldsymbol{[II^*.8]} \quad II^* \barkarrow III^* &\hspace{30pt}
 &\boldsymbol{[IV^*.1]} \quad IV^* \barkarrow I_1^* \nonumber\\
&\boldsymbol{[II^*.9]} \quad II^* \barkarrow III^* &\hspace{30pt}
 &\boldsymbol{[IV^*.2]} \quad IV^* \barkarrow I_0^* \nonumber\\
&\boldsymbol{[III.1]} \quad III \barkarrow I_2 &\hspace{30pt}
 &\boldsymbol{[IV^*.3]} \quad IV^* \barkarrow I_6 \nonumber\\
&\boldsymbol{[III.2]} \quad III \barkarrow I_1 &\hspace{30pt}
 &\boldsymbol{[IV^*.4]} \quad IV^* \barkarrow I_1^* \nonumber\\
&\boldsymbol{[III.3]} \quad III \barkarrow I_2 &\hspace{30pt}
 &\boldsymbol{[I_0^*.1]} \quad I_0^* \barkarrow I_4 \nonumber\\
&\boldsymbol{[III^*.1]} \quad III^* \barkarrow IV^* &\hspace{30pt}
 &\boldsymbol{[I_0^*.2]} \quad I_0^* \barkarrow I_3 \nonumber\\
&\boldsymbol{[III^*.2]} \quad III^* \barkarrow I_1^* &\hspace{30pt}
 &\boldsymbol{[I_n^*.1]} \quad I_n^* \barkarrow I_{n-1}^* \nonumber\\
&\boldsymbol{[III^*.3]} \quad III^* \barkarrow I_2^* &\hspace{30pt}
 &\boldsymbol{[I_n^*.2]} \quad I_n^* \barkarrow I_{n+4}. \nonumber\\
&\boldsymbol{[III^*.4]} \quad III^* \barkarrow I_0^* &\hspace{30pt}
 \nonumber
\end{align}

In a barking family, 
there appear not only the main fiber but also other singular fibers, 
which are called \emph{subordinate fibers}. 
In what follows, 
when the original singular fiber 
$X_0$ 
splits into the main fiber 
$X_{0, t}$ 
and subordinate fibers 
$X_{s_1, t}, X_{s_2, t}, \dots, X_{s_N, t} \, (s_i \neq 0)$, 
we write 
$$
X_0 \longrightarrow X_{0, t}\ +\ X_{s_1, t} + X_{s_2, t} + \dots + X_{s_N, t}
$$ 
--- 
\emph{we always put the main fiber 
$X_{0, t}$ 
on the initial term 
to distinguish it from the subordinate fibers}. 
The main fiber of a barking family is explicitly described. 
On the other hand, 
it is not clear what subordinate fibers will appear. 
The aim of this paper is 
to determine the subordinate fibers of Takamura's barking families 
for degenerations of elliptic curves.

Our results are summarized in two theorems. 
Firstly, 
the following theorem determines 
the subordinate fibers of most of the barking families in the above list 
(note: four cases remain undetermined, 
see Remark \ref{rem:introexceptionalsplitting} below):

\begin{mthm}[Theorem \ref{thm:SplittingThm}]  \label{intro:mthm:SplittingThm}
Each barking family in Takamura's list $(\ref{eq:talist})$ 
except 
$\boldsymbol{[III^*.8]}$, 
$\boldsymbol{[IV.3]}$, 
$\boldsymbol{[IV.4]}$, 
$\boldsymbol{[I_0^*.2]}$ 
splits the singular fiber as follows: 
\begin{align*}\rm
&\boldsymbol{[II.1]} \ II \longrightarrow I_1 + I_1 &
 &\boldsymbol{[III^*.2]} \ III^* \longrightarrow I_1^* + I_2 \\
&\boldsymbol{[II.2]} \ II \longrightarrow I_1 + I_1 &
 &\boldsymbol{[III^*.3]} \ III^* \longrightarrow I_2^* + I_1 \nonumber\\
&\boldsymbol{[II^*.1]} \ II^* \longrightarrow III^* + I_1 &
 &\boldsymbol{[III^*.4]} \ III^* \longrightarrow I_0^* + I_1 + I_1 + I_1 \nonumber\\
&\boldsymbol{[II^*.2]} \ II^* \longrightarrow IV^* + II &
 &\boldsymbol{[III^*.5]} \ III^* \longrightarrow I_6 + I_1 + I_1 + I_1 \nonumber\\
&\boldsymbol{[II^*.3]} \ II^* \longrightarrow I_2^* + I_1 + I_1 &
 &\boldsymbol{[III^*.6]} \ III^* \longrightarrow I_2^* + I_1 \nonumber\\
&\boldsymbol{[II^*.4]} \ II^* \longrightarrow I_5 &
 &\boldsymbol{[III^*.7]} \ III^* \longrightarrow I_7 + I_1 + I_1 \nonumber\\
&\hspace{35pt} + I_1 + I_1 + I_1 + I_1 + I_1 &
 &\boldsymbol{[III^*.9]} \ III^* \longrightarrow IV^* + I_1 \nonumber\\
&\boldsymbol{[II^*.5]} \ II^* \longrightarrow I_3^* + I_1 &
 &\boldsymbol{[IV.1]} \ IV \longrightarrow I_3 + I_1 \nonumber\\
&\boldsymbol{[II^*.6]} \ II^* \longrightarrow I_3^* + I_1 &
 &\boldsymbol{[IV.2]} \ IV \longrightarrow I_2 + I_1 + I_1 \nonumber\\
&\boldsymbol{[II^*.7]} \ II^* \longrightarrow I_8 + I_1 + I_1 &
 &\boldsymbol{[IV^*.1]} \ IV^* \longrightarrow I_1^* + I_1 \nonumber\\
&\boldsymbol{[II^*.8]} \ II^* \longrightarrow III^* + I_1 &
 &\boldsymbol{[IV^*.2]} \ IV^* \longrightarrow I_0^* + I_1 + I_1 \nonumber\\
&\boldsymbol{[II^*.9]} \ II^* \longrightarrow III^* + I_1 &
 &\boldsymbol{[IV^*.3]} \ IV^* \longrightarrow I_6 + I_1 + I_1 \nonumber\\
&\boldsymbol{[III.1]} \ III \longrightarrow I_2 + I_1 &
 &\boldsymbol{[IV^*.4]} \ IV^* \longrightarrow I_1^* + I_1 \nonumber\\
&\boldsymbol{[III.2]} \ III \longrightarrow I_1 + I_2 &
 &\boldsymbol{[I_0^*.1]} \ I_0^* \longrightarrow I_4 + I_1 + I_1 \nonumber\\
&\boldsymbol{[III.3]} \ III \longrightarrow I_2 + I_1 &
 &\boldsymbol{[I_n^*.1]} \ I_n^* \longrightarrow I_{n-1}^* + I_1 \nonumber\\
&\boldsymbol{[III^*.1]} \ III^* \longrightarrow IV^* + I_1 &
 &\boldsymbol{[I_n^*.2]} \ I_n^* \longrightarrow I_{n+4} + I_1 + I_1. \nonumber
\end{align*}
\end{mthm}

\begin{rem}  \label{rem:introexceptionalsplitting}
We have not been able to determine the subordinate fibers
of the four exceptional barking families 
$\boldsymbol{[III^*.8]}$, 
$\boldsymbol{[IV.3]}$, 
$\boldsymbol{[IV.4]}$, 
$\boldsymbol{[I_0^*.2]}$ 
(see also Remark \ref{rem:UndeterminedCases}): 
\begin{align*}
&\boldsymbol{[III^*.8]} \quad III^* \longrightarrow 
I_6 + II + I_1, \ I_6 + I_2 + I_1,\ \mathrm{or}\ I_6 + I_1 + I_1 + I_1 \\
&\boldsymbol{[IV.3]} \quad IV \longrightarrow 
I_2 + II, \ \mathrm{or}\ I_2 + I_1 + I_1 \\
&\boldsymbol{[IV.4]} \quad IV \longrightarrow 
II + II, \ II + I_2, \ \mathrm{or}\ II + I_1 + I_1 \\
&\boldsymbol{[I_0^*.2]} \quad I_0^* \longrightarrow 
I_3 + II + I_1, \ \mathrm{or}\ I_3 + I_1 + I_1 + I_1. \\
\end{align*}
\end{rem}

In contrast, 
there are splittings that never occur in a splitting family. 
If in a splitting family for a degeneration of elliptic curves 
the singular fiber 
$X_0$ 
splits into 
$N$ 
singular fibers 
$X_1, X_2, \dots, X_N$, 
then 
we have 
$e(X_0) = e(X_1) + e(X_2) + \dots + e(X_N)$, 
where 
$e(X_i)$ 
denotes the topological Euler characteristic 
of the underlying reduced curve of 
$X_i$ 
(Lemma \ref{det:lem:SplittingInv} (b)). 
However the converse does not hold. 
Even if the singular fibers satisfy this equation, 
the splitting 
$X_0 \longrightarrow X_1 + X_2 + \dots + X_N$ 
does not always occur. 
In fact:

\begin{mthm}[Theorem \ref{thm:NonsplittingThm}]  \label{intro:mthm:UnsplittingThm}
None of the following splittings occurs: 
\begin{align*}
IV \longrightarrow \ 
&I_2 + I_2, \\
II^* \longrightarrow \ 
&I_8 + II, \quad I_7 + III, \quad I_6 + IV, \\
&I_4 + I_0^*, \quad I_3 + I_1^*, \quad \\
&I_u + I_v \, (u + v = 10), \\
III^* \longrightarrow \ 
&I_7 + II, \quad I_6 + III, \quad I_5 + IV, \\
&I_3 + I_0^*, \quad I_u + I_v \, (u + v = 9), \\
IV^* \longrightarrow \ 
&I_6 + II, \quad I_5 + III, \quad I_4 + IV, \\
&I_2 + I_0^*, \quad I_u + I_v \, (u + v = 8), \\
I_n^* \, (n \geq 0) &\longrightarrow \
I_{n + 4} + II, \quad I_{n + 3} + III, \quad I_{n + 2} + IV, \\
&I_u + I_v \, (u + v = n + 6 \textrm{ and } (n, u, v) \neq (2, 4, 4)). \\
I_0^* \longrightarrow \ 
&I_3 + I_2 + I_1. 
\end{align*}
\end{mthm}

\subsection*{Organization of this paper}

This paper is organized as follows. 
In Section \ref{sec:Ta}, 
we first review Takamura's theory of barking families, 
mainly for degenerations with \emph{stellar} (star-shaped) singular fibers. 
In fact, 
most of the degenerations of elliptic curves 
may be assumed to have stellar singular fibers.

To determine the subordinate fibers of the barking families 
in Takamura's list $(\ref{eq:talist})$, 
we investigate the singular fibers in three steps: 
(1) 
In Section $\ref{sec:Euler}$, 
we first consider the Euler characteristics of the singular fibers 
and give a list of the sets of subordinate fibers 
that can appear in each of the barking families. 
(2) 
In Section $\ref{sec:Monodromies}$, 
we recall the concept of monodromies around singular fibers, 
and 
in Section $\ref{sec:Nonsplitting}$, 
by comparing the traces of monodromies, 
we prove Main Theorem \ref{intro:mthm:UnsplittingThm} 
--- 
we give a list of splittings that never occur. 
In Section $\ref{sec:Det1}$, 
based on the result of Section $\ref{sec:Nonsplitting}$, 
we determine the subordinate fibers of five of Takamura's barking families. 
(3) 
Sections $\ref{sec:PropSing}$, $\ref{sec:CoreSing}$, $\ref{sec:Number}$ 
are devoted to study of the singularities of subordinate fibers. 
We investigate 
the singularities near proportional subbranches in Section $\ref{sec:PropSing}$ and 
those near the core in Section $\ref{sec:CoreSing}$. 
In Section $\ref{sec:Number}$, 
we show useful lemmas 
which give us the number of the subordinate fibers 
and that of their singularities. 
In Section $\ref{sec:Det2}$, 
we determine the subordinate fibers of the remaining barking families, 
and complete the proof of Main Theorem \ref{intro:mthm:SplittingThm}.

In Section $\ref{sec:Decomps}$, 
we give monodromy decompositions 
corresponding to the splittings induced from 
Takamura's barking families.

In Section $\ref{sec:Talist}$, 
we provide Takamura's list of barking families for genus 
$1$ 
with figures of the singular fibers, 
which will help the reader comprehend the barking deformations.

\subsection*{Acknowledgements}

The author would like to express his deep gratitude to 
Shigeru Takamura 
for fruitful discussions. 
The author would also like to thank 
Masaaki Ue and Osamu Saeki 
for helpful comments and suggestions. 
The author also thanks 
the anonymous referee 
for giving very useful comments and warm encouragement 
on an earlier version of the paper.

\section{Takamura's theory}  \label{sec:Ta}

Let us review Takamura's theory of barking families.
For details see \cite{Ta3}.

First we recall the concept of linear degenerations.
We begin with preparation.
Let 
$\pi : M \to \Delta$ 
be a degeneration of complex curves of genus 
$g \geq 1$ 
and express its singular fiber as 
$X_0 = \sum_i m_i \Theta_i$, 
where 
$\Theta_i$ 
is an irreducible component of 
$X_0$ 
with multiplicity 
$m_i$. 
In what follows, 
we assume that 
the \emph{underlying reduced curve} 
$X_0^{\reduced} 
:= \sum_i \Theta_i$ 
of 
$X_0$ 
has at most \emph{simple normal crossings}, 
that is, 
(i) 
any singularity of 
$X_0^{\reduced}$ 
is a node and 
(ii) 
any irreducible component 
$\Theta_i$ 
is not self-intersecting 
(so 
$\Theta_i$ 
is smooth).

For an irreducible component 
$\Theta_i$ 
of 
$X_0$, 
we denote by 
$N_i$ 
the normal bundle of 
$\Theta_i$ 
in 
$M$. 
Let 
$\left\{p_i^{(1)}, p_i^{(2)}, \dots p_i^{(h)}\right\}$ 
be the set of the intersection points on 
$\Theta_i$ 
with other irreducible components of 
$X_0$ 
and 
$m^{(j)} \, (j = 1, 2, \dots, h)$ 
be the multiplicity of the irreducible component intersecting 
$\Theta_i$ 
at 
$p_i^{(j)}$. 
Then there exists a holomorphic section 
$\sigma_i$ 
of the line bundle 
$N_i^{\otimes (-m_i)}$ 
on 
$\Theta_i$ 
such that
$$
\divisor(\sigma_i) = \sum_{j=1}^h m^{(j)} p_i^{(j)},
$$
where 
$\divisor(\sigma_i)$ 
denotes the divisor defined by 
$\sigma_i$.
Here 
$\sigma_i$ 
has a zero of order 
$m^{(j)}$ 
at 
$p_i^{(j)}$. 
Note that 
$\sigma_i$ 
is uniquely determined up to multiplication by a constant. 
We call 
$\sigma_i$ 
the \emph{standard section} of 
$N_i^{\otimes (-m_i)}$ 
on 
$\Theta_i$.

Take an open covering 
$\Theta_i = \bigcup_{\alpha} U_{\alpha}$ 
such that 
$U_{\alpha} \times \C$ 
is a local trivialization of the normal bundle 
$N_i$ 
on 
$\Theta_i$. 
We denote by 
$(z_{\alpha}, \zeta_{\alpha})$ 
coordinates of 
$U_{\alpha} \times \C$. 
Now 
define holomorphic functions 
$\pi_{i, \alpha} : U_{\alpha} \times \C \to \C$ 
by
$$
\pi_{i, \alpha}(z_{\alpha}, \zeta_{\alpha}) 
:= \sigma_{i, \alpha}(z_{\alpha})\zeta_{\alpha}^{m_i}, 
$$
where 
$\sigma_{i,\alpha}$ 
is the local expression of 
$\sigma_i$ 
on 
$U_{\alpha}$.
Then 
the set 
$\left\{\pi_{i, \alpha}\right\}_{\alpha}$ 
of holomorphic functions defines a global holomorphic function 
$\pi_i : N_i \to \C$.

\begin{df} 
A degeneration 
$\pi : M \to \Delta$ 
is said to be \emph{linear} 
if for any irreducible component 
$\Theta_i$ 
of its singular fiber 
$X_0$, 
\begin{description}
\item[\rm (i)] 
a tubular neighborhood 
$N(\Theta_i)$ 
of 
$\Theta_i$ 
in 
$M$ 
is biholomorphic to a tubular neighborhood of a zero-section of the normal bundle 
$N_i$, 
and
\item[\rm (ii)] 
under the identification by the biholomorphic map of {\rm (i)}, 
the following conditions are satisfied:
\begin{itemize}
\item 
The restriction 
$\pi \big|_{N(\Theta_i)}$ 
coincides with the holomorphic function 
$\pi_i$ 
defined above. 
\item 
If 
$\Theta_i$ 
intersects 
$\Theta_j$ 
at a point 
$p$, 
then there exist local trivializations 
$U_{\alpha} \times \C$ 
of 
$N_i$ 
and 
$U_{\beta} \times \C$ 
of 
$N_j$ 
around 
$p$ 
such that neighborhoods of 
$p$ 
in 
$N(\Theta_i)$ 
and 
$N(\Theta_j)$ 
are identified by \emph{plumbing} 
$(z_{\alpha}, \zeta_{\alpha}) 
= (\zeta_{\beta}, z_{\beta})$ 
and 
$\pi$ 
is locally expressed as 
$$
\pi \big|_{N(\Theta_i)}(z_{\alpha}, \zeta_{\alpha}) 
= z_{\alpha}^{m_j} \zeta_{\alpha}^{m_i}, 
\quad 
\pi \big|_{N(\Theta_j)}(z_{\beta}, \zeta_{\beta}) 
= z_{\beta}^{m_i}\zeta_{\beta}^{m_j},  
$$
where 
$(z_{\alpha}, \zeta_{\alpha}) \in U_{\alpha} \times \C$ 
and 
$(z_{\beta}, \zeta_{\beta}) \in U_{\beta} \times \C$. 
\end{itemize}
\end{description}
\end{df}

\begin{rem}  \label{ta:rem:TopLinearizability}
Any degeneration of complex curves 
(even if the underlying reduced curve of its singular fiber 
does \emph{not} have at most simple normal crossings), 
after successive blowing up and down,
becomes a degeneration topologically equivalent to some linear degeneration.
\end{rem}

If 
$\pi : M \to \Delta$ 
is linear,
then we may express 
$M$ 
locally as a hypersurface in some space as follows: 
We first identify 
$M$ 
with the graph of 
$\pi$ 
in 
$M \times \Delta$ 
$$
\Graph(\pi) 
= \left\{(x, s) \in M \times \Delta \, : \, \pi(x) - s = 0\right\}
$$
via the natural projection 
$\Graph(\pi) \ni (x, s) \mapsto x \in M$. 
Recall that
for any irreducible component 
$\Theta_i$ 
of the singular fiber 
$X_0$, 
the map 
$\pi$ 
is expressed around 
$\Theta_i$ 
as
$$
\pi(z_i, \zeta_i) 
= \sigma_i(z_i) \zeta_i^{m_i}, 
$$
where 
$\sigma_i$ 
is the standard section of 
$N_i^{\otimes (-m_i)}$ 
on 
$\Theta_i$. 
Then 
we obtain the local expression of 
$M$ 
around 
$\Theta_i$: 
$$
\sigma_i(z_i) \zeta_i^{m_i} - s 
= 0 
\quad 
\textrm{in} \ N_i \times \Delta. 
$$
Note that 
these hypersurfaces are glued around the intersection points 
by plumbings 
$(z_j, \zeta_j, s) 
= (\zeta_i, z_i, s)$ 
where 
$(z_i, \zeta_i, s) \in N_i \times \Delta$ 
and 
$(z_j, \zeta_j, s) \in N_j \times \Delta$.

\begin{figure}[btp]
\begin{center}
\includegraphics[width=8cm,clip]{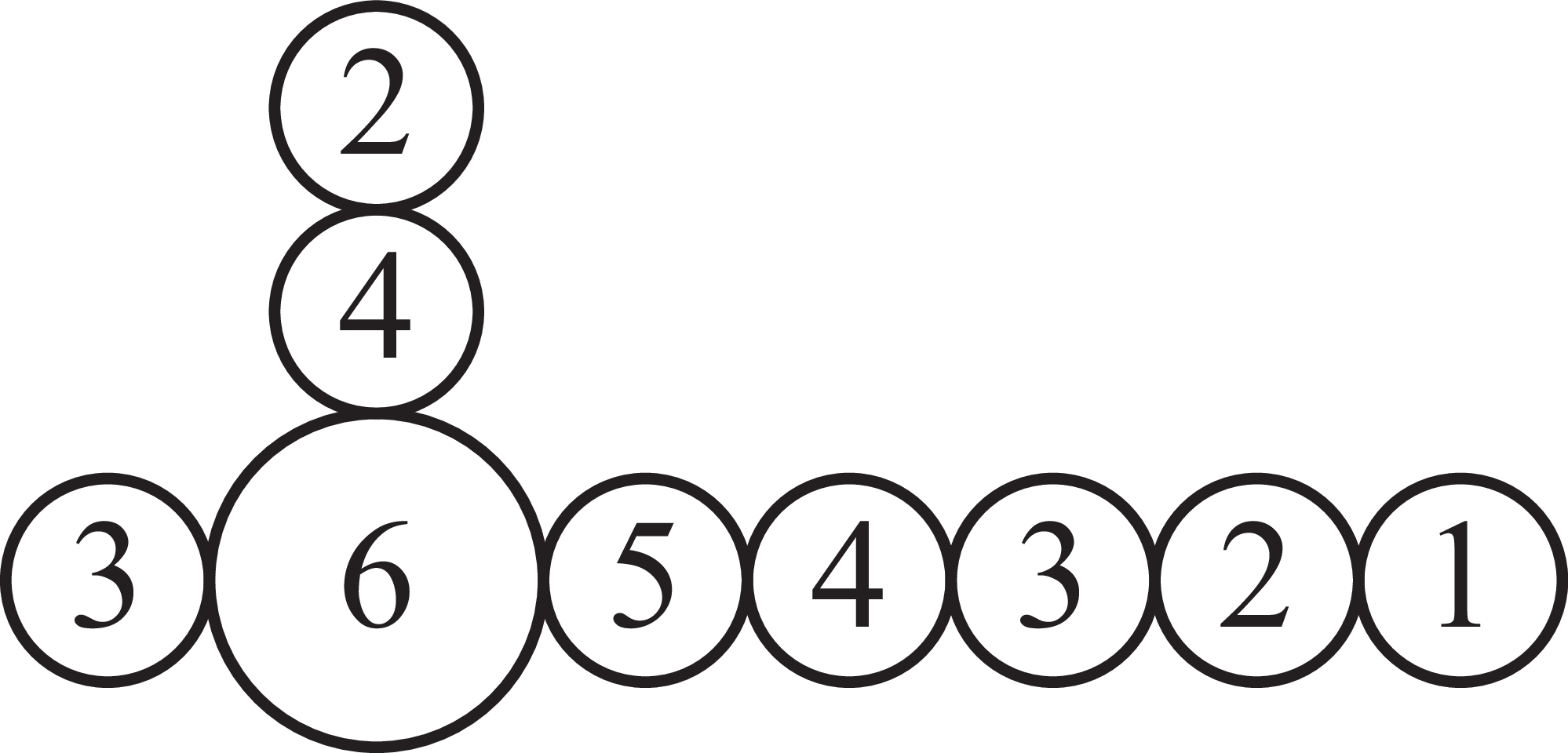}
\end{center}
\caption[]{\parbox[t]{7cm}{%
A singular fiber of type 
$II^*$ 
of a degeneration of elliptic curves is stellar. 
Each circle denotes a complex projective line, 
the number stands for its multiplicity, 
and each intersection point is a node. 
}}
\label{ta:fig:StellarExample}
\end{figure}

For a linear degeneration 
$\pi : M \to \Delta$, 
its singular fiber 
$X_0$ 
consists of three kinds of parts: 
\emph{cores}, \emph{branches} and \emph{trunks}. 
An irreducible component 
$\Theta_i$ 
of 
$X_0$ 
is called a \emph{core} 
if 
$\Theta_i$ 
intersects other irreducible components at at least three points 
or 
the genus of
$\Theta_i$ 
is positive. 
A \emph{branch} is a chain 
$\sum_i m_i \Theta_i$ 
of complex projective lines attached with a core on one hand, 
while 
a \emph{trunk} is a chain 
$\sum_i m_i \Theta_i$ 
of complex projective lines connecting other irreducible components on both hands. 
We say that 
$X_0$ is a \emph{stellar} singular fiber 
if 
$X_0$ 
consists of one core and branches emanating from the core.
See Figure \ref{ta:fig:StellarExample}. 
Otherwise 
$X_0$ 
is said to be \emph{constellar}. 
If 
$X_0$ 
is \emph{normally minimal}, that is, 
(i) 
any singularity of 
$X_0^{\reduced}$ 
is a node and 
(ii) 
any irreducible component that is a $(-1)$-curve 
(an exceptional curve of the first kind) 
intersects other irreducible component at at least three points, 
then 
all the branches and trunks of 
$X_0$ 
contain no $(-1)$-curves.

A degeneration 
whose singular fiber is a (\emph{fringed}) branch 
can be constructed explicitly and 
associated to a sequence of nonnegative integers 
(the \emph{multiplicity sequence}):

\begin{lem}  \label{ta:lem:BrExistence}
Let 
$m_0, m_1, \dots, m_{\lambda+1} \, (\lambda \geq 1)$ 
be nonnegative integers%
\footnote{
In this paper, 
by convention, 
we append 
$m_{\lambda+1} = 0$ 
to 
the sequence 
$m_0, m_1, \dots, m_\lambda$ 
of positive integers, 
so that 
$r_\lambda := (m_{\lambda-1} + m_{\lambda+1}) / m_\lambda$ 
equals 
$m_{\lambda-1} / m_\lambda$. 
See \cite{Ta3} Section 5.1. 
} 
satisfying the following conditions: 
$$
\left\{
\begin{array}{l}
m_0 > m_1 > \cdots > m_\lambda > m_{\lambda+1} = 0 \ \textrm{and} \\
r_i := \frac{m_{i-1} + m_{i+1}}{m_i} \, (i = 1, 2, \dots, \lambda) 
\textrm{ is an integer greater than } 1. 
\end{array}
\right.
$$
Then there exists a degeneration 
$\pi : M \to \Delta$
with the singular fiber
$$
X_0 
= m_0 \Delta_0 + m_1 \Theta_1 + m_2 \Theta_2 + \cdots + m_\lambda \Theta_\lambda, 
$$
where 
$\Delta_0 = \C$, 
and 
$\Theta_1, \Theta_2, \dots, \Theta_\lambda$ 
are complex projective lines, 
and each pair of 
$\Theta_i$ 
and 
$\Theta_{i+1} \, (i = 1, 2, \dots, \lambda-1)$ 
and 
$\Delta_0$ 
and 
$\Theta_1$ 
intersect transversely at one point. 
\end{lem}

\begin{proof}
We take 
$\lambda$ 
copies 
$\Theta_1, \Theta_2, \dots, \Theta_\lambda$ 
of the complex projective line. 
Let 
$\Theta_i = U_i \cup V_i$ 
be an open covering by two complex lines 
$U_i, V_i \, (= \C)$
with coordinates 
$w_i \in U_i \setminus \left\{ 0 \right\}$ 
and 
$z_i \in V_i \setminus \left\{ 0 \right\}$ 
satisfying 
$z_i = 1 / w_i$. 
Then we obtain a line bundle 
$N_i$ 
on 
$\Theta_i$ 
of degree
$-r_i$ 
from 
$U_i \times \C$ 
and 
$V_i \times \C$ 
by identifying 
$(z_i, \zeta_i) \in (V_i \setminus \left\{ 0 \right\}) \times \C$ 
with 
$(w_i, \eta_i) \in (U_i \setminus \left\{ 0 \right\}) \times \C$ 
via
$$
g_i \, : \, z_i = \frac{1}{w_i}, 
\quad \zeta_i = w_i^{r_i} \eta_i. 
$$
Now consider the hypersurface 
$W_i$ 
in 
$N_i \times \Delta$ 
defined by
$$
\begin{cases}
H_i : 
w_i^{m_{i - 1}} \eta_i^{m_i} - s = 0, 
&\quad 
\textrm{in} \ U_i \times \C \times \Delta, \\
H'_i : 
z_i^{m_{i + 1}} \zeta_i^{m_i} - s = 0, 
&\quad 
\textrm{in} \ V_i \times \C \times \Delta. 
\end{cases}
$$
Under plumbings 
$(w_{i+1}, \eta_{i+1}, s) 
= (\zeta_i, z_i, s)$ 
of 
$N_i \times \Delta$ 
and 
$N_{i+1} \times \Delta \, (i = 1, 2, \dots, \lambda-1)$, 
the hypersurfaces 
$W_1, W_2, \dots, W_\lambda$ 
are glued, 
so that 
they together define a smooth complex surface 
$M$. 
Letting 
$\pi : M \to \Delta$ 
be the natural projection, 
the central fiber is 
$$
\pi^{-1}(0) 
= m_0 \Delta_0 + m_1 \Theta_1 + m_2 \Theta_2 + \cdots + m_\lambda \Theta_\lambda, 
$$
where 
$\Delta_0 
:= \left\{ 0 \right\} \times \C 
\subset U_1 \times \C$. 
Thus the holomorphic map 
$\pi : M \to \Delta$ 
is the desired degeneration. 
\end{proof}

\begin{rem}
Precisely speaking, 
the holomorphic function 
$\pi : M \to \Delta$ 
obtained in Lemma \ref{ta:lem:BrExistence}
does not satisfy the condition to be a degeneration. 
Indeed 
$\pi$ 
is not proper. 
Note that 
we consider the restriction of a degeneration 
to a tubular neighborhood of a branch. 
\end{rem}

Next 
we define a special subdivisor of a stellar singular fiber. 
Let 
$\pi : M \to \Delta$ 
be a linear degeneration of complex curves 
with the stellar singular fiber 
$X_0 = m_0 \Theta_0 + \sum_{j=1}^h \Br^{(j)}$, 
where 
$\Theta_0$ 
is the core and 
$\Br^{(j)} \, (j = 1, 2, \dots, h)$ 
is a branch. 
Write 
$\Br^{(j)} 
= m_1^{(j)} \Theta_1^{(j)} + m_2^{(j)} \Theta_2^{(j)} + 
\cdots + m_{\lambda^{(j)}}^{(j)} \Theta_{\lambda^{(j)}}^{(j)}$ 
and let 
$\overline{\Br}^{(j)} 
= m_0 \Delta_0^{(j)} + m_1^{(j)} \Theta_1^{(j)} + 
\cdots + m_{\lambda^{(j)}}^{(j)} \Theta_{\lambda^{(j)}}^{(j)}$ 
be a fringed branch. 
Consider a connected subdivisor 
$Y = n_0 \Theta_0 + \sum_{j=1}^h \br^{(j)}$ 
of 
$X_0$, 
where 
$\br^{(j)} 
:= n_1^{(j)} \Theta_1^{(j)} + n_2^{(j)} \Theta_2^{(j)} + 
\cdots + n_{\nu^{(j)}}^{(j)} \Theta_{\nu^{(j)}}^{(j)} 
\, (j = 1, 2, \dots, h)$. 
Here 
$Y$ 
satisfies 
$0 \leq \nu^{(j)} \leq \lambda^{(j)}$ 
and 
$0 < n_i^{(j)} \leq m_i^{(j)}$ 
for each 
$i$ 
and 
$j$. 
Set 
$\overline{\br}^{(j)} 
:= n_0 \Delta_0^{(j)} + n_1^{(j)} \Theta_1^{(j)} + n_2^{(j)} \Theta_2^{(j)} 
+ \cdots + n_{\nu^{(j)}}^{(j)} \Theta_{\nu^{(j)}}^{(j)}$. 
For the time being,
we consider 
$\overline{\Br}^{(j)}$ 
and 
$\overline{\br}^{(j)}$, 
omitting the superscript 
$(j)$ 
to simplify notation. 
We call 
$\overline{\br}$ 
a \emph{subbranch} of 
$\overline{\Br}$ 
if one of the following conditions is satisfied:
\begin{itemize}
\item 
$\nu = 0, 1$, or 
\item 
$\nu \geq 2$ and $n_{i+1} = r_i n_i - n_{i-1} 
\, (i = 1, 2, \dots, \nu-1),$
\end{itemize}
where 
$r_i := (m_{i-1} + m_{i+1}) / m_i$
(see Lemma \ref{ta:lem:BrExistence}). 
Set 
$n_{\nu + 1} 
:= r_\nu n_\nu - n_{\nu - 1}$. 
If 
$\nu = 0$, 
then we set 
$n_{\nu + 1} = n_1 := 0$. 
Define the three types of subbranches for a positive integer 
$l$ 
as follows: 
\begin{description}
\item[Type $A_l$] 
A subbranch 
$\overline{\br}$ 
of 
$\overline{\Br}$ 
is of \emph{type $A_l$} 
if 
$l n_i \leq m_i$ 
for each 
$i$ 
and 
$n_{\nu+1} \leq 0$. 
\item[Type $B_l$] 
A subbranch 
$\overline{\br}$ 
of 
$\overline{\Br}$ 
is of \emph{type $B_l$} 
if 
$l n_i \leq m_i$ 
for each 
$i$, 
$n_\nu = 1$ 
and 
$m_\nu = l$. 
\item[Type $C_l$] 
A subbranch 
$\overline{\br}$ 
of 
$\overline{\Br}$ 
is of \emph{type $C_l$}
if 
$l n_i \leq m_i$ 
for each 
$i$, 
$n_\nu = n_{\nu+1}$ 
and 
$m_\nu - m_{\nu+1}$ 
divides 
$l$.
\end{description}

Now we return to a connected subdivisor 
$Y$ 
of the stellar singular fiber 
$X_0$.

\begin{df} 
Let 
$Y = n_0 \Theta_0 + \sum_{j=1}^h \br^{(j)}$ 
be a connected subdivisor of 
$X_0$ 
such that 
$n_0 < m_0$ 
and each 
$\overline{\br}^{(j)}$ 
is a subbranch of 
$\overline{\Br}^{(j)}$. 
$Y$ 
is called a \emph{crust} of 
$X_0$ 
if there exists a meromorphic section 
$\tau$ 
of the line bundle 
$N_0^{\otimes n_0}$ 
on 
$\Theta_0$ 
such that for some nonnegative divisor 
$D = \sum_{i=1}^k a_i q_i$ on $\Theta_0$, 
$$
\divisor(\tau) 
= - \sum_{j=1}^h n_1^{(j)} p^{(j)} + D,
$$
where 
$N_0$ 
denotes the normal bundle of 
$\Theta_0$ 
in 
$M$, 
$\left\{p^{(j)}\right\}$ 
is the set of the attachment points on 
$\Theta_0$ 
with the branches 
$\Br^{(j)}$. 
Moreover, 
for a positive integer 
$l$, 
if each 
$\overline{\br}^{(j)}$ 
is a subbranch of 
$\overline{\Br}^{(j)}$ 
of either type $A_l$, type $B_l$ or type $C_l$, 
then we call 
$Y$ 
a \emph{simple crust} of 
$X_0$ 
with \emph{barking multiplicity} 
$l$. 
\end{df}

We call the meromorphic section 
$\tau$ 
a \emph{core section}. 
Note that 
$\tau$ 
is not uniquely determined by 
$Y$. 
Setting 
$r_0 := \sum_{j=1}^h m_1^{(j)} / m_0$ 
and 
$r_0' := \sum_{j=1}^h n_1^{(j)} / n_0$, 
the following holds:

\begin{lem}  \label{ta:lem:CoreSecExistence}
Suppose that 
$\Theta_0$ 
is a complex projective line. 
Then 
a connected subdivisor 
$Y$ 
is a crust of 
$X_0$ 
$($equivalently, 
$Y$ 
has a core section 
$\tau)$ 
if and only if 
$r_0 \leq r_0'$. 
Moreover
$\tau$ 
has no zero, 
that is, 
$D = 0$ 
exactly when 
$r_0 = r_0'$. 
\end{lem}

\begin{figure}[btp]
\begin{center}
\includegraphics[width=10cm,clip]{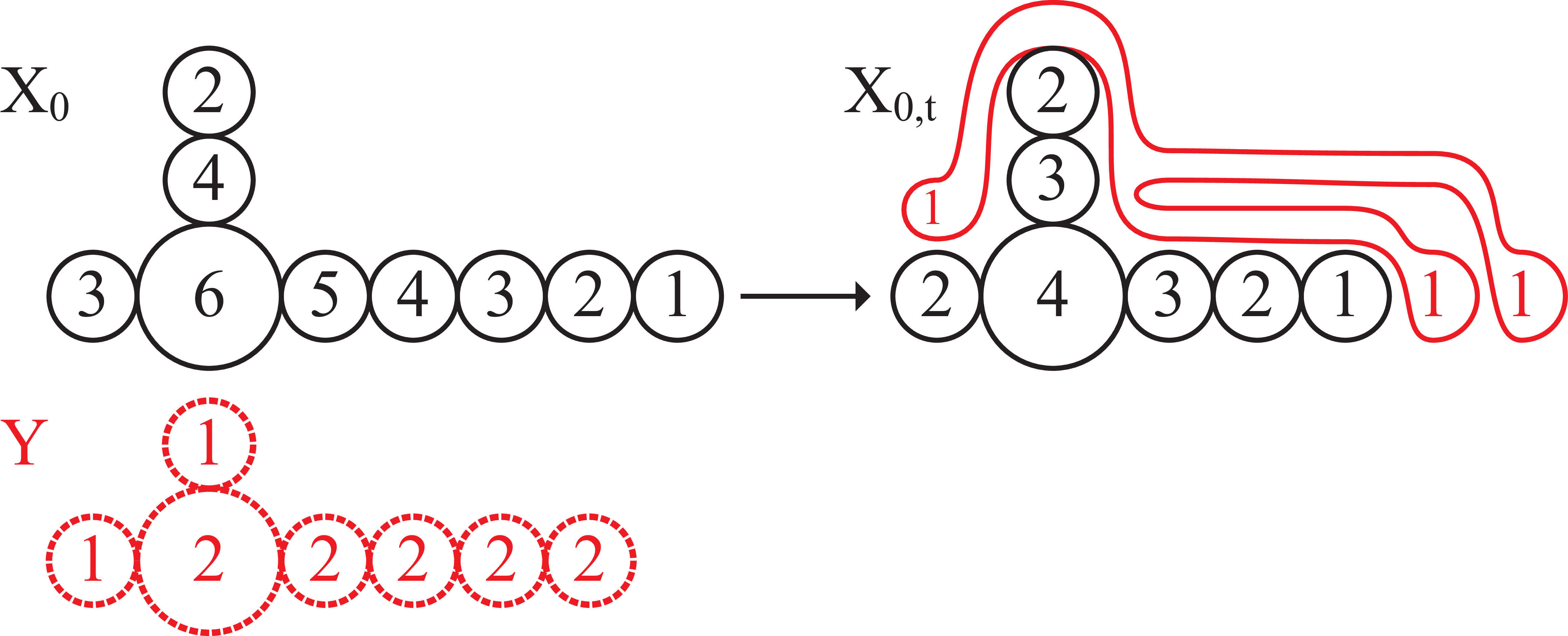}
\end{center}
\caption[]{\parbox[t]{7cm}{%
In the barking family $\boldsymbol{[II^*.1]}$, 
the singular fiber of type $II^*$ 
is deformed to the main fiber of type $III^*$. 
It seems that the simple crust $Y$ is ``barked'' (peeled) off from the original singular fiber. 
}}
\label{ta:fig:BarkingExample}
\end{figure}

Takamura constructed a deformation family of 
$\pi : M \to \Delta$ 
associated with a simple crust 
$Y$. 
We call a deformation family obtained by his method 
a \emph{barking family}. 
In a barking family, 
the original singular fiber 
$X_0$ 
is deformed to a simpler singular fiber 
in such a way that 
a part of 
$X_0$ 
looks ``barked'' off from 
$X_0$. 
The resulting singular fiber appears over the origin of 
$\Delta_t$, 
so we denote it by 
$X_{0, t}$ 
and call it 
the \emph{main fiber}. 
See Figure \ref{ta:fig:BarkingExample}.

In a barking family, 
there appear not only the main fiber 
but also other singular fibers over some points 
away from the origin of 
$\Delta_t$, 
which are called \emph{subordinate fibers}.
It is easy to see this. 
Under the deformation, 
the topological type of the singular fiber over the origin changes, 
so the local monodromy around it also changes 
(see Section \ref{sec:Monodromies} for details). 
On the other hand, 
the global monodromies before and after the deformation 
--- 
that is, 
the two monodromies 
each of which is induced by a loop in 
$\Delta$ 
(resp. 
$\Delta_t$) 
parallel to its boundary 
$\partial \Delta$ 
(resp. 
$\partial \Delta_t$) 
--- 
coincide. 
We then deduce that 
there appear other singular fibers with nontrivial monodromies. 
Thus every barking family turns out to be a splitting family. 
Therefore:

\begin{thm}[Takamura \cite{Ta3}]  \label{ta:thm:Ta3}
Let 
$\pi : M \to \Delta$ 
be a linear degeneration with the stellar singular fiber 
$X_0$. 
If 
$X_0$ 
has a simple crust 
$Y$, 
then 
$\pi : M \to \Delta$ 
admits a splitting family 
$\Psi : \mathcal{M} \to \Delta \times \Delta^{\dagger}$. 
\end{thm}

\begin{rem}
In this paper, 
for a degeneration 
which is \emph{not} necessarily relatively minimal, 
a splitting family of it is defined 
to satisfy that 
each deformation has at least two singular fibers 
(see Section \ref{sec:Intro}). 
Thus 
singular fibers of a deformation in a splitting family 
possibly become smooth fibers after blowing down. 
Such singular fibers are said to be \emph{fake}. 
\end{rem}

\begin{table}[p]\footnotesize
\begin{center}
\caption{Kodaira's notation.}
\label{ta:table:KodairaNotation} 
\begin{tabular}{|c|c|c|c|c|} \hline
&a singular fiber 
$X$ 
&$e(X)$ 
&$A_X$ 
&$\Tr(A_X)$ \\ \hline
\rule[-12pt]{0pt}{30pt}
\begin{tabular}[c]{c}
$m I_0$ \\
{\footnotesize $(m\geq 2)$}
\end{tabular} 
&\begin{tabular}[c]{c}
a multiple torus
\end{tabular} 
&$0$ 
&$\begin{pmatrix}1&0\\0&1\end{pmatrix}$ 
&$2$ \\ \hline
\rule[-12pt]{0pt}{30pt}
$m I_1$ 
&\begin{tabular}[c]{c}
a (multiple) projective line \\
with one node
\end{tabular} 
&$1$ 
&$\begin{pmatrix}1&1\\0&1\end{pmatrix}$ 
&$2$ \\ \hline
\rule[-5pt]{0pt}{32pt}
$m I_n$ 
&\begin{tabular}[c]{c}
\includegraphics[height=40pt,clip]{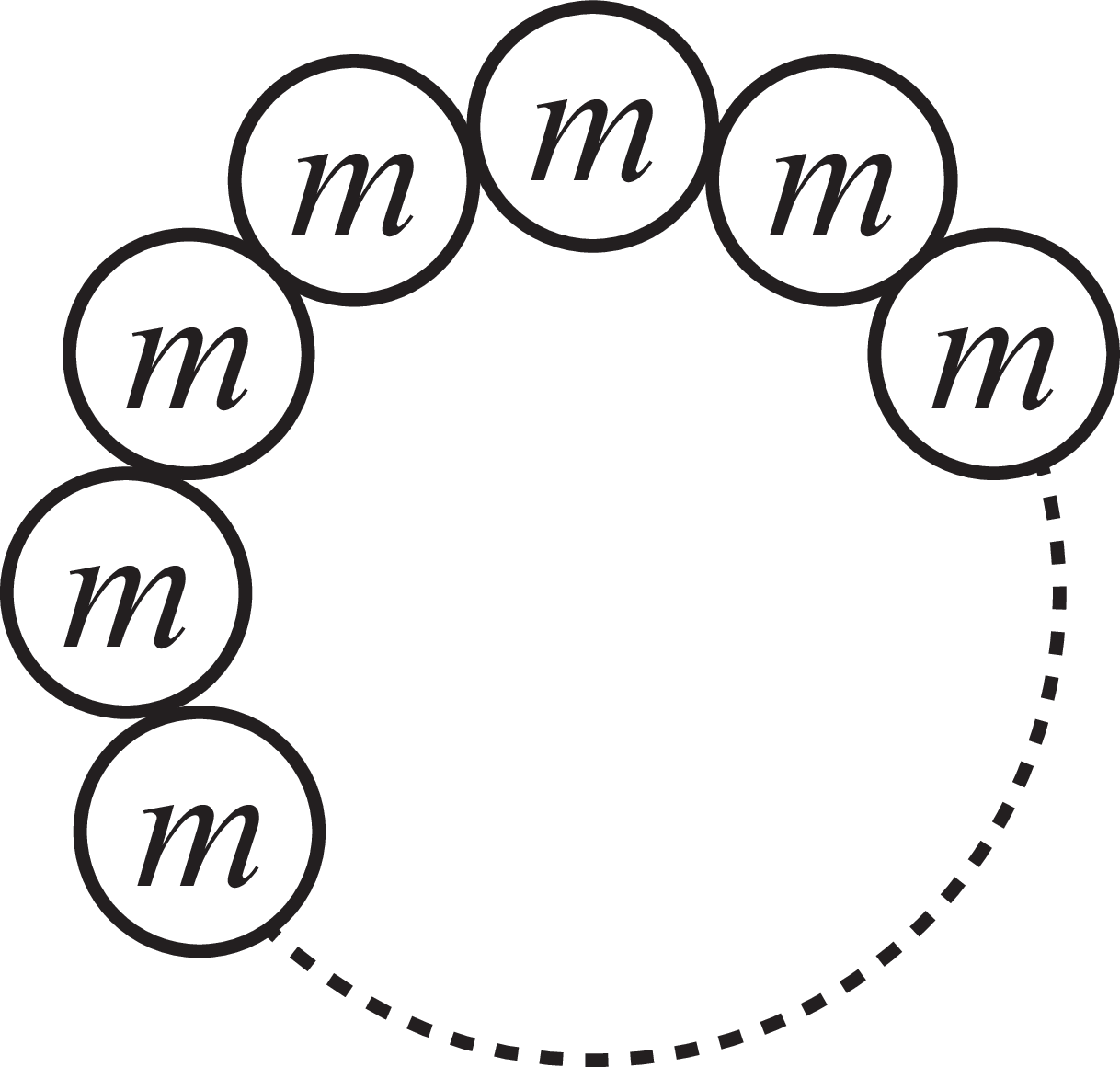}
\end{tabular} 
&$n$ 
&$\begin{pmatrix}1&n\\0&1\end{pmatrix}$ 
&$2$ \\ \hline
\rule[-12pt]{0pt}{30pt}
$II$ 
&\begin{tabular}[c]{c}
a projective line \\
with one cusp
\end{tabular}
&$2$ 
&$\begin{pmatrix}0&1\\-1&1\end{pmatrix}$ 
&$1$ \\ \hline
\rule[-12pt]{0pt}{30pt}
$III$ 
&\begin{tabular}[c]{c}
two projective lines \\
with second order contact
\end{tabular}
&$3$ 
&$\begin{pmatrix}0&1\\-1&0\end{pmatrix}$ 
&$0$\\ \hline
\rule[-12pt]{0pt}{30pt}
$IV$ 
&\begin{tabular}[c]{c}
three projective lines intersecting \\
transversally at one point
\end{tabular}
&$4$ 
&$\begin{pmatrix}-1&1\\-1&0\end{pmatrix}$ 
&$-1$ \\ \hline
\rule[-10pt]{0pt}{35pt}
$I_0^*$ 
&\begin{tabular}[c]{c}
\includegraphics[height=35pt,clip]{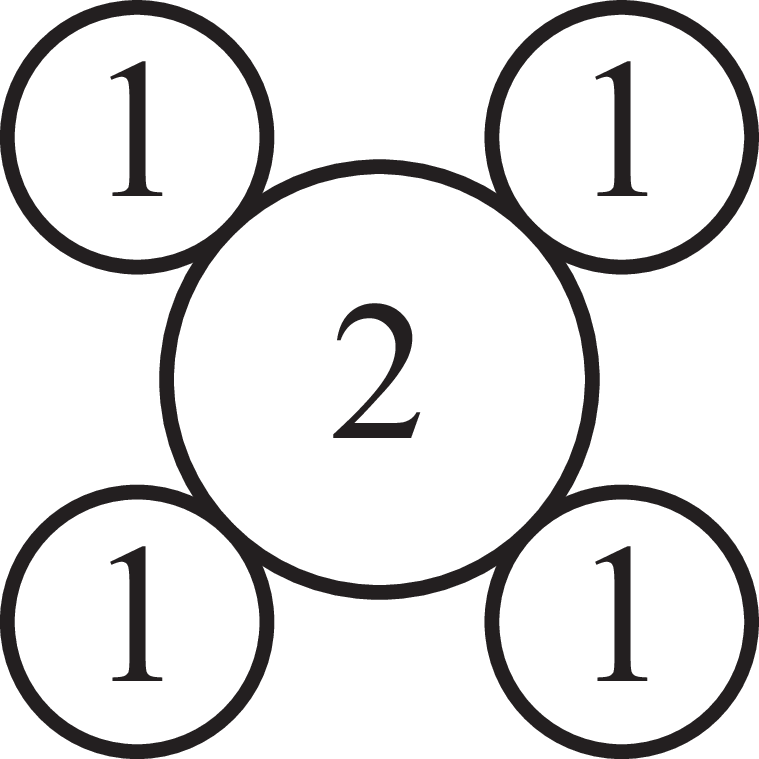}
\end{tabular}
&$6$ 
&$\begin{pmatrix}-1&0\\0&-1\end{pmatrix}$ 
&$-2$ \\ \hline
\rule[-10pt]{0pt}{35pt}
$I_n^*$ 
&\begin{tabular}[c]{c}
\includegraphics[height=35pt,clip]{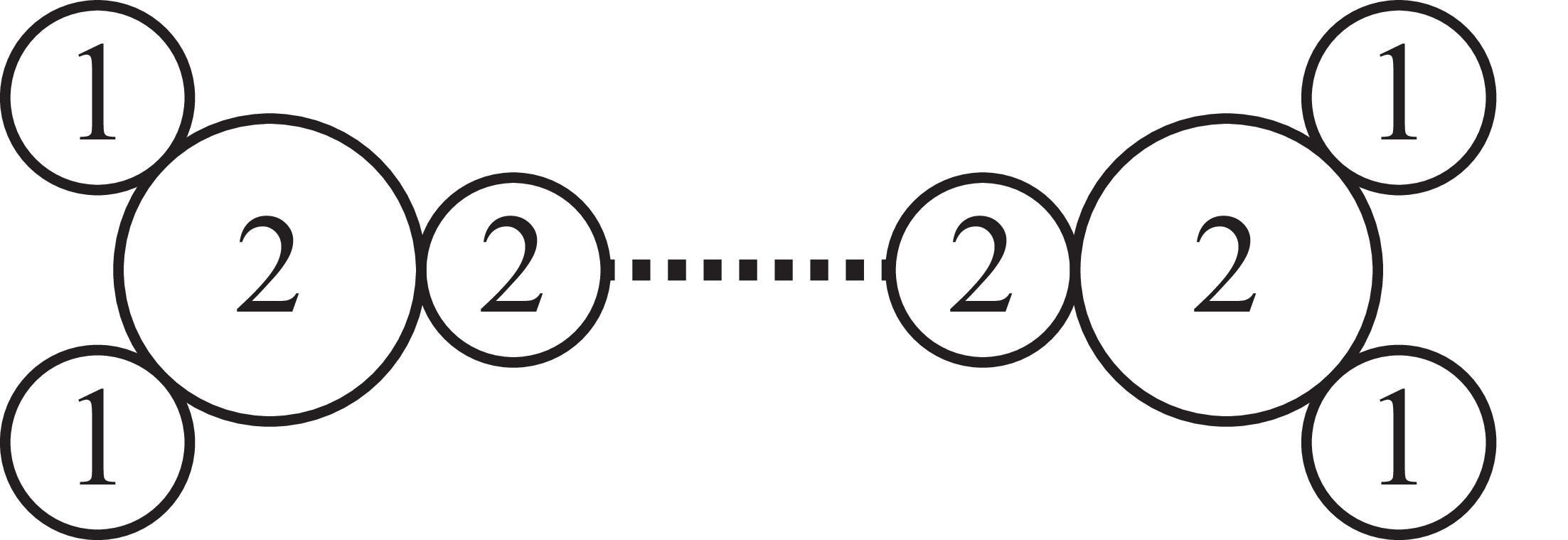}
\end{tabular}
&$6+n$ 
&$\begin{pmatrix}-1&-n\\0&-1\end{pmatrix}$ 
&$-2$ \\ \hline
\rule[-10pt]{0pt}{40pt}
$II^*$ 
&\begin{tabular}[c]{c}
\includegraphics[height=45pt,clip]{II+.eps}
\end{tabular}
&$10$ 
&$\begin{pmatrix}1&-1\\1&0\end{pmatrix}$ 
&$1$ \\ \hline
\rule[-10pt]{0pt}{35pt}
$III^*$ 
&\begin{tabular}[c]{c}
\includegraphics[height=35pt,clip]{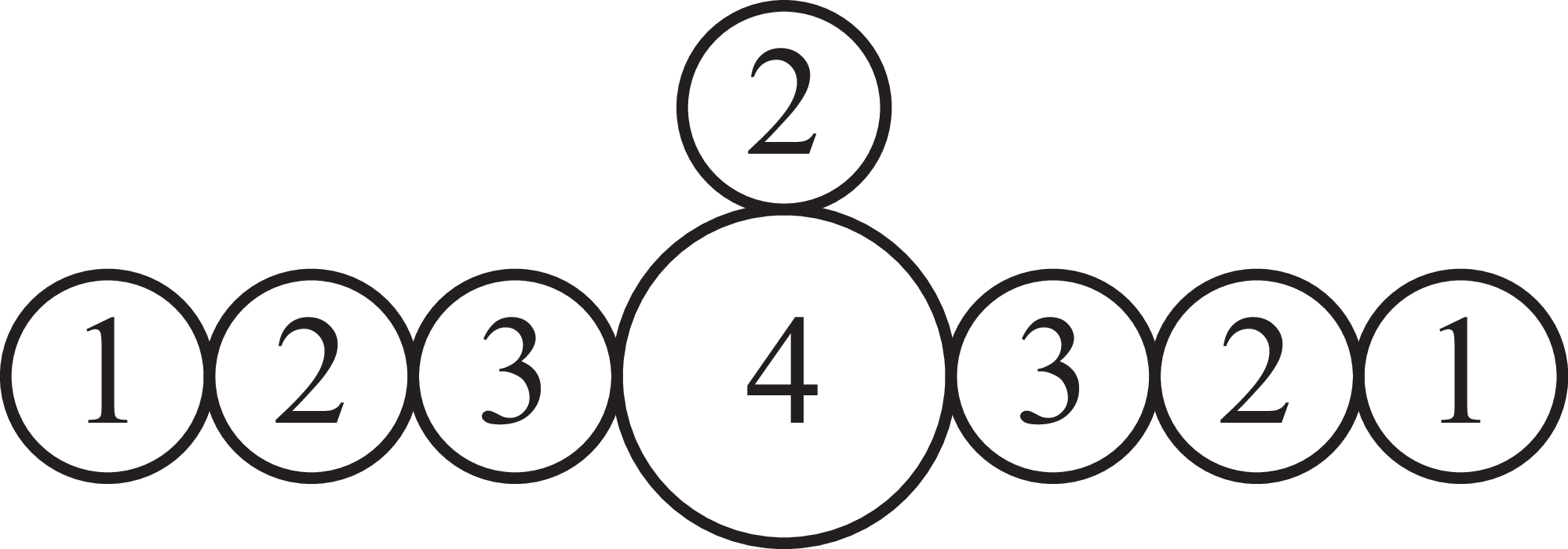}
\end{tabular} 
&$9$ 
&$\begin{pmatrix}0&-1\\1&0\end{pmatrix}$ 
&$0$ \\ \hline
\rule[-10pt]{0pt}{40pt}
$IV^*$ 
&\begin{tabular}[c]{c}
\includegraphics[height=45pt,clip]{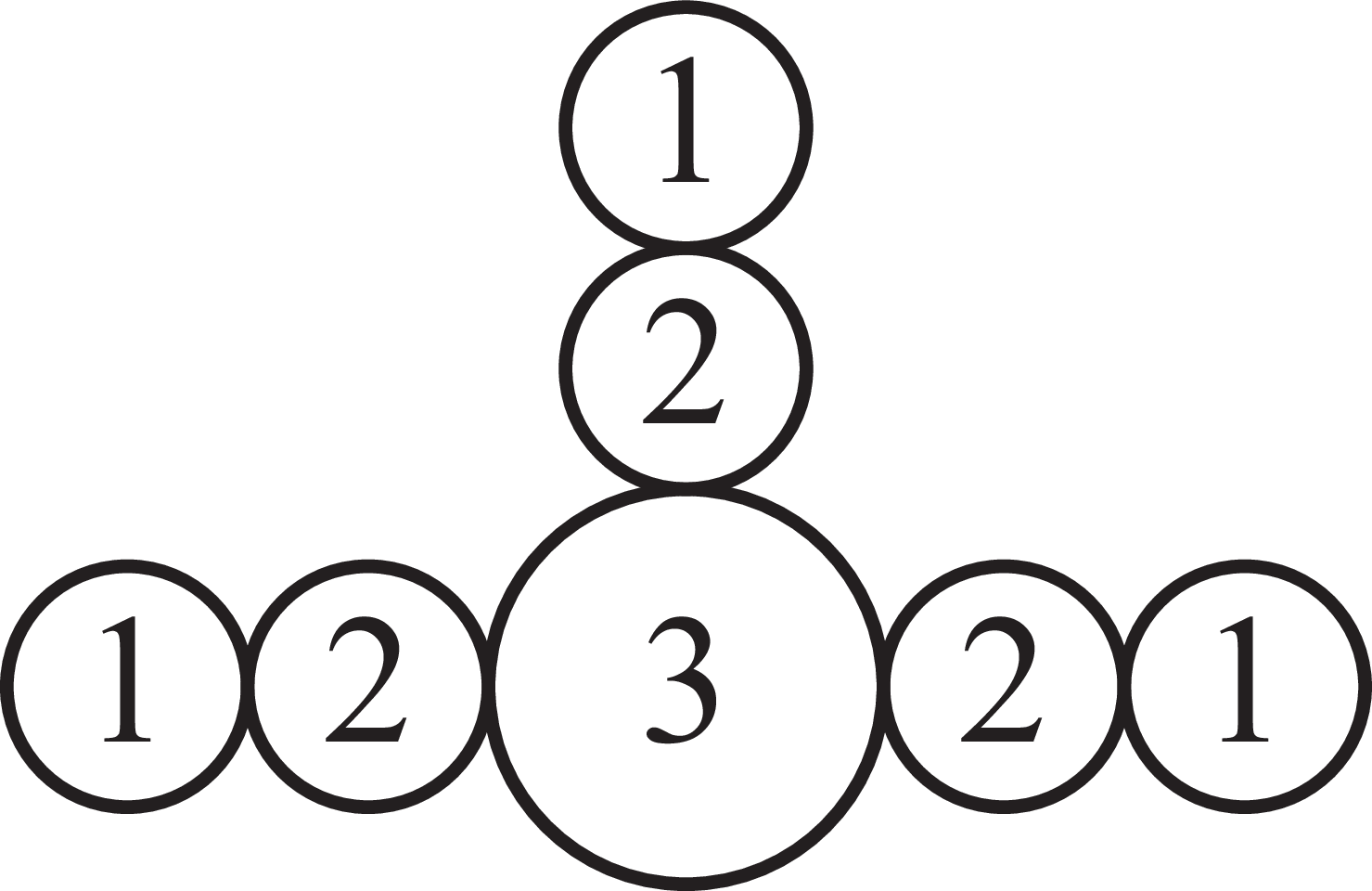}
\end{tabular}
&$8$ 
&$\begin{pmatrix}0&-1\\1&-1\end{pmatrix}$ 
&$-1$ \\ \hline
\end{tabular}
\end{center}
\end{table}%

\subsection*{Kodaira's notation}

Before proceeding,
we supply 
Kodaira's list of singular fibers 
of (relatively) minimal degenerations of elliptic curves \cite{Ko}. 
See Table \ref{ta:table:KodairaNotation}. 
For a singular fiber $X$, 
we denote by $e(X)$ 
the topological Euler characteristic of the underlying reduced curve $X^{\reduced}$ of $X$.
$A_X \in SL(2, \Z)$ 
is 
the standard monodromy matrix of $X$ 
and its trace is denoted by $\Tr(A_X)$. 

Note that 
minimal singular fibers of type $I_n^*$, $II^*$, $III^*$ and $IV^*$ in this table 
are normally minimal 
and their underlying reduced curves have 
at most simple normal crossings. 
In contrast, 
minimal singular fibers of type $II$, $III$ and $IV$ 
have a singularity that is not a node. 
However, 
after successive blowing up, 
they become 
normally minimal degenerations 
such that 
$X^{\reduced}$ 
has at most simple normal crossings. 
In this paper, 
such degenerations are also referred 
to be of type $II$, $III$ and $IV$.

\section{Constraints from Euler characteristics}  \label{sec:Euler}

In \cite{Ta3}, 
Takamura listed 
thirty five barking families 
for degenerations of complex curves of genus 
$g = 1$, 
that is, 
for degenerations of elliptic curves, 
and determined the type of the main fiber of each of them 
as follows 
(see also the list in Section \ref{sec:Talist}): 
\begin{align*}
&\boldsymbol{[II.1]} \quad II \barkarrow I_1 &\hspace{30pt}
 &\boldsymbol{[III^*.5]} \quad III^* \barkarrow I_6 \nonumber\\
&\boldsymbol{[II.2]} \quad II \barkarrow I_1 &\hspace{30pt}
 &\boldsymbol{[III^*.6]} \quad III^* \barkarrow I_2^* \nonumber\\
&\boldsymbol{[II^*.1]} \quad II^* \barkarrow III^* &\hspace{30pt}
 &\boldsymbol{[III^*.7]} \quad III^* \barkarrow I_7 \nonumber\\
&\boldsymbol{[II^*.2]} \quad II^* \barkarrow IV^* &\hspace{30pt}
 &\boldsymbol{[III^*.8]} \quad III^* \barkarrow I_6 \nonumber\\
&\boldsymbol{[II^*.3]} \quad II^* \barkarrow I_2^* &\hspace{30pt}
 &\boldsymbol{[III^*.9]} \quad III^* \barkarrow IV^* \nonumber\\
&\boldsymbol{[II^*.4]} \quad II^* \barkarrow I_5 &\hspace{30pt}
 &\boldsymbol{[IV.1]} \quad IV \barkarrow I_3 \nonumber\\
&\boldsymbol{[II^*.5]} \quad II^* \barkarrow I_3^* &\hspace{30pt}
 &\boldsymbol{[IV.2]} \quad IV \barkarrow I_2 \nonumber\\
&\boldsymbol{[II^*.6]} \quad II^* \barkarrow I_3^* &\hspace{30pt}
 &\boldsymbol{[IV.3]} \quad IV \barkarrow I_2 \nonumber\\
&\boldsymbol{[II^*.7]} \quad II^* \barkarrow I_8 &\hspace{30pt}
 &\boldsymbol{[IV.4]} \quad IV \barkarrow II \nonumber\\
&\boldsymbol{[II^*.8]} \quad II^* \barkarrow III^* &\hspace{30pt}
 &\boldsymbol{[IV^*.1]} \quad IV^* \barkarrow I_1^* \nonumber\\
&\boldsymbol{[II^*.9]} \quad II^* \barkarrow III^* &\hspace{30pt}
 &\boldsymbol{[IV^*.2]} \quad IV^* \barkarrow I_0^* \nonumber\\
&\boldsymbol{[III.1]} \quad III \barkarrow I_2 &\hspace{30pt}
 &\boldsymbol{[IV^*.3]} \quad IV^* \barkarrow I_6 \nonumber\\
&\boldsymbol{[III.2]} \quad III \barkarrow I_1 &\hspace{30pt}
 &\boldsymbol{[IV^*.4]} \quad IV^* \barkarrow I_1^* \nonumber\\
&\boldsymbol{[III.3]} \quad III \barkarrow I_2 &\hspace{30pt}
 &\boldsymbol{[I_0^*.1]} \quad I_0^* \barkarrow I_4 \nonumber\\
&\boldsymbol{[III^*.1]} \quad III^* \barkarrow IV^* &\hspace{30pt}
 &\boldsymbol{[I_0^*.2]} \quad I_0^* \barkarrow I_3 \nonumber\\
&\boldsymbol{[III^*.2]} \quad III^* \barkarrow I_1^* &\hspace{30pt}
 &\boldsymbol{[I_n^*.1]} \quad I_n^* \barkarrow I_{n-1}^* \nonumber\\
&\boldsymbol{[III^*.3]} \quad III^* \barkarrow I_2^* &\hspace{30pt}
 &\boldsymbol{[I_n^*.2]} \quad I_n^* \barkarrow I_{n+4}. \nonumber\\
&\boldsymbol{[III^*.4]} \quad III^* \barkarrow I_0^* &\hspace{30pt}
 \nonumber
\end{align*}
The aim of this paper is 
to determine the subordinate fibers of the above barking families. 
In this section, 
we give a list of the sets of subordinate fibers 
that can appear in each of the barking families, 
using results on Euler characteristics of singular fibers of degenerations.

For a singular fiber 
$X$, 
we denote by 
$e(X)$ 
the \emph{topological Euler characteristic} 
of the underlying reduced curve of 
$X$.

\begin{lem}  \label{det:lem:SplittingInv}
Let 
$\pi : M \to \Delta$ 
be a degeneration of complex curves of genus 
$g \geq 1$ 
with the singular fiber 
$X_0$ 
and let 
$\Psi : \mathcal{M} \to \Delta \times \Delta^{\dagger}$ 
be a splitting family of 
$\pi : M \to \Delta$,
say, 
$X_0$ 
splits into singular fibers 
$X_1, X_2, \dots, X_N \, 
(N \geq 2)$ 
of a deformation 
$\pi_t : M_t \to \Delta_t$. 
\begin{description}
\item[\rm(a)] 
Then the following formula holds: 
$$
e(X_0) - 2(1 - g) 
= \sum_{i=1}^N \left\{e(X_i) - 2(1 - g)\right\}.
$$
\item[\rm(b)] 
In particular, 
if 
$g = 1$, 
then the following holds: 
\begin{equation}  \label{eq:Eulercharsum}
e(X_0) = e(X_1) + e(X_2) + \dots + e(X_N). 
\end{equation}
\end{description}
\end{lem}

\begin{proof}
(a) 
The left hand side equals the Euler characteristic 
$e(M)$ 
of 
$M$, 
while 
the right hand side equals 
$e(M_t)$ 
(see \cite[p.~97]{BPV}). 
Since 
$M_t$ 
is diffeomorphic to 
$M$, 
we have 
$e(M) = e(M_t)$, 
which confirms the assertion.

(b) clearly follows from (a). 
\end{proof}

Consider a barking family 
$\Psi : \mathcal{M} \to \Delta \times \Delta^{\dagger}$ 
of the degeneration 
$\pi : M \to \Delta$ 
of elliptic curves. 
Recall that
for a singular fiber 
$X_{s, t} := \Psi^{-1}(s, t) \, (t \neq 0)$, 
we call 
$X_{s, t}$ 
the main fiber if 
$s = 0$, 
and a subordinate fiber if 
$s \neq 0$. 
Suppose that 
$\Psi : \mathcal{M} \to \Delta \times \Delta^{\dagger}$ 
splits the original singular fiber 
$X_0$ 
into the main fiber 
$X_{0, t}$ 
and subordinate fibers 
$X_{s_1, t}, X_{s_2, t}, \dots, X_{s_N, t} 
\, (N \geq 1)$. 
In these notations, 
we restate $(\ref{eq:Eulercharsum})$ in Lemma \ref{det:lem:SplittingInv} as 
\begin{equation}
e(X_0) = e(X_{0, t}) + \sum_{i=1}^N e(X_{s_i, t}). 
\end{equation}
This confirms (a) of the following:

\begin{lem}  \label{det:lem:SingleLefshetzFiber}
Let 
$\pi : M \to \Delta$ 
be a degeneration of elliptic curves with the singular fiber 
$X_0$. 
Suppose that 
a barking family 
$\Psi : \mathcal{M} \to \Delta \times \Delta^{\dagger}$ 
splits the original singular fiber 
$X_0$ 
into the main fiber 
$X_{0, t}$ 
and subordinate fibers 
$X_{s_1, t}, X_{s_2, t}, \dots, X_{s_N, t} 
\, (N \geq 1)$. 
Then: 
\begin{description}
 \item[\rm(a)] 
The sum of the Euler characteristics of the subordinate fibers is 
$e(X_0) - e(X_{0, t})$: 
$$
\sum_{i=1}^N e(X_{s_i, t}) = e(X_0) - e(X_{0, t}). 
$$
 \item[\rm(b)] 
If 
$e(X_0) - e(X_{0, t}) = 1$ 
holds, then 
$\Psi$ 
splits 
$X_0$ 
into the main fiber 
$X_{0, t}$ 
and one subordinate fiber
$I_1$: 
$$
X_0 \longrightarrow X_{0, t} + I_1. 
$$
\end{description}
\end{lem}

\begin{proof}
It remains 
to show the second statement (b). 
From the assumption 
$e(X_0) - e(X_{0, t}) = 1$ 
together with (a), 
we have 
$$
e(X_{s_1, t}) + e(X_{s_2, t}) + \cdots + e(X_{s_N, t}) 
= 1. 
$$
Note that 
every subordinate fiber of any barking family 
is a reduced curve only with $A$-singularities 
(Lemma \ref{det:lem:SubordSingularity}). 
In particular, 
each subordinate fiber 
$X_{s_i, t}$ 
is not a multiple torus 
(whose Euler characteristic is 
$0$), 
thus 
$e(X_{s_i, t}) \geq 1$. 
Hence we have 
$N = 1$ 
(that is, 
$X_{s_1, t}$ 
is the unique subordinate fiber) 
and 
$e(X_{s_1, t}) = 1$. 
This equality holds exactly when 
$X_{s_1, t}$ 
is 
$m I_1 \, 
(m \geq 1)$. 
By Lemma \ref{det:lem:SubordSingularity} again, 
$X_{s_1, t}$ 
is a reduced curve, 
so 
$m = 1$. 
Accordingly
$X_{s_1, t}$ 
is 
$I_1$. 
\end{proof}

Lemma \ref{det:lem:SingleLefshetzFiber} (b) 
immediately yields the following:

\begin{prop}[Case: $e(X_0) - e(X_{0, t}) = 1$]  \label{det:prop:CaseSubordEuler1}
In each of the following barking families, 
the subordinate fiber is 
$I_1$. 
\begin{align*}\rm
&\boldsymbol{[II.1]} \quad II \barkarrow I_1 &\hspace{30pt}
 &\boldsymbol{[III^*.1]} \quad III^* \barkarrow IV^* \\ 
&\boldsymbol{[II.2]} \quad II \barkarrow I_1 &\hspace{30pt}
 &\boldsymbol{[III^*.3]} \quad III^* \barkarrow I_2^* \\
&\boldsymbol{[II^*.1]} \quad II^* \barkarrow III^* &\hspace{30pt}
 &\boldsymbol{[III^*.6]} \quad III^* \barkarrow I_2^* \\
&\boldsymbol{[II^*.5]} \quad II^* \barkarrow I_3^* &\hspace{30pt}
 &\boldsymbol{[III^*.9]} \quad III^* \barkarrow IV^* \\
&\boldsymbol{[II^*.6]} \quad II^* \barkarrow I_3^* &\hspace{30pt}
 &\boldsymbol{[IV.1]} \quad IV \barkarrow I_3 \\
&\boldsymbol{[II^*.8]} \quad II^* \barkarrow III^* &\hspace{30pt}
 &\boldsymbol{[IV^*.1]} \quad IV^* \barkarrow I_1^* \\
&\boldsymbol{[II^*.9]} \quad II^* \barkarrow III^* &\hspace{30pt}
 &\boldsymbol{[IV^*.4]} \quad IV^* \barkarrow I_1^* \\
&\boldsymbol{[III.1]} \quad III \barkarrow I_2 &\hspace{30pt}
 &\boldsymbol{[I_n^*.1]} \quad I_n^* \barkarrow I_{n-1}^*. \\
&\boldsymbol{[III.3]} \quad III \barkarrow I_2 
\end{align*}
\end{prop}

If 
$e(X_0) - e(X_{0, t}) \geq 2$, 
then 
we need another criterion to determine the subordinate fibers. 
However 
by Lemma \ref{det:lem:SingleLefshetzFiber} (a) 
we can narrow down candidates.

\begin{lem}[Case: $e(X_0) - e(X_{0, t}) = 2$]  \label{det:lem:CaseSubordEuler2}
In each of the following barking families, 
the set of subordinate fibers 
is one of 
$\left\{
II
\right\}$, 
$\left\{
I_2
\right\}$, 
and 
$\left\{
I_1, I_1
\right\}$. 
\begin{align*}\rm
&\boldsymbol{[II^*.2]} \quad II^* \barkarrow IV^* &\hspace{30pt}
 &\boldsymbol{[IV.3]} \quad IV \barkarrow I_2 \\
&\boldsymbol{[II^*.3]} \quad II^* \barkarrow I_2^* &\hspace{30pt}
 &\boldsymbol{[IV.4]} \quad IV \barkarrow II \\
&\boldsymbol{[II^*.7]} \quad II^* \barkarrow I_8 &\hspace{30pt}
 &\boldsymbol{[IV^*.2]} \quad IV^* \barkarrow I_0^* \\
&\boldsymbol{[III.2]} \quad III \barkarrow I_1 &\hspace{30pt}
 &\boldsymbol{[IV^*.3]} \quad IV^* \barkarrow I_6 \\
&\boldsymbol{[III^*.2]} \quad III^* \barkarrow I_1^* &\hspace{30pt}
 &\boldsymbol{[I_0^*.1]} \quad I_0^* \barkarrow I_4 \\
&\boldsymbol{[III^*.7]} \quad III^* \barkarrow I_7 &\hspace{30pt}
 &\boldsymbol{[I_n^*.2]} \quad I_n^* \barkarrow I_{n+4}. \\
&\boldsymbol{[IV.2]} \quad IV \barkarrow I_2 &\hspace{30pt}
\end{align*}
\end{lem}

\begin{lem}[Case: $e(X_0) - e(X_{0, t}) = 3$]  \label{det:lem:CaseSubordEuler3}
In each of the following barking families, 
the set of subordinate fibers 
is one of 
$\left\{
III
\right\}$, 
$\left\{
I_3
\right\}$, 
$\left\{
II, I_1
\right\}$, 
$\left\{
I_2, I_1
\right\}$, 
and 
$\left\{
I_1, I_1, I_1
\right\}$. 
\begin{align*}
&\boldsymbol{[III^*.4]} \quad III^* \barkarrow I_0^* &\hspace{30pt}
 &\boldsymbol{[I_0^*.2]} \quad I_0^* \barkarrow I_3. \\
&\boldsymbol{[III^*.5]} \quad III^* \barkarrow I_6 \\
&\boldsymbol{[III^*.8]} \quad III^* \barkarrow I_6 &\hspace{30pt}
\end{align*}
\end{lem}

\begin{lem}[Case: $e(X_0) - e(X_{0, t}) = 5$]  \label{det:lem:CaseSubordEuler5}
 The sum of the Euler characteristics of the subordinate fibers of the following barking family 
 is $5$: 
$$
\boldsymbol{[II^*.4]} \quad II^* \barkarrow I_5. 
$$
\end{lem}

\section{Monodromies around singular fibers}  \label{sec:Monodromies}

Next 
we consider 
the monodromies around singular fibers of splitting families 
(\emph{not necessarily} barking families).

Let 
$\pi : M \to \Delta$ 
be a (relatively) minimal degeneration of elliptic curves 
with the singular fiber 
$X_0$. 
We take a base point 
$s_0$
in 
$\Delta \setminus \left\{0\right\}$ 
and a loop (simple closed curve) 
$l_0$ 
in 
$\Delta \setminus \left\{0\right\}$ 
passing through the base point 
$s_0$
and 
circuiting around the origin with the counterclockwise orientation. 
Then 
$\pi^{-1}(l_0)$ 
is a real 3-manifold and the restriction 
$\pi : \pi^{-1}(l_0) \to l_0$ 
is a $\Sigma$-bundle over 
$S^1$, 
where 
$\Sigma$ 
is an elliptic curve. 
Here 
$\pi^{-1}(l_0)$ 
is obtained from 
$\Sigma \times [0, 1]$ 
by the identification of the boundaries 
$\Sigma \times \left\{0\right\}$ 
and 
$\Sigma \times \left\{1\right\}$ 
via an orientation-preserving homeomorphism 
$f$ 
of 
$\Sigma$.
The isotopy class 
$[f]$ 
of 
$f$ 
is called the \emph{topological monodromy} around 
$X_0$. 
Then 
$f$ 
induces an automorphism 
$\rho := f_\ast$ 
on 
$H_1(\Sigma, \Z)$, 
which is called 
the (\emph{homological}) \emph{monodromy} around 
$X_0$. 
Under an identification of 
$\Sigma$ 
and 
$\R^2 / \Z^2$, 
fixing a basis of 
$H_1(\Sigma, \Z)$, 
we obtain an isomorphism 
$$
\Aut(H_1(\Sigma, \Z)) \rightarrow SL(2,\Z). 
$$
In the subsequent discussion, 
we consider 
$\rho$ 
as an element of $SL(2, \Z)$.

Next suppose that 
$\Psi : \mathcal{M} \to \Delta \times \Delta^{\dagger}$ 
is a splitting family 
of the degeneration 
$\pi : M \to \Delta$, 
that is, 
the deformation 
$\pi_t : M_t \to \Delta_t$ 
of 
$\pi : M \to \Delta$ 
for a fixed 
$t \neq 0$ 
has singular fibers 
$X_1, X_2, \dots, X_N \, (N \geq 2)$. 
Then we say that
$X_0$ 
splits into 
$X_1, X_2, \dots, X_N$ 
and express 
$X_0 \longrightarrow X_1 + X_2 + \cdots + X_N$. 
Now we define the local monodromies 
around the singular fibers 
$X_k 
\, 
(k = 1, 2, \dots, N)$ 
as follows: 
Set 
$s_k := \pi_t(X_k)$. 
We take a base point 
$s'_0$ 
in 
$\Delta_t \setminus \left\{s_1, s_2, \dots, s_N\right\}$ 
(so the fiber 
$X_{s'_0} = \pi_t^{-1}(s'_0)$ 
is smooth). 
For each 
$k = 1, 2, \dots, N$, 
we take a loop 
$l_k$ 
in 
$\Delta_t \setminus \left\{s_1, s_2, \dots, s_N\right\}$ 
passing through the base point 
$s'_0$ 
and circuiting around 
$s_k$ 
with the counterclockwise orientation. 
Then the loop 
$l_k$ 
induces an orientation-preserving homeomorphism 
$f_k$ 
of 
$\Sigma$, 
which defines the local topological monodromy 
$[f_k]$ 
and 
the local (homological) monodromy 
$\rho_k \in SL(2, \Z)$ 
around 
$X_k$.

The following is known 
(see \cite{U}):

\begin{lem}   \label{det:lem:MonodromyFoundation}

The monodromy 
$\rho$ 
around 
$X_0$ 
$($%
resp. the local monodromy 
$\rho_k$ 
around 
$X_k$ 
for each 
$k = 1, 2, \dots, N$%
$)$ 
is conjugate to the standard monodromy matrix%
\footnote{
See Table \ref{ta:table:KodairaNotation} in Section \ref{sec:Ta}. 
}
corresponding to the singular fiber 
$X_0$ 
$($resp. 
$X_k)$. 

\end{lem}

Possibly after renumbering, 
we may assume that 
$l_1 \circ l_2 \circ \cdots \circ l_N$ 
is homotopic to 
a loop rounding all the singular values 
$s_1, s_2, \dots, s_N$ 
with the counterclockwise orientation. 
Let 
$\mathcal{D} \subset \Delta \times \Delta^{\dagger}$ 
be the set of singular values of 
$\Psi$. 
We now take a path 
$l$ 
in 
$\left(
\Delta \times \Delta^{\dagger}
\right) 
\setminus \mathcal{D}$ 
connecting 
$s_0 \in \Delta_0$ 
and 
$s'_0 \in \Delta_t$. 
Note that 
for any point 
$(s, t) \in l$, 
the fiber 
$X_{s, t} = \Psi^{-1}(s,t)$ 
is smooth. 
Since 
the loop 
$l^{-1} \circ l_1 \circ l_2 \circ \cdots \circ l_N \circ l$ 
is homotopic to the loop 
$l_0$, 
the topological monodromy 
$[f]$ 
is conjugate to 
the composition of the local topological monodromies 
$[f_1] \circ [f_2] \circ \cdots \circ [f_N]$. 
Similarly:

\begin{lem}   \label{det:lem:MonodromyDecompositon}

The monodromy 
$\rho$ 
is conjugate to the composition 
of the local monodromies 
$\rho_1, \rho_2, \dots, \rho_N$. 

\end{lem}

We prepare notation. 
$SL(2, \Z) 
= \left\langle a, b \, \big| \, a^3 = b^2 = -E \right\rangle$ 
is generated by
$$
a 
:= \left(
\begin{array}{cc}
0  & 1\\
-1 & 1
\end{array}
\right) 
\quad 
\textrm{and}
\quad 
b 
:= \left(
\begin{array}{cc}
0  & 1\\
-1 & 0
\end{array}
\right). 
$$
Setting%
\footnote{
The notations 
$s_0$ 
and 
$s_2$ 
are used in \cite{FM} Section 2.4, 
where 
`$s_1$' 
is defined as 
$s_1 := a b a$. 
}
$$
s_0 
:= a^{-1}b 
= \left(
\begin{array}{cc}
1 & 1\\
0 & 1
\end{array}
\right) 
\quad 
\textrm{and} 
\quad 
s_2 
:= b a^{-1} 
= \left(
\begin{array}{cc}
1  & 0\\
-1 & 1
\end{array}
\right), 
$$
then 
$s_0$ 
and 
$s_2$ 
are also generators of 
$SL(2, \Z)$: 
indeed we have 
$a 
= s_0 s_2$ 
and 
$b 
= s_0 s_2 s_0 
= s_2 s_0 s_2$. 
Since 
$s_2 
= (s_0 s_2) s_0 (s_0 s_2)^{-1}$, 
$s_2$ 
is conjugate to 
$s_0$.

Next 
we express the standard monodromy matrices of singular fibers 
as a product of 
$s_0$ 
and 
$s_2$ 
as follows (see \cite{U}): 
$$
\begin{array}{rcl}
A_{I_n} &=& \left(s_0\right)^n \, (n \geq 1) \\
A_{II} &=& s_0s_2 \\
A_{III} &=& s_0s_2s_0 = s_2s_0s_2 \\
A_{IV} &=& s_0s_2s_0s_2 \\
A_{II^*} &=& \left(s_0s_2\right)^5 \\
A_{III^*} &=& \left(s_0s_2\right)^4s_0 \\
A_{IV^*} &=& \left(s_0s_2\right)^4 \\
A_{I_n^*} &=& \left(s_0s_2\right)^3\left(s_0\right)^n  \, (n \geq 0). 
\end{array}
$$
The number of 
$s_0, s_2$ 
contained in each product coincides with 
the Euler characteristic of the corresponding singular fiber. 
Note that 
$s_0$ 
is the standard monodromy matrix 
$A_{I_1}$ 
of the singular fiber 
$I_1$. 
It is known that 
for any degeneration of elliptic curves 
except with 
$m I_0 \, (m \geq 2)$, 
the singular fiber splits into 
singular fibers of type 
$I_1$ 
(whose Euler characteristic 
$e(I_1)$ 
is equal to 1) 
after successive deformations. 
See \cite{Ka}, \cite{M}.

\begin{ex}
The barking family 
$\boldsymbol{[III.1]}$ 
splits the singular fiber 
$III$ 
into 
the main fiber 
$I_2$ 
and 
a subordinate fiber 
$I_1$: 
$$
III \longrightarrow I_2 + I_1. 
$$
Lemma \ref{det:lem:MonodromyDecompositon} states that, 
if $X_0$ 
splits into 
$X_1, X_2, \dots, X_N$, 
then a monodromy matrix of 
$X_0$ 
is conjugate to the composition of monodromy matrices of 
$X_1, X_2, \dots, X_N$ 
(that is, 
conjugacies of the standard monodromy matrices corresponding to 
$X_1, X_2, \dots, X_N$ 
respectively). 
In this case, 
the standard monodromy matrix 
$A_{III}$ 
of 
$III$ 
is decomposed into 
conjugacies of the standard monodromy matrices corresponding to 
$I_2$ 
and 
$I_1$: 
\begin{align*}
A_{III} 
&= s_0 s_2 s_0 &\\
&= s_0^2 (s_0^{-1} s_2 s_0) &\\
&= s_0^2 (s_2 s_0 s_2^{-1}) 
&(s_0 s_2 s_0 
= s_2 s_0 s_2)\\
&= A_{I_2} \cdot (s_2 A_{I_1} s_2^{-1}). &
\end{align*}
\end{ex}

In Section \ref{sec:Decomps}, 
we will give decompositions of the standard monodromy matrix 
corresponding to the splittings induced from 
Takamura's barking families.

\section{Constraints from monodromies}  \label{sec:Nonsplitting}

From Lemmas \ref{det:lem:MonodromyFoundation} and \ref{det:lem:MonodromyDecompositon}, 
it is a necessary condition for a singular fiber 
$X_0$ 
to split into singular fibers 
$X_1, X_2, \dots, X_N 
\, (N \geq 2)$ 
that 
some monodromy matrix of 
$X_0$ 
is conjugate to the composition of monodromy matrices of 
$X_1, X_2, \dots, X_N$, 
which means that 
monodromies give some constraints to splittings. 
In this section, 
we prove that 
none of the following splittings occurs 
(Theorem \ref{thm:NonsplittingThm}): 
\begin{align*}
IV \longrightarrow \ 
&I_2 + I_2, \\
II^* \longrightarrow \ 
&I_8 + II, \quad I_7 + III, \quad I_6 + IV, \\
&I_4 + I_0^*, \quad I_3 + I_1^*, \\
&I_u + I_v \, (u + v = 10), \\
III^* \longrightarrow \ 
&I_7 + II, \quad I_6 + III, \quad I_5 + IV, \quad I_3 + I_0^*, \\
&I_u + I_v \, (u + v = 9), \\
IV^* \longrightarrow \ 
&I_6 + II, \quad I_5 + III, \quad I_4 + IV, \quad I_2 + I_0^*, \\
&I_u + I_v \, (u + v = 8), \\
I_n^* \, (n \geq 0) &\longrightarrow \
I_{n + 4} + II, \quad I_{n + 3} + III, \quad I_{n + 2} + IV, \\
&I_u + I_v \, (u + v = n + 6 \textrm{ and } (n, u, v) \neq (2, 4, 4)). \\
I_0^* \longrightarrow \ 
&I_3 + I_2 + I_1. 
\end{align*}

We begin with preparation.

\begin{lem}  \label{lem:TrInv}
If matrices 
$A_1, A_2 \in SL(2, \Z)$ 
are conjugate, 
then 
$\Tr(A_1) = \Tr(A_2)$, 
where 
$\Tr(A_i)$ 
denotes the trace of 
$A_i$. 
\end{lem}

\begin{proof}
By assumption, 
we may write 
$A_1 = P A_2 P^{-1}$ 
for some 
$P \in SL(2, \Z)$. 
Hence 
$$
\Tr(A_1) 
= \Tr((P A_2) P^{-1}) 
= \Tr(P^{-1} (P A_2)) 
= \Tr(A_2). 
$$
\end{proof}

The following is useful:

\begin{lem}  \label{lem:TrCong1}
Suppose that 
a singular fiber 
$X$ 
splits into two singular fibers 
$I_n 
\, (n \geq 1)$ 
and 
$Y$: 
$$
X \longrightarrow I_n + Y. 
$$
Then 
$$
\Tr(A_X) \equiv \Tr(A_Y) \bmod n, 
$$
where 
$A_X$ 
and 
$A_Y$ 
are the standard monodromy matrices of 
$X$ 
and 
$Y$. 
\end{lem}

\begin{proof}
If 
$X$ 
splits into 
$I_n$ 
and 
$Y$, 
then for some monodromy matrix 
$C = 
\displaystyle 
\left(
\begin{array}{cc}
a & b \\
c & d
\end{array}
\right)$ 
of 
$Y$, 
$$
B 
:= A_{I_n} C 
= \left(
\begin{array}{cc}
1 & n \\
0 & 1
\end{array}
\right) 
\left(
\begin{array}{cc}
a & b \\
c & d
\end{array}
\right) 
= \left(
\begin{array}{cc}
a + n c & b + n d \\
c & d 
\end{array}
\right) 
$$
is a monodromy matrix of 
$X$. 
Then we have 
$$
\Tr(B) = a + n c + d 
= \Tr(C) + n c. 
$$ 
Thus 
$$
\Tr(B) \equiv \Tr(C) \bmod n.
$$
Where 
$A_X$ 
and 
$A_Y$ 
denote the standard monodromy matrices corresponding to 
$X$ 
and 
$Y$ 
respectively, 
$B$ 
is conjugate to 
$A_X$, 
while 
$C$ 
is conjugate to 
$A_Y$. 
By Lemma \ref{lem:TrInv} 
we have 
$\Tr(B) = \Tr(A_X)$ 
and 
$\Tr(C) = \Tr(A_Y)$. 
Accordingly 
$$
\Tr(A_X) \equiv \Tr(A_Y) \bmod n. 
$$
\end{proof}

We now consider the singular fiber 
$IV$. 
Since the Euler characteristic of 
$IV$ 
is 
$4$ 
and that of 
$I_2$ 
is 
$2$, 
$$
e(IV) = e(I_2) + e(I_2) 
$$
holds. 
Note that, 
if a singular fiber 
$X_0$ 
splits into two singular fibers 
$X_1$ 
and 
$X_2$, 
then 
$e(X_0) = e(X_1) + e(X_2)$ 
(Lemma \ref{det:lem:SplittingInv} (b)). 
So it is plausible that 
some deformation family splits the singular fiber 
$IV$ 
into two 
$I_2$. 
However this is not the case. 
If 
$IV$ 
splits into two 
$I_2$, 
by Lemma \ref{lem:TrCong1}, 
we have 
$$
\Tr(A_{IV}) \equiv \Tr(A_{I_2}) \bmod 2, 
$$
which contradicts that 
$\Tr(A_{IV}) = -1$ 
and 
$\Tr(A_{I_2}) = 2$. 
Thus the splitting 
$$
IV \longrightarrow I_2 + I_2 
$$
does not occur. 
We have shown 
the first statement of the following lemma, 
and we can show the others by the same argument:

\begin{lem}  \label{lem:Nonsplitting1}
\begin{description}
 \item[\rm(a)] 
The singular fiber 
$IV$ 
never splits as follows: 
$$
IV \longrightarrow I_2 + I_2. 
$$

 \item[\rm(b)] 
The singular fiber 
$II^*$ 
never splits as follows: 
\begin{align*}
II^* \longrightarrow \
&I_8 + II, \quad I_7 + III, \quad I_6 + IV, \\
&I_4 + I_0^*, \quad I_3 + I_1^*, \quad \\
&I_u + I_v \, (u + v = 10). 
\end{align*}

 \item[\rm(c)] 
The singular fiber 
$III^*$ 
never splits as follows: 
\begin{align*}
III^* \longrightarrow \
&I_7 + II, \quad I_6 + III, \quad I_5 + IV, \\
&I_3 + I_0^*, \quad I_u + I_v \, (u + v = 9). 
\end{align*}

 \item[\rm(d)] 
The singular fiber 
$IV^*$ 
never splits as follows: 
\begin{align*}
IV^* \longrightarrow \
&I_6 + II, \quad I_5 + III, \quad I_4 + IV, \\
&I_2 + I_0^*, \quad I_u + I_v \, (u + v = 8). 
\end{align*}

 \item[\rm(e)] 
The singular fiber 
$I_n^* \, (n \geq 1)$ 
never splits as follows: 
\begin{align*}
I_n^* \longrightarrow \
&I_{n + 4} + II, \quad I_{n + 3} + III, \quad I_{n + 2} + IV, \\
&I_u + I_v 
\quad (u + v = n + 6, \ (n, u, v) \neq (2, 4, 4)). 
\end{align*}
\end{description}
\end{lem}

Next 
we consider splittings of
$I_0^*$. 
The standard monodromy matrix of 
$I_0^*$ 
is 
$A_{I_0^*} = -E$, 
where 
$E$ 
is the identity matrix.

\begin{lem}  \label{lem:TrEq}
Suppose that 
the singular fiber 
$I_0^*$ 
splits into two singular fibers 
$X$ 
and 
$Y$: 
$$
I_0^* \longrightarrow X + Y. 
$$
Then 
$$
\Tr(A_X) + \Tr(A_Y) = 0. 
$$
\end{lem}

\begin{proof}
If 
$I_0^*$ 
splits into 
$X$ 
and 
$Y$, 
then for monodromy matrices 
$B$ 
and 
$C$ 
of 
$X$ 
and 
$Y$, 
we have 
$A_{I_0^*} 
= B C$,
where 
$A_{I_0^*}$ 
is the standard monodromy matrix of 
$I_0^*$. 
Since 
$A_{I_0^*} 
= -E$, 
we have 
$-E 
= B C$,
that is, 
$B 
= -C^{-1}$. 
In particular, 
$$
\Tr(B) = - \Tr(C). 
$$ 
Since 
$B$ 
(resp. 
$C$) 
is conjugate to 
$A_X$ 
(resp. 
$A_Y$), 
by Lemma \ref{lem:TrInv} 
we have 
$\Tr(B) = \Tr(A_X)$ 
and 
$\Tr(C) = \Tr(A_Y)$. 
Thus 
$$
\Tr(A_X) + \Tr(A_Y) = 0. 
$$
\end{proof}

\begin{lem}  \label{lem:TrCong2}
Suppose that 
the singular fiber 
$I_0^*$ 
splits into three singular fibers 
$I_n 
\, (n \geq 1)$, 
$X$ 
and 
$Y$: 
$$
I_0^* \longrightarrow I_n + X + Y. 
$$
Then 
$$
\Tr(A_X) + \Tr(A_Y) \equiv 0 \bmod n. 
$$
\end{lem}

\begin{proof}
If 
$I_0^*$ 
splits into 
$I_3$, $X_1$ 
and 
$X_2$, 
then for monodromy matrices 
$B$ 
and 
$C$ 
of 
$X$ 
and 
$Y$, 
we have 
$A_{I_0^*} 
= A_{I_n} B C$. 
Since 
$A_{I_0^*} = -E$ 
and 
$A_{I_n} = 
\displaystyle 
\left(
\begin{array}{cc}
1 & n \\
0 & 1
\end{array}
\right)$, 
writing 
$B = 
\displaystyle 
\left(
\begin{array}{cc}
a & b \\
c & d
\end{array}
\right)$, 
$$
-C^{-1} 
= A_{I_n} B 
= \left(
\begin{array}{cc}
1 & n \\
0 & 1
\end{array}
\right) 
\left(
\begin{array}{cc}
a & b \\
c & d
\end{array}
\right) 
= \left(
\begin{array}{cc}
a + n c & b + n d \\
c & d 
\end{array}
\right). 
$$
Then we have 
$$
-\Tr(C) 
= a + n c + d 
= \Tr(B) + n c. 
$$ 
Thus 
$$
\Tr(B) + \Tr(C) \equiv 0 \bmod n.
$$
Since 
$B$ 
(resp. 
$C$) 
is conjugate to 
$A_X$ 
(resp. 
$A_Y$), 
by Lemma \ref{lem:TrInv} 
we have 
$\Tr(B) = \Tr(A_X)$ 
and 
$\Tr(C) = \Tr(A_Y)$. 
Accordingly 
$$
\Tr(A_X) + \Tr(A_Y) \equiv 0 \bmod n.
$$
\end{proof}

\begin{lem}  \label{lem:Nonsplitting2}
\begin{description}
 \item[\rm(a)] 
The singular fiber 
$I_0^*$ 
never splits as follows: 
\begin{align*}
I_0^* \longrightarrow \
&I_4 + II, \quad I_3 + III, \quad I_2 + IV, \\
&I_5 + I_1, \quad I_4 + I_2, \quad I_3 + I_3. 
\end{align*}
 \item[\rm(b)] 
The singular fiber 
$I_0^*$ 
never splits as follows: 
$$
I_0^* \longrightarrow I_3 + I_2 + I_1. 
$$
\end{description}
\end{lem}

\begin{proof}
(a) we only show that 
the splitting 
$I_0^* \longrightarrow \ I_4 + II$ 
does not occur, 
because we can give the proof for the other splittings 
by the same argument. 
If 
$I_0^*$ 
splits into 
$I_4$ 
and 
$II$, 
by Lemma \ref{lem:TrEq}, 
we have 
$$
\Tr(A_{I_4}) + \Tr(A_{II}) = 0, 
$$
which contradicts that 
$\Tr(A_{I_4}) = 2$ 
and 
$\Tr(A_{II}) = 1$. 
Thus the splitting 
$$
I_0^* \longrightarrow I_4 + II
$$
does not occur.

(b) 
If 
$I_0^*$ 
splits into 
$I_3$, 
$I_2$ 
and 
$I_1$, 
by Lemma \ref{lem:TrCong2}, 
we have 
$$
\Tr(A_{I_2}) + \Tr(A_{I_1}) \equiv 0 \bmod 3. 
$$
which contradicts that 
$\Tr(A_{I_2}) = \Tr(A_{I_1}) = 2$. 
Thus the splitting 
$$
I_0^* \longrightarrow I_3 + I_2 + I_1. 
$$
does not occur. 
\end{proof}

\begin{rem}
We can give an alternative proof of Lemma \ref{lem:Nonsplitting2} (a) 
\emph{except for the splitting} 
$I_0^* \longrightarrow \ I_4 + I_2$ 
as follows; 
For instance, 
suppose that 
$I_0^*$ 
splits into 
$I_4$ 
and 
$II$. 
By Lemma \ref{lem:TrCong1}, 
we then have 
$$
\Tr(A_{I_0^*}) \equiv \Tr(A_{II}) \bmod 4, 
$$
which contradicts that 
$\Tr(A_{I_0^*}) = -2$ 
and 
$\Tr(A_{II}) = 1$. 
Thus the splitting 
$$
I_0^* \longrightarrow \ I_4 + II 
$$
does not occur. 
\end{rem}

We summarize Lemmas \ref{lem:Nonsplitting1} and \ref{lem:Nonsplitting2} 
as follows:

\begin{thm}  \label{thm:NonsplittingThm}
None of the following splittings occurs: 
\begin{align*}
IV \longrightarrow \ 
&I_2 + I_2, \\
II^* \longrightarrow \ 
&I_8 + II, \quad I_7 + III, \quad I_6 + IV, \\
&I_4 + I_0^*, \quad I_3 + I_1^*, \\
&I_u + I_v \, (u + v = 10), \\
III^* \longrightarrow \ 
&I_7 + II, \quad I_6 + III, \quad I_5 + IV, \quad I_3 + I_0^*, \\
&I_u + I_v \, (u + v = 9), \\
IV^* \longrightarrow \ 
&I_6 + II, \quad I_5 + III, \quad I_4 + IV, \quad I_2 + I_0^*, \\
&I_u + I_v \, (u + v = 8), \\
I_n^* \, (n \geq 0) &\longrightarrow \
I_{n + 4} + II, \quad I_{n + 3} + III, \quad I_{n + 2} + IV, \\
&I_u + I_v \, (u + v = n + 6 \textrm{ and } (n, u, v) \neq (2, 4, 4)). \\
I_0^* \longrightarrow \ 
&I_3 + I_2 + I_1. 
\end{align*}
\end{thm}

\section{Determination of subordinate fibers, $1$}  \label{sec:Det1}

In this section, 
based on the result of the previous section, 
we determine the subordinate fibers of Takamura's barking families 
$\boldsymbol{[II^*.7]}$, 
$\boldsymbol{[III^*.7]}$, 
$\boldsymbol{[IV^*.3]}$, 
$\boldsymbol{[I_0^*.1]}$, 
$\boldsymbol{[I_n^*.2]}$.

\begin{prop}  \label{det:prop:Splitting2.7}
The barking family $\boldsymbol{[II^*.7]}$ splits the singular fiber 
$II^*$ 
as follows: 
$$
II^* \longrightarrow I_8 \, + \, I_1 + I_1, 
$$
where 
$I_8$ 
is the main fiber and the two 
$I_1$ 
are subordinate fibers. 
\end{prop}

\begin{proof}
In the barking family $\boldsymbol{[II^*.7]}$, 
$II^*$ 
is deformed to 
$I_8$: 
$$
II^* \barkarrow I_8. 
$$
By Lemma \ref{det:lem:CaseSubordEuler2}, 
the set of subordinate fibers 
is one of 
(i) 
$\left\{
II
\right\}$, 
(ii)
$\left\{
I_2
\right\}$, 
and 
(iii) 
$\left\{
I_1, I_1
\right\}$. 
Now Lemma \ref{lem:Nonsplitting1} (b) eliminates the cases (i) and (ii). 
Thus the subordinate fibers are two 
$I_1$. 
\end{proof}

\begin{prop}  \label{det:prop:Splitting4.7}
The barking family $\boldsymbol{[III^*.7]}$ splits the singular fiber 
$III^*$ 
as follows: 
$$
III^* \longrightarrow I_7 \, + \, I_1 + I_1, 
$$
where 
$I_7$ 
is the main fiber and the two 
$I_1$ 
are subordinate fibers.
\end{prop}

\begin{proof}
In the barking family $\boldsymbol{[III^*.7]}$, 
$III^*$ 
is deformed to 
$I_7$: 
$$
III^* \barkarrow I_7. 
$$
By Lemma \ref{det:lem:CaseSubordEuler2}, 
the set of subordinate fibers 
is one of 
(i) 
$\left\{
II
\right\}$, 
(ii)
$\left\{
I_2
\right\}$, 
and 
(iii) 
$\left\{
I_1, I_1
\right\}$. 
Now Lemma \ref{lem:Nonsplitting1} (c) eliminates the cases (i) and (ii). 
Thus the subordinate fibers are two 
$I_1$. 
\end{proof}

\begin{prop}  \label{det:prop:Splitting6.3}
The barking family $\boldsymbol{[IV^*.3]}$ splits the singular fiber 
$IV^*$ 
as follows: 
$$
IV^* \longrightarrow I_6 \, + \, I_1 + I_1, 
$$
where 
$I_6$ 
is the main fiber and the two 
$I_1$ 
are subordinate fibers.
\end{prop}

\begin{proof}
In the barking family $\boldsymbol{[IV^*.3]}$, 
$IV^*$ 
is deformed to 
$I_6$: 
$$
IV^* \barkarrow I_6. 
$$
By Lemma \ref{det:lem:CaseSubordEuler2}, 
the set of subordinate fibers 
is one of 
(i) 
$\left\{
II
\right\}$, 
(ii)
$\left\{
I_2
\right\}$, 
and 
(iii) 
$\left\{
I_1, I_1
\right\}$. 
Now Lemma \ref{lem:Nonsplitting1} (d) eliminates the cases (i) and (ii). 
Thus the subordinate fibers are two 
$I_1$. 
\end{proof}

\begin{prop}  \label{det:prop:Splitting7.1}
The barking family $\boldsymbol{[I_0^*.1]}$ splits the singular fiber 
$I_0^*$ 
as follows: 
$$
I_0^* \longrightarrow I_4 \, + \, I_1 + I_1, 
$$
where 
$I_4$ 
is the main fiber and the two 
$I_1$ 
are subordinate fibers.
\end{prop}

\begin{proof}
In the barking family $\boldsymbol{[I_0^*.1]}$,
$I_0^*$ 
is deformed to 
$I_4$: 
$$
I_0^* \longrightarrow I_4. 
$$
By Lemma \ref{det:lem:CaseSubordEuler2}, 
the set of subordinate fibers 
is one of 
(i) 
$\left\{
II
\right\}$, 
(ii)
$\left\{
I_2
\right\}$, 
and 
(iii) 
$\left\{
I_1, I_1
\right\}$. 
Now Lemma \ref{lem:Nonsplitting2} (a) eliminates the cases (i) and (ii). 
Thus the subordinate fibers are two 
$I_1$. 
\end{proof}

\begin{prop}  \label{det:prop:Splitting8.2}
The barking family $\boldsymbol{[I_n^*.2]}$ splits the singular fiber 
$I_n^*$ 
as follows: 
$$
I_n^* \longrightarrow I_{n + 4} \, + \, I_1 + I_1, 
$$
where 
$I_{n + 4}$ 
is the main fiber and the two 
$I_1$ 
are subordinate fibers. 
\end{prop}

\begin{proof}
In the barking family $\boldsymbol{[I_n^*.2]}$,
$I_n^*$ 
is deformed to 
$I_{n+4}$: 
$$
I_n^* \barkarrow I_{n+4}. 
$$
By Lemma \ref{det:lem:CaseSubordEuler2}, 
the set of subordinate fibers 
is one of 
(i) 
$\left\{
II
\right\}$, 
(ii)
$\left\{
I_2
\right\}$, 
and 
(iii) 
$\left\{
I_1, I_1
\right\}$. 
Now Lemma \ref{lem:Nonsplitting1} (e) eliminates the cases (i) and (ii). 
Thus the subordinate fibers are two 
$I_1$. 
\end{proof}

\begin{rem}  \label{rem:UndeterminedCases}
For the barking families 
$\boldsymbol{[IV.3]}$, $\boldsymbol{[III^*.8]}$, $\boldsymbol{[I_0^*.2]}$, 
we cannot determine the subordinate fibers 
but we can narrow down candidates: 
\begin{itemize}
 \item 
The splitting of $IV$ induced from the barking family $\boldsymbol{[IV.3]}$ 
is one of the following: 
\begin{align*}
&IV \longrightarrow I_2 \, + \, II, \\
&IV \longrightarrow I_2 \, + \, I_1 + I_1. 
\end{align*}
In fact, 
by Lemma \ref{det:lem:CaseSubordEuler2}, 
the set of subordinate fibers 
is one of 
(i) 
$\left\{
II
\right\}$, 
(ii)
$\left\{
I_2
\right\}$, 
and 
(iii) 
$\left\{
I_1, I_1
\right\}$, 
and Lemma \ref{lem:Nonsplitting1} (a) eliminates the case (ii). 

 \item 
The splitting of $III^*$ induced from the barking family $\boldsymbol{[III^*.8]}$ 
is one of the following: 
\begin{align*}
&III^* \longrightarrow I_6 \, + \, II + I_1, \\
&III^* \longrightarrow I_6 \, + \, I_2 + I_1, \\
&III^* \longrightarrow I_6 \, + \, I_1 + I_1 + I_1. 
\end{align*}
In fact, 
by Lemma \ref{det:lem:CaseSubordEuler3}, 
the set of subordinate fibers 
is one of 
(i) 
$\left\{
III
\right\}$, 
(ii) 
$\left\{
I_3
\right\}$, 
(iii) 
$\left\{
II, I_1
\right\}$, 
(iv) 
$\left\{
I_2, I_1
\right\}$, 
and 
(v) 
$\left\{
I_1, I_1, I_1
\right\}$, 
and Lemma \ref{lem:Nonsplitting1} (c) eliminates the cases (i) and (ii). 

 \item 
The splitting of $I_0^*$ induced from the barking family $\boldsymbol{[I_0^*.2]}$ 
is one of the following: 
\begin{align*}
&I_0^* \longrightarrow I_3 \, + \, II + I_1, \\
&I_0^* \longrightarrow I_3 \, + \, I_1 + I_1 + I_1. 
\end{align*}
In fact, 
by Lemma \ref{det:lem:CaseSubordEuler3}, 
the set of subordinate fibers 
is one of 
(i) 
$\left\{
III
\right\}$, 
(ii) 
$\left\{
I_3
\right\}$, 
(iii) 
$\left\{
II, I_1
\right\}$, 
(iv) 
$\left\{
I_2, I_1
\right\}$, 
and 
(v) 
$\left\{
I_1, I_1, I_1
\right\}$, 
and Lemma \ref{lem:Nonsplitting2} eliminates the cases (i), (ii) and (iv). 
\end{itemize}
\end{rem}

\section{Singularities near proportional subbranches}  \label{sec:PropSing}

Let 
$\pi : M \to \Delta$ 
be a linear degeneration of complex curves 
with a stellar singular fiber 
$X_0 = m_0 \Theta_0 + \sum_{j = 1}^h \Br^{(j)}$. 
If there exists a simple crust 
$Y$ 
of 
$X_0$, 
then we can construct a splitting family of 
$\pi : M \to \Delta$, 
which is called a barking family associated with 
$Y$ 
(Theorem \ref{ta:thm:Ta3}). 
Suppose that 
$Y = n_0 \Theta_0 + \sum_{j = 1}^h \br^{(j)}$ 
is a simple crust of 
$X_0$ 
with barking multiplicity 
$l$.

Recall that
each subbranch of 
$Y$ 
is of type $A_l$, $B_l$ or $C_l$. 
A subbranch 
$\overline{\br}^{(j)}$ 
is said to be \emph{proportional} 
if 
$m_0 n_1^{(j)} = n_0 m_1^{(j)}$ 
(equivalently 
$n_0 / m_0 = n_1^{(j)} / m_1^{(j)} 
= \cdots = n_\nu^{(j)}/m_\nu^{(j)}$). 
Note that 
every proportional subbranch of simple crusts is of type $A_l$. 
Indeed, 
any proportional subbranch of type $B_l$ is of type $A_l$,
and no proportional subbranch is of type $C_l$. 
Moreover 
every proportional subbranch 
$\overline{\br}^{(j)}$ 
has the same length as that of 
$\overline{\Br}^{(j)}$ 
(that is, 
$\nu^{(j)} = \lambda^{(j)}$) 
and 
satisfies 
$n_{\lambda^{(j)} +1} = 0$.

The following lemma is important 
(\cite{Ta3} Proposition 16.2.6):

\begin{lem}  \label{det:lem:SubordSingularity}
Suppose that
$\Psi : \mathcal{M} \to \Delta \times \Delta^{\dagger}$ 
is a barking family of the degeneration 
$\pi : M \to \Delta$ 
associated with a simple crust 
$Y$. 
Then any subordinate fiber of 
$\Psi$ 
is a reduced curve only with $A$-singularities%
\footnote{
An $A$-singularity 
is a singularity analytically equivalent to 
$y^2 = x^{\mu + 1}$ 
for some positive integer 
$\mu$. 
}. 
Moreover 
these singularities lie 
{\rm(i)} near the core or 
{\rm(ii)} near the edge%
\footnote{
To be precise, 
near the `terminal' irreducible component 
$\Theta_{\lambda^{(j)}}$ 
of the branch 
$\Br^{(j)}$ 
corresponding to each proportional subbranch 
$\br^{(j)}$. 
} 
of each proportional subbranch 
if it exists. 
\end{lem}

\begin{rem}
By Lemma \ref{det:lem:SubordSingularity}, 
every subordinate fiber in barking families is 
a reduced curve only with isolated singularities. 
In particular, 
for degenerations of elliptic curves, 
none of 
$m I_n \, (m \geq 2)$, $IV^*$, $III^*$, $II^*$, $m I_n^* \, (m \geq 2)$ 
appears 
as a subordinate fiber.
\end{rem}

The rest of this section investigates 
the singularities of subordinate fibers 
near a proportional subbranch. 
Let 
$\pi : M \to \Delta$ 
be a linear degeneration of complex curves 
with a stellar singular fiber 
$X_0$ 
and 
$\Psi : \mathcal{M} \to \Delta \times \Delta^{\dagger}$ 
be a barking family
associated with a simple crust 
$Y$ 
with barking multiplicity 
$l$. 
Suppose that 
$Y$ 
has a proportional subbranch 
$\br$ 
of a branch 
$\Br$ 
of 
$X_0$. 
First recall that
near the branch 
$\Br$, 
$\mathcal{M}$ 
is given by the following data (see \cite{Ta3} Chapter 7): 
for 
$i = 1, 2, \dots, \lambda$, 
$$
\begin{cases}
\mathcal{H}_i : 
w_i^{m_{i - 1} - l n_{i - 1}} \eta_i^{m_i - l n_i} 
\left( w_i^{n_{i - 1}} \eta_i^{n_i} + t^d f_i \right)^l - s = 0, &\\
&\hspace{-40pt} \textrm{in} \ U_i \times \C \times \Delta \times \Delta^{\dagger}, \\
\mathcal{H}'_i : 
z_i^{m_{i + 1} - l n_{i + 1}} \zeta_i^{m_i - l n_i} 
\left( z_i^{n_{i + 1}} \zeta_i^{n_i} + t^d \hat{f}_i \right)^l - s = 0, &\\
&\hspace{-40pt} \textrm{in} \ V_i \times \C \times \Delta \times \Delta^{\dagger}. 
\end{cases}
$$
Note that, 
substituting 
$t = 0$ 
into these equations, 
we obtain 
$$
\begin{cases}
\mathcal{H}_i \big|_{t = 0} : 
w_i^{m_{i - 1}} \eta_i^{m_i} - s = 0, \\
\mathcal{H}'_i \big|_{t = 0} : 
z_i^{m_{i+1}} \zeta_i^{m_i} - s = 0, \\
\end{cases}
$$
which are the local expressions of 
$M$ 
near 
$\Br$. 
See the proof of Lemma \ref{ta:lem:BrExistence}. 
For a fixed 
$(s, t) \in \Delta \times \Delta^{\dagger}$, 
we consider the fiber 
$X_{s, t} = \Psi^{-1}(s, t)$ 
of 
$\Psi$. 
The following is required 
(\cite{Ta3} Section 7.2):

\begin{lem}  \label{det:lem:HypersurfSing} 
Let 
$m, n, l$ 
be positive integers satisfying 
$m - l n > 0$ 
and 
$m', n'$ 
be nonnegative integers satisfying 
$m' - l n' \geq 0$. 
Set 
$h(z, \zeta) := f(z^{p'}\zeta^p)$ 
for a non-vanishing holomorphic function 
$f$ 
and positive integers 
$p, p' \, (p < p')$. 
Then 
a complex curve 
$C_{s, t}$ 
in 
$\C^2$ 
defined by 
$$
C_{s, t} \, : \, 
z^{m' - l n'} \zeta^{m - l n} \left( z^{n'} \zeta^n + t h \right)^l - s 
= 0 
$$
is singular if and only if
\begin{description}
\item[\rm(i)] 
$s = 0$ 
or
\item[\rm(ii)] 
$m' = n' = 0$ 
and 
$\displaystyle \left(\frac{l n - m}{l n}\right)^{l\bar{n}} s^{\bar{n}} 
= \left(\frac{l n - m}{m} t c\right)^{\bar{m}}$, 
\end{description}
where 
$c := h(0, 0)$ 
and 
$\bar{m}$ 
and 
$\bar{n}$ 
are the relatively prime integers satisfying 
$\bar{n} / \bar{m} = n / m$. 
In the case {\rm(ii)}, 
$(z, \zeta) \in C_{s, t}$ 
is a singularity exactly when 
$$
z = 0 \quad \textrm{ and } \quad \zeta^n = \frac{l n - m}{m} t c. 
$$
\end{lem}

Since 
$\br$ 
is proportional, 
we have 
$m_{\lambda +1} = n_{\lambda +1} = 0$, 
so 
$$
\mathcal{H}'_\lambda\big|_{s, t} : 
\zeta_\lambda^{m_\lambda - l n_\lambda} 
\left( \zeta_\lambda^{n_\lambda} + t^d \hat{f}_\lambda \right)^l - s 
= 0. 
$$
Lemma \ref{det:lem:HypersurfSing} ensures that 
for some 
$(s, t) \, (s, t \neq 0)$, 
the curve
$\mathcal{H}'_\lambda\big|_{s, t}$ 
has singularities. 
In what follows, 
we write 
$m := m_\lambda$ 
and 
$n := n_\lambda$, 
and denote by 
$\bar{m}$ 
and 
$\bar{n}$ 
the relatively prime integers 
satisfying 
$\bar{n} / \bar{m} = n / m$.

For a fixed 
$t \neq 0$, 
the equation 
$$
\left(\frac{l n - m}{l n}\right)^{l \bar{n}} s^{\bar{n}} 
= \left(\frac{l n - m}{m} t^d c\right)^{\bar{m}} 
$$
for 
$s$ 
has 
$\bar{n}$ 
solutions, 
say, 
$s_1, s_2, \dots, s_{\bar{n}}$. 
Since 
$(0, \zeta)$ 
satisfying 
$\zeta^n 
= \frac{l n - m}{m} t^d c$ 
is a singularity of 
$\mathcal{H}'_\lambda\big|_{s_k, t}$ 
for some 
$s_k$, 
each 
$\mathcal{H}'_\lambda\big|_{s_k, t}$ 
has 
$n / \bar{n} \, (= \gcd(m, n))$ 
singularities.

The above result is summarized as follows:

\begin{prop}  \label{det:prop:PropSingularityNum}
Let 
$\pi : M \to \Delta$ 
be a linear degeneration of complex curves 
with a stellar singular fiber 
$X_0$
and let 
$\Psi : \mathcal{M} \to \Delta \times \Delta^{\dagger}$ 
be a barking family 
associated with a simple crust 
$Y$ 
with barking multiplicity 
$l$. 
Suppose that 
$Y$ 
has a proportional subbranch 
$\br^{(j)}$ 
of a branch 
$\Br^{(j)}$ 
of 
$X_0$. 
Write 
$\overline{\Br}^{(j)} 
:= m_0 \Delta_0 + m_1 \Theta_1 + m_2 \Theta_2 + \cdots + m_\lambda \Theta_\lambda$ 
and 
$\overline{\br}^{(j)} 
:= n_0 \Delta_0 + n_1 \Theta_1 + n_2 \Theta_2 + \cdots + n_\lambda \Theta_\lambda$ 
and let 
$\bar{m}$ 
and 
$\bar{n}$ 
be 
the relatively prime positive integers satisfying 
$\bar{n} / \bar{m} = n_\lambda / m_\lambda$. 
Then 
in the deformation 
$\pi_t : M_t \to \Delta_t$ 
for a fixed 
$t \neq 0$, 
there exist 
$\bar{n}$ 
subordinate fibers 
that have singularities 
near the edge of 
$\Br^{(j)}$. 
Moreover, 
each of these subordinate fibers has 
$n / \bar{n} \, (= \gcd(m, n))$ 
singularities near the edge of 
$\Br^{(j)}$. 
\end{prop}

\section{Singularities near the core}  \label{sec:CoreSing}

We next investigate the singularities of subordinate fibers near the core.

Let 
$\pi : M \to \Delta$ 
be a linear degeneration of complex curves 
with a stellar singular fiber 
$X_0 = m_0 \Theta_0 + \sum_{j = 1}^h \Br^{(j)}$ 
and let 
$\Psi : \mathcal{M} \to \Delta \times \Delta^{\dagger}$ 
be a barking family of the degeneration 
$\pi : M \to \Delta$ 
associated with a simple crust 
$Y = n_0 \Theta_0 + \sum_{j = 1}^h \br^{(j)}$. 
Write 
$\overline{\Br}^{(j)} 
= m_0 \Delta_0^{(j)} + m_1^{(j)} \Theta_1^{(j)} + 
\cdots + m_{\lambda^{(j)}}^{(j)} \Theta_{\lambda^{(j)}}^{(j)}$, 
$\overline{\br}^{(j)} 
= n_0 \Delta_0^{(j)} + n_1^{(j)} \Theta_1^{(j)} + 
\cdots + n_{\nu^{(j)}}^{(j)} \Theta_{\nu^{(j)}}^{(j)}$ 
and let 
$p^{(j)}$ 
be the attachment point on 
$\Theta_0$ 
with 
$\Br^{(j)}$. 
For brevity, 
we assume that 
\emph{the subbranches 
$\overline{\br}^{(1)}, \overline{\br}^{(2)}, \dots, \overline{\br}^{(v)}$ 
are proportional 
and 
$\overline{\br}^{(v + 1)}, \overline{\br}^{(v + 2)}, \dots, \overline{\br}^{(h)}$ 
are not}.

Let 
$N_0$ 
be the normal bundle 
of 
$\Theta_0$ 
in 
$M$. 
Recall that 
the local expression of 
$\mathcal{M}$ 
near the core 
$\Theta_0$ 
is given by 
$$
\sigma(z) \zeta^{m_0} - s 
+ \sum_{k = 1}^l {}_l\mathrm{C}_k t^{k d} \sigma(z) \tau(z)^k \zeta^{m_0 - k n_0} = 0 
\quad \textrm{in} \ N_0 \times \Delta \times \Delta^{\dagger}, 
$$
equivalently 
$$
\sigma(z) \zeta^{m_0 - l n_0} \left(\zeta^{n_0} + t^d \tau(z)\right)^l - s = 0, 
$$
where 
$\sigma$ 
is the standard section of 
$N_0^{\otimes (-m_0)}$ 
and 
$\tau$ 
is a core section of 
$N_0^{\otimes n_0}$ 
for 
$Y$ 
(see \cite{Ta3} Chapter 16). 
Substituting 
$t = 0$ 
into this equation, 
we obtain 
$$
\sigma(z) \zeta^{m_0} - s = 0 
\quad \textrm{in} \ N_0 \times \Delta \times \left\{0\right\}, 
$$
which is the local expression of 
$M$ 
around 
$\Theta_0$. 
See the paragraph subsequent to Remark \ref{ta:rem:TopLinearizability}. 
Note that
$\sigma$ 
has a zero of order 
$m_1^{(j)}$ 
at 
$p^{(j)}$, 
while 
$\tau$ 
has a pole of order 
$n_1^{(j)}$ 
at 
$p^{(j)}$. 
Suppose that 
$\tau$ 
has a zero of order 
$a_i$ 
at 
$q_i \, (i = 1, 2, \dots, k)$ 
on 
$\Theta_0$.

Fixing 
$s, t \neq 0$, 
consider a fiber 
$X_{s, t} := \Psi^{-1}(s, t)$ 
of 
$\Psi : \mathcal{M} \to \Delta \times \Delta^{\dagger}$. 
Set 
$F := 
\sigma(z) \zeta^{m_0 - l n_0} \left(\zeta^{n_0} + t^d \tau(z)\right)^l$. 
Then 
$(z, \zeta) \in X_{s, t}$ 
is a singularity 
if and only if 
$$
\frac{\partial}{\partial z}F(z, \zeta) = \frac{\partial}{\partial \zeta}F(z, \zeta) = 0, 
$$
equivalently
$$
\begin{cases}
\zeta^{m_0 - l n_0} 
\left(\zeta^{n_0} + t^d \tau(z)\right)^{l - 1} 
\left\{\sigma_z(z) \zeta^{n_0} 
+ t^d \left(\sigma_z(z) \tau(z) + l \sigma(z) \tau_z(z)\right) \right\} &\\
&\hspace{-25pt}= 0, \\
\zeta^{m_0 - l n_0 - 1} 
\left(\zeta^{n_0} + t^d \tau(z)\right)^{l - 1} \sigma(z) 
\left(m_0 \zeta^{n_0} + \left(m_0 - l n_0\right) t^d \tau(z)\right) 
&\hspace{-25pt}= 0, 
\end{cases}
$$
where 
$\sigma_z := \frac{d}{dz}\sigma$ 
and 
$\tau_z = \frac{d}{dz}\tau$. 
Set 
$K(z) 
:= n_0 \sigma_z(z) \tau(z) + m_0 \sigma(z) \tau_z(z)$, 
which is called the \emph{plot function}%
\footnote{
Note that 
$K(z)$ 
is \emph{not} a function on 
$\Theta_0$ 
but a meromorphic section of a line bundle 
$N_0^{\otimes(n - m)} \otimes \Omega_{\Theta_0}^1$ 
on 
$\Theta_0$, 
where 
$\Omega_{\Theta_0}^1$ 
is the cotangent bundle of 
$\Theta_0$. 
}. 
Then 
the above equations hold 
precisely when 
$$
\begin{cases}
K(z) 
= 0, 
\quad 
\sigma(z) 
\neq 0, 
\quad 
\tau(z) 
\neq 0, \\
\zeta^{n_0} 
= \frac{l n_0 - m_0}{m_0} t^d \tau(z). 
\end{cases}
$$
In particular, 
whether 
$(z, \zeta) \in X_{s, t}$ 
is a singularity 
does \emph{not} depend on 
$s$. 
Noting that 
every point 
$(z, \zeta)$ 
in 
$X_{s, t}$ 
satisfies 
$$
\sigma(z) \zeta^{m_0 - l n_0} \left(\zeta^{n_0} + t^d \tau(z)\right)^l - s = 0, 
$$
$s$ 
is given by 
\begin{align*}
s 
&= \sigma(z) \zeta^{m_0 - l n_0} \left(\zeta^{n_0} + t^d \tau(z)\right)^l \\
&= \sigma(z) \zeta^{m_0 - l n_0} 
\left\{\zeta^{n_0} + \left(\frac{m_0}{l n_0 - m_0} \zeta^{n_0}\right)\right\}^l \\
&= \left(\frac{l n_0}{l n_0 - m_0}\right)^l \sigma(z) \zeta^{m_0}. 
\end{align*}
Hence:

\begin{lem}  \label{det:lem:CoreSingularityCondition}

Fix 
$t \neq 0$. 
A point 
$(z, \zeta) \in N_0$ 
is a singularity of some subordinate fiber 
$X_{s, t}$ 
of 
the deformation 
$\pi_t : M_t \to \Delta_t$ 
if and only if 
the following condition is satisfied: 
$$
\begin{cases}
K(z) 
= 0, 
\quad 
\sigma(z) 
\neq 0, 
\quad 
\tau(z) 
\neq 0, \\
\zeta^{n_0} 
= \frac{l n_0 - m_0}{m_0} t^d \tau(z). 
\end{cases}
$$
In this case, 
the following holds: 
$$
s 
= 
\left(\frac{l n_0}{l n_0 - m_0}\right)^l \sigma(z) \zeta^{m_0}. 
$$
\end{lem}

We call a zero 
$\alpha$ 
of the plot function 
$K(z)$ 
an \emph{essential zero} 
if 
$\sigma(\alpha) \neq 0$ 
and 
$\tau(\alpha) \neq 0$. 
For an essential zero 
$\alpha$ 
of 
$K(z)$, 
Lemma \ref{det:lem:CoreSingularityCondition} implies that 
$(\alpha, \beta) \in N_0$ 
is a singularity of a subordinate fiber 
$X_{s, t}$ 
if and only if 
$$
\begin{cases}
\beta^{n_0} 
= \frac{l n_0 - m_0}{m_0} t^d \tau(\alpha), \\
s 
= \left(\frac{l n_0}{l n_0 - m_0}\right)^l \sigma(\alpha) \beta^{m_0}. 
\end{cases}
$$
Eliminating 
$\beta$, 
we have 
$$
s^{\bar{n}_0} 
= \left(\frac{l n_0}{l n_0 - m_0}\right)^{l\bar{n}_0} 
\left(\frac{l n_0 - m_0}{m_0}\right)^{\bar{m}_0} 
t^{d \bar{m}_0} \sigma(\alpha)^{\bar{n}_0} \tau(\alpha)^{\bar{m}_0}, 
$$
where 
$\bar{m}_0$ 
and 
$\bar{n}_0$ 
are the relatively prime integers satisfying 
$\bar{n}_0 / \bar{m}_0 = {n_0} / {m_0}$. 
This equation for 
$s$ 
has 
$\bar{n}_0$ 
solutions, 
say, 
$s_1, s_2, \dots, s_{\bar{n}_0}$. 
Observe that 
the equation 
$$
\beta^{n_0} = \frac{l n_0 - m_0}{m_0} t^d \tau(\alpha) 
$$
for 
$\beta$ 
has 
$n_0$ 
solutions, 
say 
$\beta_1, \beta_2, \dots, \beta_{n_0}$. 
Then 
$n_0 / \bar{n}_0 \, (= \gcd(m_0, n_0))$ 
points among 
$(\alpha, \beta_1), (\alpha, \beta_2), \dots, (\alpha, \beta_{n_0})$ 
lie on one of the subordinate fibers 
$X_{s_1, t}, X_{s_2, t}, \dots, X_{s_{\bar{n}_0}, t}$.

\begin{lem}  \label{det:lem:CoreSingularityNum}
Let 
$\alpha$ 
be an essential zero of 
$K(z)$. 
Then: 
\begin{description}
 \item[\rm(a)] 
 There exist 
 $\bar{n}_0$ 
 subordinate fibers 
 $X_{s_1, t}, X_{s_2, t}, \dots, X_{s_{\bar{n}_0}, t}$ 
 that have singularities with $z$-coordinate 
 $\alpha$. 
 $($In fact, 
 $s_1, s_2, \dots, s_{\bar{n}_0}$ 
 are given as the solutions of the following equation for 
 $s$: \\
 $
 s^{\bar{n}_0} 
 = \left(\frac{l n_0}{l n_0 - m_0}\right)^{l\bar{n}_0} 
 \left(\frac{l n_0 - m_0}{m_0}\right)^{\bar{m}_0} 
 t^{d \bar{m}_0} \sigma(\alpha)^{\bar{n}_0} \tau(\alpha)^{\bar{m}_0}.)
 $ 
 \item[\rm(b)] 
 Moreover 
 the number of such singularities on each of these subordinate fibers is 
 $n_0 / \bar{n}_0$. 
\end{description}
\end{lem}

Next 
we write 
$K(z) = \sigma \tau \omega$, 
where 
$\displaystyle 
\omega(z) := \frac{d \log (\sigma^{n_0} \tau^{m_0})}{dz}$. 
Here 
$\omega$ 
is a meromorphic section of the cotangent bundle 
$\Omega_{\Theta_0}^1$ 
on 
$\Theta_0$. 
Recall the assumption that 
the subbranches 
$\overline{\br}^{(j)} \, (j = 1, 2, \dots, v)$ 
are proportional 
(so 
$m_0 n_1^{(j)} - m_1^{(j)} n_0 = 0$) 
and the others are not. 
Then 
$\omega(z)$ 
is holomorphic at 
$p^{(1)}, p^{(2)}, \dots, p^{(v)}$, 
whereas 
$\omega(z)$ 
has a pole of order 
$1$ 
at 
$p^{(v + 1)}, p^{(v + 2)}, \dots, p^{(h)}$. 
On the other hand, 
$\omega(z)$ 
has a pole of order 
$1$ 
at 
$q_1, q_2, \dots, q_k$ 
(which are zeros of the core section 
$\tau$). 
Moreover 
$$
\begin{cases}
K(z) = 0, \\
\sigma(z) \neq 0, \\
\tau(z) \neq 0 
\end{cases}
\Longleftrightarrow 
\quad 
\begin{cases}
\omega(z) = 0, \\
z \notin \left\{p^{(1)}, p^{(2)}, \dots, p^{(v)}\right\}. 
\end{cases}
$$

\begin{lem}[\cite{Ta3} Lemma 21.3.5]  \label{det:lem:CoreInvFomula}
Let 
$g_0$ 
denote the genus of the core
$\Theta_0$. 
Then
$$
\sum_{K(\alpha) = 0, 
\sigma(\alpha) \neq 0, 
\tau(\alpha) \neq 0} 
\ord_{\alpha}(K(z)) 
= (h - v) + k + (2 g_0 - 2) 
- \sum_{j = 1}^{v} \ord_{p^{(j)}}(\omega). 
$$
\end{lem}

We set 
$\chi 
:= (h - v) + k + (2 g_0 - 2) - \sum_{j = 1}^{v} \ord_{p^{(j)}}(\omega)$, 
which is called the \emph{core invariant}.

\begin{cor}  \label{det:cor:CoreInvIneq}

Let 
$\kappa$ 
denote the number of essential zeros of 
$K(z)$. 
Then we have 
$$
\kappa \leq \chi,
$$
where the equality holds 
precisely when
the order of any essential zero of 
$K(z)$ 
equals 
$1$. 
\end{cor}

\begin{proof}
For any essential zero 
$\alpha$ 
of the plot function 
$K(z)$ 
we have 
$$
\ord_{\alpha}(K(z)) \geq 1, 
$$ 
thus 
$$
\sum_{K(\alpha) = 0, 
\sigma(\alpha) \neq 0, 
\tau(\alpha) \neq 0} 
\ord_{\alpha}(K(z)) 
\geq \kappa. 
$$ 
From Lemma $\ref{det:lem:CoreInvFomula}$, 
the left hand side of this inequality is equal to the core invariant 
$\chi$, 
which confirms the assertion. 

\end{proof}

Let 
$\alpha_1, \alpha_2, \dots, \alpha_\kappa$ 
be the essential zeros of $K(z)$, 
where 
$\kappa$ 
is the number of essential zeros of 
$K(z)$. 
By Lemma \ref{det:lem:CoreSingularityNum} (a), 
for each 
$\alpha_i$, 
there exist 
$\bar{n}_0$ 
subordinate fibers that have singularities with $z$-coordinate 
$\alpha_i$, 
and 
their singular values are given by
$$
s^{\bar{n}_0} 
= \left(\frac{l n_0}{l n_0 - m_0}\right)^{l\bar{n}_0} 
\left(\frac{l n_0 - m_0}{m_0}\right)^{\bar{m}_0} 
t^{d \bar{m}_0} \sigma(\alpha_i)^{\bar{n}_0} \tau(\alpha_i)^{\bar{m}_0}.
$$ 
Thus, 
if 
$\alpha_i$ 
and 
$\alpha_j$ 
satisfy 
$$
\sigma(\alpha_i)^{\bar{n}_0} \tau(\alpha_i)^{\bar{m}_0} 
= \sigma(\alpha_j)^{\bar{n}_0} \tau(\alpha_j)^{\bar{m}_0}, 
$$
then 
the the singularities with $z$-coordinate 
$\alpha_i$ 
and 
$\alpha_j$ 
lie on the \emph{same} subordinate fiber. 
We denote by 
$\bar{\kappa}$ 
the number of the distinct values of the set 
$\left\{\sigma(\alpha_i)^{\bar{n}_0} \tau(\alpha_i)^{\bar{m}_0} 
\, : \, i = 1, 2, \dots, \kappa \right\}$. 
Then 
for a fixed 
$t \neq 0$, 
the deformation 
$\pi_t : M_t \to \Delta_t$ 
has exactly 
$\bar{n}_0 \bar{\kappa}$ 
subordinate fibers that have singularities near the core. 
This result together with Lemma \ref{det:lem:CoreSingularityNum} 
and Corollary \ref{det:cor:CoreInvIneq} confirms the following:

\begin{prop}  \label{det:prop:CoreSingularityNum}

Let us consider 
the deformation 
$\pi_t : M_t \to \Delta_t$ 
of 
$\pi : M \to \Delta$ 
for a fixed 
$t \neq 0$. 
Then 
we have the following. 

\begin{description}

\item[\rm(a)] 
$\left(\begin{tabular}[c]{c}
\textrm{The number of subordinate fibers in} 
$M_t$\\
\textrm{that have singularities near} 
$\Theta_0$
\end{tabular}\right) 
\leq \bar{n}_0 \chi$. 

\noindent
Here the equality holds 
precisely when 
the order of any essential zero equals 
$1$ 
and 
$\bar{\kappa} = \kappa$. 

\item[\rm(b)] 
$\left(\begin{tabular}[c]{c}
\textrm{The number of singularities near} $\Theta_0$ \\
\textrm{on each subordinate fiber in} 
$M_t$
\end{tabular}\right) 
\leq \displaystyle \frac{n_0}{\bar{n}_0} \chi$. \\
Here the equality holds 
precisely when 
the order of any essential zero equals 
$1$ 
and 
$\bar{\kappa} = 1$. 

\end{description}

\end{prop}

\section{Constraints from the numbers of singularities}  \label{sec:Number}

In this section, 
we show two useful lemmas 
which give us the number of the subordinate fibers 
and that of their singularities. 
See Lemmas \ref{det:lem:NoProp3Br} and \ref{det:lem:PropChi0}.

Let 
$\pi : M \to \Delta$ 
be a linear degeneration of complex curves 
with a stellar singular fiber 
$X_0 = m_0 \Theta_0 + \sum_{j = 1}^h \Br^{(j)}$. 
Suppose that 
$X_0$ 
has a simple crust 
$Y = n_0 \Theta_0 + \sum_{j = 1}^h \br^{(j)}$ 
of with barking multiplicity 
$l$. 
For brevity, 
we assume that 
\emph{the subbranches 
$\overline{\br}^{(1)}, \overline{\br}^{(2)}, \dots, \overline{\br}^{(v)}$ 
are proportional 
and 
the others are not} 
(so 
$v$ 
is the number of the proportional subbranches). 
Let 
$\Psi : \mathcal{M} \to \Delta \times \Delta^{\dagger}$ 
be a barking family of 
$\pi : M \to \Delta$ 
associated with 
$Y$. 
We define the core invariant of 
$Y$ 
as 
$$
\chi 
:= (h - v) + k + (2 g_0 - 2) - \sum_{j = 1}^{v} \ord_{p^{(j)}}(\omega),
$$
where 
$g_0$ 
is the genus of the core 
$\Theta_0$ 
and
$\omega := \frac{d}{dz}\log (\sigma^{n_0} \tau^{m_0})$.

First we assume that 
$Y$ 
has \emph{no} proportional subbranches. 
Since 
$v = 0$, 
we have 
$\chi = h + k + (2 g_0 - 2)$. 
Then 
Lemma \ref{det:lem:SubordSingularity} ensures that
the subordinate fibers have singularities only near the core.

\begin{lem}  \label{det:lem:NoPropChi1}
Suppose that 
$Y$ 
has no proportional subbranch. 
Set 
$c := \gcd(m_0, n_0)$ 
and 
$\bar{n}_0 := n_0 / c$. 
If 
$\chi = 1$, 
then for a fixed 
$t \neq 0$, 
we have the following. 

\begin{description}

\item[\rm(a)] 
$\pi_t : M_t \to \Delta_t$ 
has exactly 
$\bar{n}_0$ 
subordinate fibers. 

\item[\rm(b)] 
Each subordinate fiber of 
$\pi_t : M_t \to \Delta_t$ 
has 
$c$ 
singularities. 

\item[\rm(c)] 
The number of singularities of all the subordinate fibers of 
$\pi_t : M_t \to \Delta_t$ 
is 
$n_0$. 

\end{description}

\end{lem}

\begin{proof}
First note that 
the plot function 
$K(z)$ 
has at least one essential zero. 
Otherwise, 
from Lemma \ref{det:lem:CoreSingularityCondition}, 
there would exist no singularities around the core, 
which implies that 
$\pi_t : M_t \to \Delta_t$ 
has no subordinate fibers. 
Accordingly 
$$
1 \leq \bigl(\textrm{the number of essential zeros of } K(z) \bigr). 
$$
On the other hand, 
Corollary \ref{det:cor:CoreInvIneq} states that 
$$
\bigl(\textrm{the number of essential zeros of } K(z) \bigr) \leq \chi. 
$$
From the assumption 
$\chi = 1$, 
we obtain 
$$
\bigl(\textrm{the number of essential zeros of } K(z) \bigr) = 1. 
$$
Namely 
$K(z)$ 
has exactly one zero of order 
$1$. 
By Proposition \ref{det:prop:CoreSingularityNum}, 
we have 
\begin{eqnarray*}
&&\bigl(\textrm{the number of subordinate fibers of\ } \pi_t : M_t \to \Delta_t \bigr) 
= \bar{n}_0,\\
&&\bigl(\textrm{the number of singularities on each subordinate fiber} \bigr) 
= c,
\end{eqnarray*}
confirming (a) and (b).

(c) clearly follows from (a) and (b). 
\end{proof}

In particular:

\begin{lem}  \label{det:lem:NoProp3Br}
Suppose that
{\rm (i)} 
$\Theta_0$ 
is a complex projective line, 
{\rm (ii)} 
$X_0$ 
has three branches, 
{\rm (iii)} 
the core section 
$\tau$ 
has no zero and 
{\rm (iv)} 
$Y$ 
has no proportional subbranches. 
Set 
$c:= \gcd(m_0, n_0)$ 
and 
$\bar{n}_0 := n_0 / c$. 
Then for a fixed 
$t \neq 0$, 
we have the following. 
\begin{description}
\item[\rm (a)] 
$\pi_t : M_t \to \Delta_t$ 
has exactly 
$\bar{n}_0$ 
subordinate fibers. 

\item[\rm (b)] 
Each subordinate fiber of 
$\pi_t : M_t \to \Delta_t$ 
has 
$c$ 
singularities. 

\item[\rm (c)] 
The number of singularities of all the subordinate fibers of 
$\pi_t : M_t \to \Delta_t$ 
is 
$n_0$. 
\end{description}
\end{lem}

\begin{proof}
By assumption, 
we have 
$g_0 = 0$, $h=3$, $k=0$, 
and so 
$\chi = 1$. 
Hence 
Lemma \ref{det:lem:NoPropChi1} confirms the assertion. 
\end{proof}

\begin{rem}
By Lemma \ref{ta:lem:CoreSecExistence}, 
we can restate the condition (iii) of Lemma \ref{det:lem:NoProp3Br} as 
``$r_0 = r'_0$,'' 
where 
$r_0 := \sum_{j=1}^h m_1^{(j)} / m_0$ 
and 
$r'_0 := \sum_{j=1}^h n_1^{(j)} / n_0$. 
\end{rem}

Next we assume that 
$Y$ 
has a proportional subbranch.

\begin{lem}  \label{det:lem:PropChi0}
Suppose that
{\rm (i)} 
$\Theta_0$ 
is a complex projective line, 
{\rm (ii)} 
$X_0$ 
has three branches, 
{\rm (iii)} 
the core section 
$\tau$ 
has no zero and 
{\rm (iv)} 
$Y$ 
has a proportional subbranch 
$\overline{\br}^{(1)} 
= n_0 \Delta_0 + n_1 \Theta_1 + n_2 \Theta_2 + 
\cdots + n_\lambda \Theta_\lambda$ of $\overline{\Br}^{(1)}$. 
Then 
$\overline{\br}^{(1)}$ 
is the unique proportional subbranch of 
$Y$ 
$($that is, 
$v = 1)$. 
Moreover 
for a fixed 
$t \neq 0$, 
we have the following. 

\begin{description}

\item[\rm (a)] 
$\pi_t : M_t \to \Delta_t$ 
has exactly 
$\bar{n}_\lambda$ 
subordinate fibers. 

\item[\rm (b)] 
Each subordinate fiber of 
$\pi_t : M_t \to \Delta_t$ 
has 
$c$ 
singularities. 

\item[\rm (c)] 
The number of singularities of all the subordinate fibers of 
$\pi_t : M_t \to \Delta_t$ 
is 
$n_\lambda$. 

\end{description}

\noindent 
Here 
$c:= \gcd(m_\lambda, n_\lambda)$ 
and 
$\bar{n}_\lambda := n_\lambda / c$. 
\end{lem}

\begin{proof}
By assumption,
we have 
$g_0 = 0$, $h=3$, $k=0$. 
Thus
$$
\chi 
= 1 - v - \sum_{j = 1}^{v} \ord_{p^{(j)}}(\omega), 
$$
so 
$$
\chi + \sum_{j = 1}^{v} \left(\ord_{p^{(j)}}(\omega) + 1 \right) 
= 1. 
$$
Recall that 
$\omega(z)$ 
is holomorphic at 
$p^{(j)}$ 
for 
$j 
= 1, 2, \dots, v$, 
that is, 
$\ord_{p^{(j)}}(\omega) \geq 0$. 
Noting that 
$\chi \geq 0$ 
and 
$v \geq 1$, 
we deduce that
$\chi = 0$, 
$v = 1$ 
and 
$\ord_{p^{(1)}}(\omega) = 0$. 
Hence 
$\overline{\br}^{(1)}$ 
is the unique proportional subbranch. 
Since 
$\chi = 0$, 
from Proposition \ref{det:prop:CoreSingularityNum}, 
every subordinate fiber of 
$\pi_t : M_t \to \Delta$ 
has no singularities near the core 
$\Theta_0$. 
Therefore Proposition \ref{det:prop:PropSingularityNum}
confirms (a), (b) and (c).
\end{proof}

\section{Determination of the subordinate fibers, $2$}  \label{sec:Det2}

We now determine the subordinate fibers of the remaining barking families.

We first consider barking families whose simple crust has \emph{no} proportional subbranches. 
In the barking family $\boldsymbol{[III.2]}$, 
$III$ 
is deformed to 
$I_1$: 
$$
III \barkarrow I_1. 
$$
By Lemma \ref{det:lem:CaseSubordEuler2}, 
the set of subordinate fibers 
is one of 
(i) 
$\left\{
II
\right\}$, 
(ii)
$\left\{
I_2
\right\}$, 
and 
(iii) 
$\left\{
I_1, I_1
\right\}$. 
Note that 
the simple crust for this family has no proportional subbranches. 
See $\boldsymbol{[III.2]}$ of the list in Section \ref{sec:Talist}. 
Applying Lemma \ref{det:lem:NoProp3Br}, 
since 
$c = 2$ 
and 
$\bar{n}_0 = 1$, 
we deduce that 
there appears exactly one subordinate fiber 
and it has two singularities. 
This condition is satisfied only for the case (ii). 
Hence:

\begin{prop}  \label{det:prop:Splitting3.2}
The barking family $\boldsymbol{[III.2]}$ splits the singular fiber 
$III$ 
as follows: 
$$
III \longrightarrow I_1 \, + \, I_2, 
$$
where 
$I_1$ 
is the main fiber and 
$I_2$ 
is a subordinate fiber. 
\end{prop}

Similarly:

\begin{prop}  \label{det:prop:Splitting4.2}
The barking family $\boldsymbol{[III^*.2]}$ splits the singular fiber 
$III^*$ 
as follows: 
$$
III^* \longrightarrow I_1^* \, + \, I_2, 
$$
where 
$I_1$ 
is the main fiber and 
$I_2$ 
is a subordinate fiber. 
\end{prop}

In the barking family $\boldsymbol{[IV.2]}$, 
$IV$ 
is deformed to 
$I_2$: 
$$
IV \barkarrow I_2. 
$$
By Lemma \ref{det:lem:CaseSubordEuler2}, 
the set of subordinate fibers 
is one of 
(i) 
$\left\{
II
\right\}$, 
(ii)
$\left\{
I_2
\right\}$, 
and 
(iii) 
$\left\{
I_1, I_1
\right\}$. 
Applying Lemma \ref{det:lem:NoProp3Br}, 
since 
$c = 1$ 
and 
$\bar{n}_0 = 2$, 
we deduce that 
there appear two subordinate fibers 
and each of them has one singularity. 
This condition is satisfied only for the case (iii). 
Hence:

\begin{prop}  \label{det:prop:Splitting5.2}
The barking family $\boldsymbol{[IV.2]}$ splits the singular fiber 
$IV$ 
as follows: 
$$
IV \longrightarrow I_2 \, + \, I_1 + I_1, 
$$
where 
$I_2$ 
is the main fiber and the two 
$I_1$ 
are subordinate fibers. 
\end{prop}

Similarly:

\begin{prop}  \label{det:prop:Splitting6.2}
The barking family $\boldsymbol{[IV^*.2]}$ splits the singular fiber 
$IV^*$ 
as follows: 
$$
IV^* \longrightarrow I_0^* \, + \, I_1 + I_1, 
$$
where 
$I_0^*$ 
is the main fiber and the two 
$I_1$ 
are subordinate fibers. 
\end{prop}

In the barking family $\boldsymbol{[III^*.4]}$, 
$III^*$ 
is deformed to 
$I_0^*$: 
$$
III^* \barkarrow I_0^*. 
$$
By Lemma \ref{det:lem:CaseSubordEuler3}, 
the set of subordinate fibers 
is one of 
(i) 
$\left\{
III
\right\}$, 
(ii) 
$\left\{
I_3
\right\}$, 
(iii) 
$\left\{
II, I_1
\right\}$, 
(iv) 
$\left\{
I_2, I_1
\right\}$, 
and 
(v) 
$\left\{
I_1, I_1, I_1
\right\}$. 
Applying Lemma \ref{det:lem:NoProp3Br}, 
since 
$c = 1$ 
and 
$\bar{n}_0 = 3$, 
we deduce that 
there appear three subordinate fibers 
and each of them has one singularity. 
This condition is satisfied only for the case (v). 
Hence:

\begin{prop}  \label{det:prop:Splitting4.4}
The barking family $\boldsymbol{[III^*.4]}$ splits the singular fiber 
$III^*$ 
as follows: 
$$
III^* \longrightarrow I_0^* \, + \, I_1 + I_1 + I_1, 
$$
where 
$I_0^*$ 
is the main fiber and the three 
$I_1$ 
are subordinate fibers. 
\end{prop}

Similarly:

\begin{prop}  \label{det:prop:Splitting4.5}
The barking family $\boldsymbol{[III^*.5]}$ splits the singular fiber 
$III^*$ 
as follows: 
$$
III^* \longrightarrow I_6 \, + \, I_1 + I_1 + I_1, 
$$
where 
$I_6$ 
is the main fiber and the three 
$I_1$ 
are subordinate fibers. 
\end{prop}

In the barking family $\boldsymbol{[II^*.4]}$, 
$II^*$ 
is deformed to 
$I_5$: 
$$
II^* \barkarrow I_5. 
$$
By Lemma \ref{det:lem:CaseSubordEuler5}, 
the sum of the Euler characteristics of the subordinate fibers is 
$5$. 
Applying Lemma \ref{det:lem:NoProp3Br}, 
since 
$c = 1$ 
and 
$\bar{n}_0 = 5$, 
we deduce that 
there appear five subordinate fibers 
and each of them has one singularity. 
Hence:

\begin{prop}  \label{det:prop:Splitting2.4}
The barking family $\boldsymbol{[II^*.4]}$ splits the singular fiber 
$II^*$ 
as follows: 
$$
II^* \longrightarrow I_5 \, + \, I_1 + I_1 + I_1 + I_1 + I_1, 
$$
where 
$I_5$ 
is the main fiber and the five 
$I_1$ 
are subordinate fibers. 
\end{prop}

In the following cases, 
the simple crust has a proportional subbranch. 
In the barking family $\boldsymbol{[II^*.2]}$, 
$II^*$ 
is deformed to 
$IV^*$: 
$$
II^* \barkarrow IV^*. 
$$
By Lemma \ref{det:lem:CaseSubordEuler2}, 
the set of subordinate fibers 
is one of 
(i) 
$\left\{
II
\right\}$, 
(ii)
$\left\{
I_2
\right\}$, 
and 
(iii) 
$\left\{
I_1, I_1
\right\}$. 
Note that 
the simple crust for this family 
has a proportional subbranch of length 
$2$. 
See $\boldsymbol{[II^*.2]}$ of the list in Section \ref{sec:Talist}. 
Applying Lemma \ref{det:lem:PropChi0}, 
since 
$c = 1$ 
and 
$\bar{n}_2 = 1$, 
we deduce that 
there appears exactly one subordinate fiber 
and it has one singularity. 
This condition is satisfied only for the case (i). 
Hence:

\begin{prop}  \label{det:prop:Splitting2.2}
The barking family $\boldsymbol{[II^*.2]}$ splits the singular fiber 
$II^*$ 
as follows: 
$$
II^* \longrightarrow IV^* \, + \, II, 
$$
where 
$IV^*$ 
is the main fiber and 
$II$ 
is a subordinate fiber. 
\end{prop}

In the barking family $\boldsymbol{[II^*.3]}$, 
$II^*$ 
is deformed to 
$I_2^*$: 
$$
II^* \barkarrow I_2^*. 
$$
By Lemma \ref{det:lem:CaseSubordEuler2}, 
the set of subordinate fibers 
is one of 
(i) 
$\left\{
II
\right\}$, 
(ii)
$\left\{
I_2
\right\}$, 
and 
(iii) 
$\left\{
I_1, I_1
\right\}$. 
Note that 
the simple crust for this family 
has a proportional subbranch of length 
$1$. 
See $\boldsymbol{[II^*.3]}$ of the list in Section \ref{sec:Talist}. 
Applying Lemma \ref{det:lem:PropChi0}, 
since 
$c = 1$ 
and 
$\bar{n}_1 = 2$, 
we deduce that 
there appear two subordinate fibers 
and each of them has one singularity. 
This condition is satisfied only for the case (iii). 
Hence:

\begin{prop}  \label{det:prop:Splitting2.3}
The barking family $\boldsymbol{[II^*.3]}$ splits the singular fiber 
$II^*$ 
as follows: 
$$
II^* \longrightarrow I_2^* \, + \, I_1 + I_1, 
$$
where 
$I_2^*$ 
is the main fiber and the two 
$I_1$ 
are subordinate fibers. 
\end{prop}

We summarize Propositions \ref{det:prop:CaseSubordEuler1}, 
\ref{det:prop:Splitting2.7}--\ref{det:prop:Splitting8.2}, 
\ref{det:prop:Splitting3.2}--\ref{det:prop:Splitting2.3} as follows:

\begin{thm}  \label{thm:SplittingThm}
Each barking family in Takamura's list $(\ref{eq:talist})$ 
except 
$\boldsymbol{[III^*.8]}$, 
$\boldsymbol{[IV.3]}$, 
$\boldsymbol{[IV.4]}$, 
$\boldsymbol{[I_0^*.2]}$ 
splits the singular fiber as follows: 
\begin{align*}\rm
&\boldsymbol{[II.1]} \ II \longrightarrow I_1 + I_1 &
 &\boldsymbol{[III^*.2]} \ III^* \longrightarrow I_1^* + I_2 \\
&\boldsymbol{[II.2]} \ II \longrightarrow I_1 + I_1 &
 &\boldsymbol{[III^*.3]} \ III^* \longrightarrow I_2^* + I_1 \nonumber\\
&\boldsymbol{[II^*.1]} \ II^* \longrightarrow III^* + I_1 &
 &\boldsymbol{[III^*.4]} \ III^* \longrightarrow I_0^* + I_1 + I_1 + I_1 \nonumber\\
&\boldsymbol{[II^*.2]} \ II^* \longrightarrow IV^* + II &
 &\boldsymbol{[III^*.5]} \ III^* \longrightarrow I_6 + I_1 + I_1 + I_1 \nonumber\\
&\boldsymbol{[II^*.3]} \ II^* \longrightarrow I_2^* + I_1 + I_1 &
 &\boldsymbol{[III^*.6]} \ III^* \longrightarrow I_2^* + I_1 \nonumber\\
&\boldsymbol{[II^*.4]} \ II^* \longrightarrow I_5 &
 &\boldsymbol{[III^*.7]} \ III^* \longrightarrow I_7 + I_1 + I_1 \nonumber\\
&\hspace{35pt} + I_1 + I_1 + I_1 + I_1 + I_1 &
 &\boldsymbol{[III^*.9]} \ III^* \longrightarrow IV^* + I_1 \nonumber\\
&\boldsymbol{[II^*.5]} \ II^* \longrightarrow I_3^* + I_1 &
 &\boldsymbol{[IV.1]} \ IV \longrightarrow I_3 + I_1 \nonumber\\
&\boldsymbol{[II^*.6]} \ II^* \longrightarrow I_3^* + I_1 &
 &\boldsymbol{[IV.2]} \ IV \longrightarrow I_2 + I_1 + I_1 \nonumber\\
&\boldsymbol{[II^*.7]} \ II^* \longrightarrow I_8 + I_1 + I_1 &
 &\boldsymbol{[IV^*.1]} \ IV^* \longrightarrow I_1^* + I_1 \nonumber\\
&\boldsymbol{[II^*.8]} \ II^* \longrightarrow III^* + I_1 &
 &\boldsymbol{[IV^*.2]} \ IV^* \longrightarrow I_0^* + I_1 + I_1 \nonumber\\
&\boldsymbol{[II^*.9]} \ II^* \longrightarrow III^* + I_1 &
 &\boldsymbol{[IV^*.3]} \ IV^* \longrightarrow I_6 + I_1 + I_1 \nonumber\\
&\boldsymbol{[III.1]} \ III \longrightarrow I_2 + I_1 &
 &\boldsymbol{[IV^*.4]} \ IV^* \longrightarrow I_1^* + I_1 \nonumber\\
&\boldsymbol{[III.2]} \ III \longrightarrow I_1 + I_2 &
 &\boldsymbol{[I_0^*.1]} \ I_0^* \longrightarrow I_4 + I_1 + I_1 \nonumber\\
&\boldsymbol{[III.3]} \ III \longrightarrow I_2 + I_1 &
 &\boldsymbol{[I_n^*.1]} \ I_n^* \longrightarrow I_{n-1}^* + I_1 \nonumber\\
&\boldsymbol{[III^*.1]} \ III^* \longrightarrow IV^* + I_1 &
 &\boldsymbol{[I_n^*.2]} \ I_n^* \longrightarrow I_{n+4} + I_1 + I_1. \nonumber
\end{align*}
\end{thm}

\section{Supplement: Monodromy decompositions}  \label{sec:Decomps}

In this section, 
we give decompositions of the standard monodromy matrices 
corresponding to the splittings of the singular fibers 
induced by Takamura's barking families. 
Recall that 
$SL(2, \Z)$ 
is generated by 
$$
s_0 
= \left(
\begin{array}{cc}
1 & 1\\
0 & 1
\end{array}
\right) 
\quad 
\textrm{and}
\quad
s_2 
= \left(
\begin{array}{cc}
1  & 0\\
-1 & 1
\end{array}
\right). 
$$
Note that, 
since 
$s_0 s_2 s_0 
= s_2 s_0 s_2$, 
we have 
$$
s_2 
= \left(s_0 s_2\right) s_0 \left(s_0 s_2\right)^{-1}. 
$$

\paragraph{\bf Decomposition of $A_{II}$}

\quad
The standard monodromy matrix of 
$II$ 
is 
$A_{II} 
= s_0 s_2$. 
$A_{II}$ 
is decomposed into two conjugacies of 
$A_{I_1}$ 
as follows: 
$$
A_{II} 
= s_0 s_2 = A_{I_1} \cdot \left(s_0 s_2\right) A_{I_1} \left(s_0 s_2\right)^{-1}. 
$$
In fact, 
the splitting 
$II \longrightarrow I_1 + I_1$ 
occurs in the barking families 
$\boldsymbol{[II.1]}$ 
and 
$\boldsymbol{[II.2]}$.

\paragraph{\bf Decomposition of $A_{III}$}

\quad
The standard monodromy matrix of 
$III$ 
is 
$A_{III} 
= s_0 s_2 s_0$. 
$A_{III}$ 
is decomposed into 
$A_{I_2}$ 
and a conjugacy of 
$A_{I_1}$: 
$$
A_{III} 
= s_0 s_2 s_0 
= s_0^2 \left(s_0^{-1} s_2 s_0\right) 
= A_{I_2} \cdot s_2 A_{I_1} s_2^{-1}. 
$$
In fact, 
the splitting 
$III \longrightarrow I_2 + I_1$ 
occurs in the barking families 
$\boldsymbol{[III.1]}$, 
$\boldsymbol{[III.2]}$, 
$\boldsymbol{[III.3]}$.

$A_{III}$ 
has other monodromy decompositions as follows 
(but we have not found barking families that admit the corresponding splittings): 
\begin{align*}
A_{III} 
&= \left(s_0 s_2\right) s_0 = A_{II} \cdot A_{I_1}, \\
&\quad
\bigl(III \longrightarrow II + I_1\bigr) \\
A_{III} 
&= s_0 s_2 s_0 = A_{I_1} \cdot \left(s_0 s_2\right) A_{I_1} \left(s_0 s_2\right)^{-1} \cdot A_{I_1}. \\
&\quad
\bigl(III \longrightarrow I_1 + I_1 + I_1\bigr)
\end{align*}

\paragraph{\bf Decomposition of $A_{IV}$}

\quad
The standard monodromy matrix of 
$IV$ 
is 
$A_{IV} 
= s_0 s_2 s_0 s_2$. 
$A_{IV}$ 
is decomposed into 
$A_{I_3}$ 
and a conjugacy of 
$A_{I_1}$: 
\begin{align*}
A_{IV} 
&= s_0 s_2 s_0 s_2 
= s_0^3 \left(s_0^{-1} s_2 s_0\right) \\
&= A_{I_3} \cdot s_2 A_{I_1} s_2^{-1}. 
\end{align*}
In fact, 
the splitting 
$IV \longrightarrow I_3 + I_1$ 
occurs in the barking family 
$\boldsymbol{[IV.1]}$.

$A_{IV}$ 
has another monodromy decomposition 
\begin{align*}
A_{IV} 
&= s_0 s_2 s_0 s_2 
= s_0^2 s_2 s_0 \\
&= A_{I_2} 
\cdot \left(s_0 s_2\right) A_{I_1} \left(s_0 s_2\right)^{-1} 
\cdot A_{I_1}, 
\end{align*}
while the barking family 
$\boldsymbol{[IV.2]}$ 
induces the splitting 
$IV \longrightarrow I_2 + I_1 + I_1$.

We have 
other monodromy decompositions of 
$A_{IV}$ 
as follows 
(but we have not found splitting families that admit the corresponding splittings): 
\begin{align*}
A_{IV} 
&= \left(s_0 s_2 s_0\right) s_2 
= A_{III} 
\cdot \left(s_0 s_2\right) A_{I_1} \left(s_0 s_2\right)^{-1}, \\
&\quad
\bigl(IV \longrightarrow III + I_1\bigr) \\
A_{IV} 
&= \left(s_0 s_2\right)^2 
= A_{II} 
\cdot A_{II}, \\
&\quad
\bigl(IV \longrightarrow II + II\bigr) \\
A_{IV} 
&= \left(s_0 s_2\right) s_0 s_2 
= A_{II} 
\cdot A_{I_1} 
\cdot \left(s_0 s_2\right) A_{I_1} \left(s_0 s_2\right)^{-1}, \\
&\quad
\bigl(IV \longrightarrow II + I_1 + I_1\bigr) \\
A_{IV} 
&= s_0^2 s_2 \left(s_0 s_2\right) s_2^{-1} 
= A_{I_2} 
\cdot s_2 A_{II} s_2^{-1}, \\
&\quad
\bigl(IV \longrightarrow I_2 + II\bigr) \\
A_{IV} 
&= s_0 s_2 s_0 s_2 \\
&= A_{I_1} 
\cdot \left(s_0 s_2\right) A_{I_1} \left(s_0 s_2\right)^{-1} 
\cdot A_{I_1} 
\cdot \left(s_0 s_2\right) A_{I_1} \left(s_0 s_2\right)^{-1}. \\
&\quad
\bigl(IV \longrightarrow I_1 + I_1 + I_1 + I_1\bigr)
\end{align*}

\paragraph{\bf Decomposition of $A_{II^*}$}

\quad
The standard monodromy matrix of 
$II^*$ 
is 
$A_{II^*} 
= \left(s_0 s_2\right)^5$. 
$A_{II^*}$ 
is decomposed into 
$A_{III^*}$ 
and a conjugacy of 
$A_{I_1}$: 
$$
A_{II^*} 
= \left(s_0 s_2\right)^4 s_0 s_2 
= A_{III^*} \cdot (s_0 s_2) A_{I_1} (s_0 s_2)^{-1}. 
$$
In fact, 
the splitting 
$II^* \longrightarrow III^* + I_1$ 
occurs in the barking families 
$\boldsymbol{[II^*.1]}$, 
$\boldsymbol{[II^*.8]}$, 
$\boldsymbol{[II^*.9]}$.

$A_{II^*}$ 
is also decomposed into 
$A_{I_3^*}$ 
and a conjugacy of 
$A_{I_1}$: 
\begin{align*}
A_{II^*} 
&= \left(s_0 s_2\right)^3 s_0 s_2 s_0 s_2 
= \left(s_0 s_2\right)^3 s_0^3 \left(s_0^{-1} s_2 s_0\right) \\
&= A_{I_3^*} \cdot s_2 A_{I_1} s_2^{-1}. 
\end{align*}
Note that 
the barking families 
$\boldsymbol{[II^*.5]}$ 
and 
$\boldsymbol{[II^*.6]}$ 
induce the splitting 
$II^*\longrightarrow I_3^* + I_1$.

We have 
other monodromy decompositions of 
$A_{II^*}$ 
which respectively correspond to 
the splittings induced by Takamura's barking families 
as follows: 
\begin{align*}
A_{II^*} 
&= \left(s_0 s_2\right)^4 \left(s_0 s_2\right) 
= A_{IV^*} 
\cdot A_{II}, \\
&\quad
\bigl(\boldsymbol{[II^*.2]} \ II^* \longrightarrow IV^* + II\bigr) \\
A_{II^*} 
&= \left(s_0 s_2\right)^3 s_0^2 \left(s_0^{-1} s_2 s_0\right) s_2 
= A_{I_2^*} 
\cdot s_2 A_{I_1} s_2^{-1} 
\cdot \left(s_0 s_2\right) A_{I_1} \left(s_0 s_2\right)^{-1}, \\
&\quad
\bigl(\boldsymbol{[II^*.3]} \ II^* \longrightarrow I_2^* + I_1 + I_1\bigr) \\
A_{II^*} 
&= s_0^5 \left(s_0^{-1} s_2 s_0\right) s_0 s_2 s_2 s_0 \\
&\hspace{-12pt}= A_{I_5} 
\cdot s_2 A_{I_1} s_2^{-1} 
\cdot A_{I_1} 
\cdot \left(s_0 s_2\right) A_{I_1} \left(s_0 s_2 \right)^{-1} 
\cdot \left(s_0 s_2\right) A_{I_1} \left(s_0 s_2\right)^{-1} \cdot A_{I_1}, \\
&\quad
\bigl(\boldsymbol{[II^*.4]} \ II^* \longrightarrow I_5 + I_1 + I_1 + I_1 + I_1 + I_1\bigr) \\
A_{II^*} 
&= s_0^8 \left(s_0^{-2} s_2 s_0^2\right) \left(s_0^{-1} s_2^{-2} s_0 s_2^2 s_0\right) \\
&\hspace{-12pt}= A_{I_8} 
\cdot \left(s_0^{-1} s_2\right) A_{I_1} \left(s_0^{-1} s_2\right)^{-1} 
\cdot \left(s_0^{-1} s_2^{-2}\right) A_{I_1} \left(s_0^{-1} s_2^{-2}\right)^{-1}. \\
&\quad
\bigl(\boldsymbol{[II^*.7]} \ II^* \longrightarrow I_8 + I_1 + I_1\bigr)
\end{align*}

\paragraph{\bf Decomposition of $A_{III^*}$}

\quad
The standard monodromy matrix of 
$III^*$ 
is 
$A_{III^*} 
= \left(s_0s_2\right)^4 s_0$. 
$A_{III^*}$ 
is decomposed into
$A_{IV^*}$ 
and 
$A_{I_1}$: 
$$
A_{III^*} 
= \left(s_0 s_2\right)^4 s_0 
= A_{IV^*} \cdot A_{I_1}. 
$$
In fact, 
the splitting 
$III^* \longrightarrow IV^* + I_1$ 
occurs in the barking families 
$\boldsymbol{[III^*.1]}$ 
and 
$\boldsymbol{[III^*.9]}$.

$A_{III^*}$ 
is also decomposed into 
$A_{I_2^*}$ 
and a conjugacy of 
$A_{I_1}$: 
\begin{align*}
A_{III^*} 
&= \left(s_0 s_2\right)^3 s_0 s_2 s_0 
= \left(s_0 s_2\right)^3 s_0^2 \left(s_0^{-1} s_2 s_0\right) \\
&= A_{I_2^*} \cdot s_2 A_{I_1} s_2^{-1}. 
\end{align*}
Note that 
the barking families 
$\boldsymbol{[III^*.3]}$ 
and 
$\boldsymbol{[III^*.6]}$ 
induce 
the splitting 
$III^* \longrightarrow I_2^* + I_1$.

We have 
other monodromy decompositions of 
$A_{III^*}$ 
which respectively correspond to 
the splittings induced by Takamura's barking families 
as follows: 
\begin{align*}
A_{III^*} 
&= s_2^{-1} \left(s_0 s_2\right)^3 s_0 s_2 s_0^2 
= s_2^{-1} A_{I_1^*} s_2 
\cdot A_{I_2}, \\
&\quad
\bigl(\boldsymbol{[III^*.2]} \ III^* \longrightarrow I_1^* + I_2\bigr) \\
A_{III^*} 
&= \left(s_0 s_2\right)^3 s_0 s_2 s_0 
= A_{I_0^*} 
\cdot A_{I_1} 
\cdot \left(s_0 s_2\right) A_{I_1} \left(s_0 s_2\right)^{-1} 
\cdot A_{I_1}, \\
&\quad
\bigl(\boldsymbol{[III^*.4]} \ III^* \longrightarrow I_0^* + I_1 + I_1 + I_1\bigr) \\
A_{III^*} 
&= s_0^6 \left(s_0^{-3} s_2 s_0^3\right) \left(s_0^{-1} s_2 s_0\right) \left(s_0^{-1} s_2 s_0\right) \\
&= A_{I_6} 
\cdot \left(s_0^{-2} s_2\right) A_{I_1} \left(s_0^{-2} s_2\right)^{-1} 
\cdot s_2 A_{I_1} s_2^{-1} 
\cdot s_2 A_{I_1} s_2^{-1}, \\
&\quad
\bigl(\boldsymbol{[III^*.5]} \ III^* \longrightarrow I_6 + I_1 + I_1 + I_1\bigr) \\
A_{III^*} 
&= s_0^7 \left(s_0^{-5} s_2 s_0^5\right) \left(s_0^{-2} s_2 s_0^2\right) \\
&= A_{I_7} 
\cdot \left(s_0^{-4} s_2\right) A_{I_1} \left(s_0^{-4} s_2\right)^{-1} 
\cdot \left(s_0^{-1} s_2\right) A_{I_1} \left(s_0^{-1} s_2\right)^{-1}. \\
&\quad
\bigl(\boldsymbol{[III^*.7]} \ III^* \longrightarrow I_7 + I_1 + I_1\bigr)
\end{align*}

\paragraph{\bf Decomposition of $A_{IV^*}$}

\quad
The standard monodromy matrix of 
$IV^*$ 
is 
$A_{IV^*} 
= \left(s_0 s_2\right)^4$. 
$A_{IV^*}$ 
is decomposed into 
$A_{I_1^*}$ 
and a conjugacy of 
$A_{I_1}$: 
$$
A_{IV^*} 
= \left(s_0 s_2\right)^3 s_0 s_2 
= A_{I_1^*} \cdot \left(s_0 s_2\right) A_{I_1} \left(s_0 s_2\right)^{-1}. 
$$
In fact, 
the splitting 
$IV^*\longrightarrow I_1^* + I_1$
occurs in the barking families 
$\boldsymbol{[IV^*.1]}$ 
and 
$\boldsymbol{[IV^*.4]}$.

We have 
other monodromy decompositions of 
$A_{IV^*}$ 
which respectively correspond to 
the splittings induced by Takamura's barking families 
as follows: 
\begin{align*}
A_{IV^*} 
&= \left(s_0 s_2\right)^3 s_0 s_2 
= A_{I_0^*} 
\cdot A_{I_1} 
\cdot \left(s_0 s_2\right) A_{I_1} \left(s_0 s_2\right)^{-1}, \\
&\quad
\bigl(\boldsymbol{[IV^*.2]} \ IV^* \longrightarrow I_0^* + I_1 + I_1\bigr) \\
A_{IV^*} 
&= s_0^6 \left(s_0^{-4} s_2 s_0^4\right) \left(s_0^{-1} s_2 s_0\right) \\
&= A_{I_6} 
\cdot \left(s_0^{-3} s_2\right) A_{I_1} \left(s_0^{-3} s_2\right)^{-1} 
\cdot s_2 A_{I_1} s_2^{-1}. \\
&\quad
\bigl(\boldsymbol{[IV^*.3]} \ IV^* \longrightarrow I_6 + I_1 + I_1\bigr)
\end{align*}

\paragraph{\bf Decomposition of $A_{I_n^*} \, (n \geq 0)$}

\quad
The standard monodromy matrix of 
$I_n^* \, (n \geq 0)$ 
is 
$A_{I_n^*} 
= \left(s_0 s_2\right)^3 s_0^n$. 
$A_{I_n^*}$ 
is decomposed into 
$A_{I_{n+4}}$ 
and 
two conjugacies of 
$A_{I_1}$ 
as follows: 
\begin{align*}
A_{I_{n+4}} 
&= s_0 s_2 s_0 s_2 s_0 s_2 s_0^n
= s_2 \left(s_0^2 s_2 s_0^{-2}\right) s_0^4 s_0^n \\
&= \left(s_0 s_2\right) A_{I_1} \left(s_0 s_2\right)^{-1} 
\cdot \left(s_0^{3} s_2\right) A_{I_1} \left(s_0^{3} s_2\right)^{-1} 
\cdot A_{I_{n + 4}}. 
\end{align*}
In fact, 
the splitting 
$IV^* \longrightarrow I_1^* + I_1$ 
occurs in the barking families 
$\boldsymbol{[I_0^*.1]}$ 
and 
$\boldsymbol{[I_n^*.2]}$.

For 
$n \geq 1$, 
note that 
the barking family 
$\boldsymbol{[I_n^*.1]}$ 
induces the splitting 
$I_n^* \longrightarrow I_{n - 1}^* + I_1$. 
Then 
$A_{I_n^*}$ 
is also decomposed into 
$A_{I_{n - 1}^*}$ 
and 
$A_{I_1}$ as follows: 
$$
A_{I_n^*} 
= \left(s_0 s_2\right)^3 s_0^{n - 1} s_0
= A_{I_{n - 1}^*} \cdot A_{I_1}. 
$$

\section{Appendix: Takamura's list for genus $g = 1$}  \label{sec:Talist}

In \cite{Ta3}, 
for genera up to $5$,
Takamura made a list of barking families 
--- precisely speaking, 
a list of simple crusts 
(and weighted crustal sets) 
for constructing barking families --- 
which enables him to show that 
a degeneration is absolutely atomic 
if and only if 
its singular fiber is either a Lefschetz fiber 
or a multiple of a smooth curve. 
Recall that 
in a barking family, 
for a fixed 
$t \neq 0$, 
the singular fiber 
$X_{0, t}$ 
over the origin is called the main fiber 
and other singular fibers 
$X_{s, t} \, (s \neq 0)$ 
are called subordinate fibers. 
As we saw in Section \ref{sec:Ta}, 
the main fibers of barking families 
are explicitly described. 
In this paper,
when the original singular fiber 
$X_0$ 
is deformed to the main fiber 
$X_{0, t}$, 
we express 
$X_0 \barkarrow\ X_{0, t}$.

For the convenience of the reader, 
we provide Takamura's list of barking families for genus 
$1$ 
with figures of the singular fibers:

\vspace{2.0em}

$\boldsymbol{[II.1]} \ II \barkarrow I_1$
\begin{center}\includegraphics[width=5.5cm,clip]{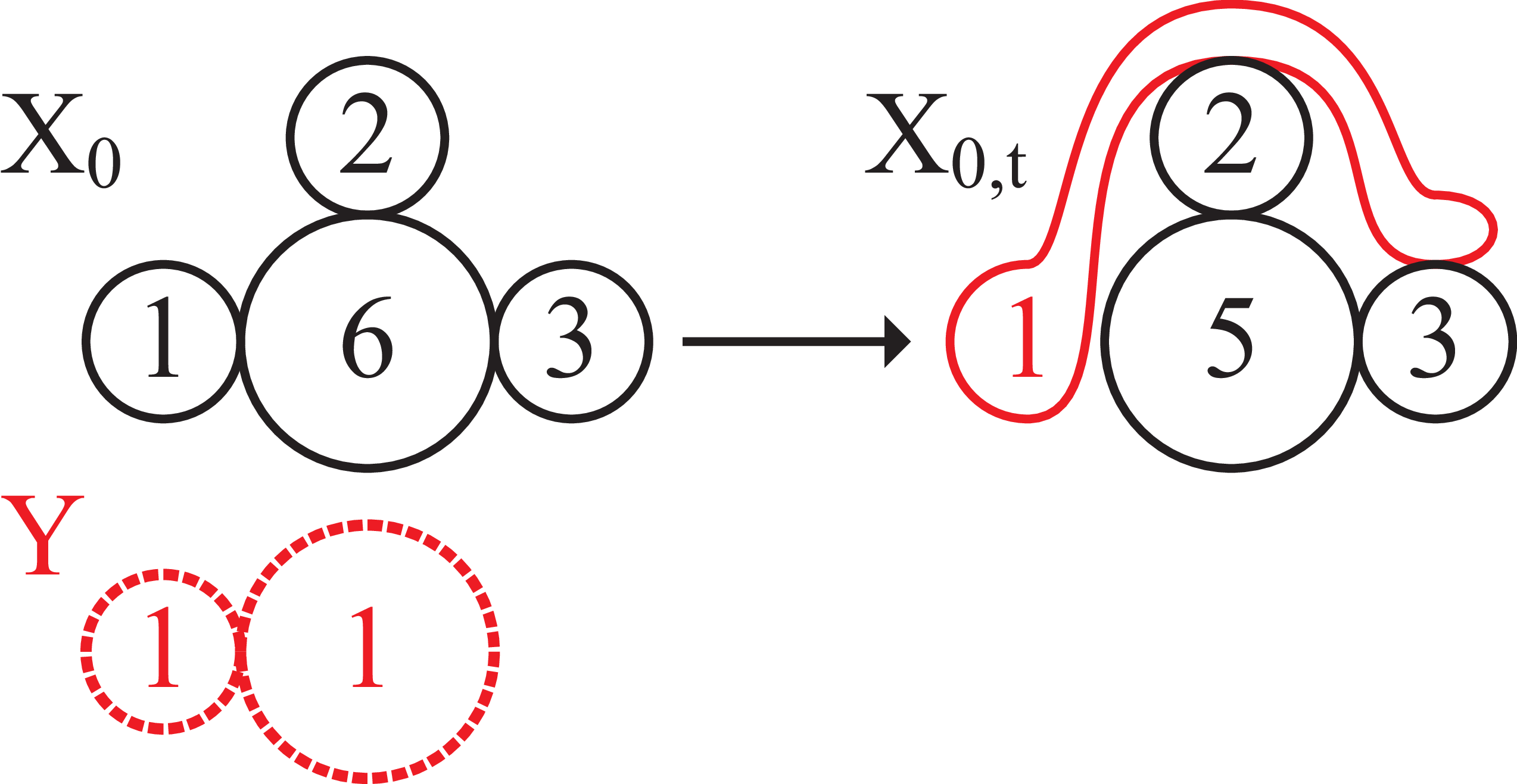}\end{center}

$\boldsymbol{[II.2]} \ II \barkarrow I_1$
\begin{center}\includegraphics[width=5.5cm,clip]{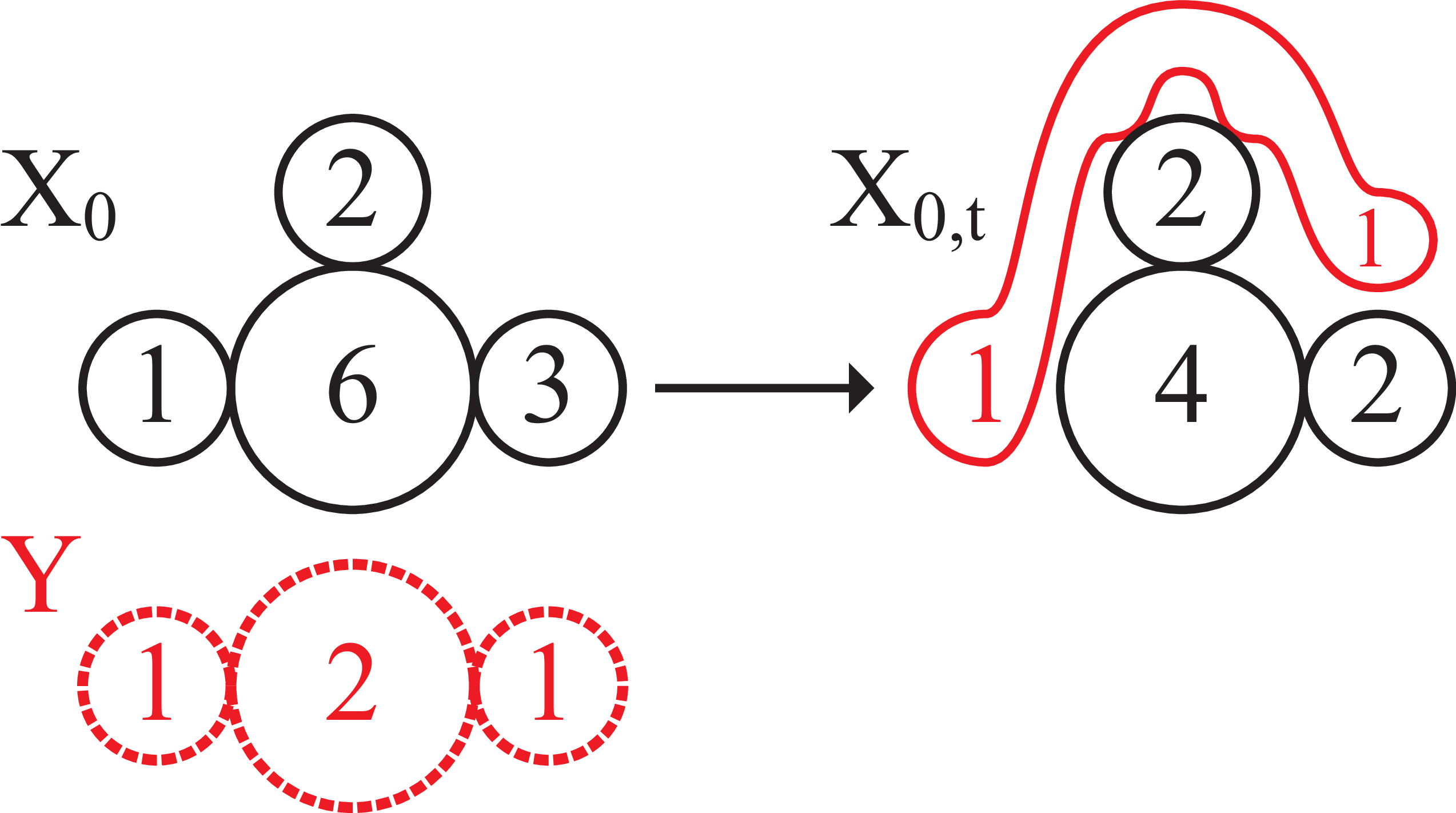}\end{center}

\vspace{2.0em}

$\boldsymbol{[II^*.1]} \ II^* \barkarrow III^*$
\begin{center}\includegraphics[width=10cm,clip]{II+1.eps}\end{center}
\newpage

$\boldsymbol{[II^*.2]} \ II^* \barkarrow IV^*$
\begin{center}\includegraphics[width=10cm,clip]{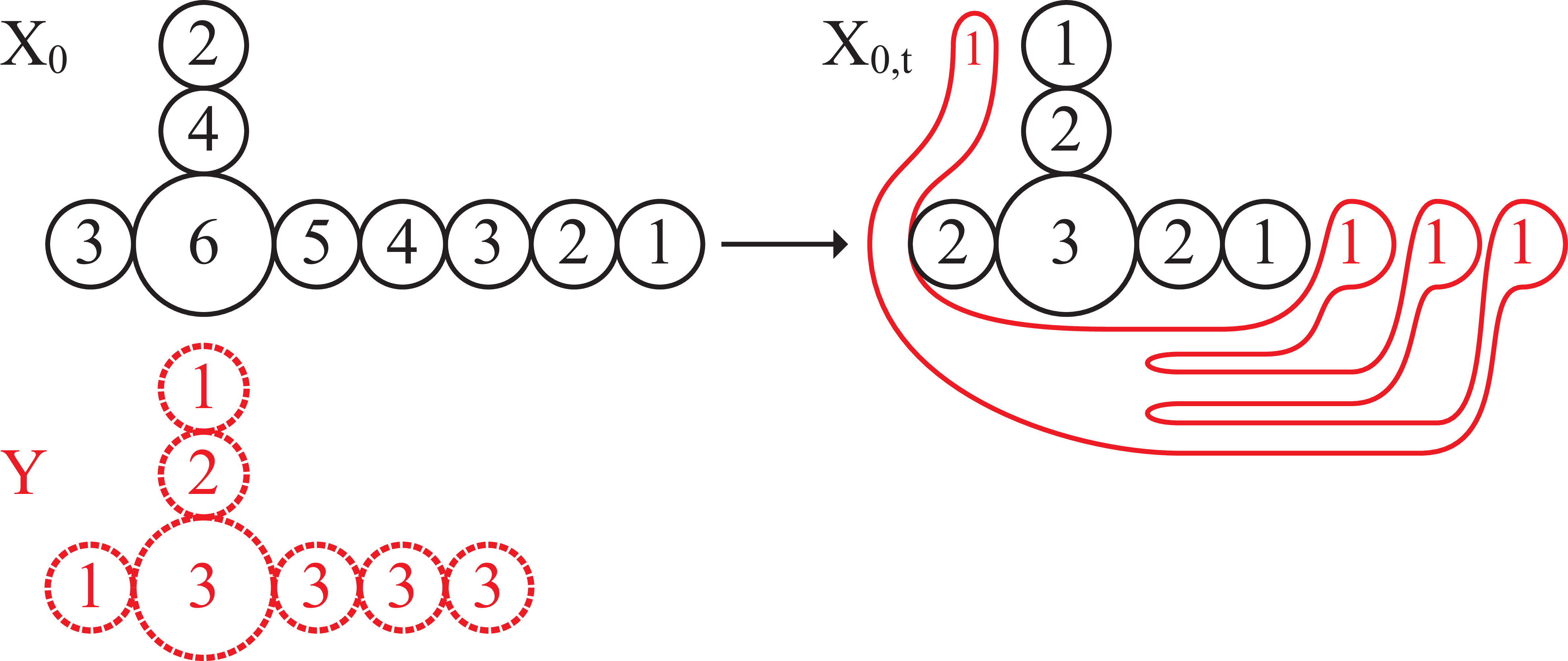}\end{center}

$\boldsymbol{[II^*.3]} \ II^* \barkarrow I_2^*$
\begin{center}\includegraphics[width=10cm]{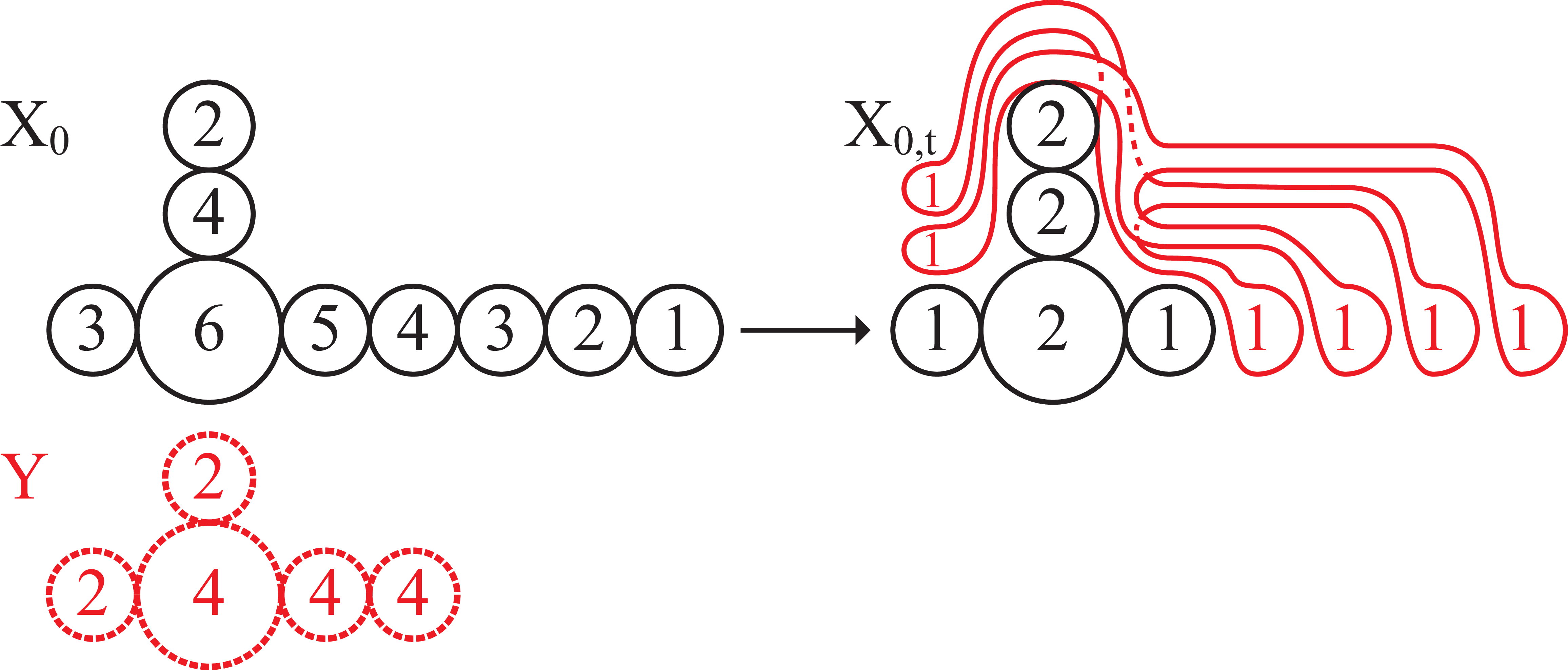}\end{center}

$\boldsymbol{[II^*.4]} \ II^* \barkarrow I_5$
\begin{center}\includegraphics[width=10cm,clip]{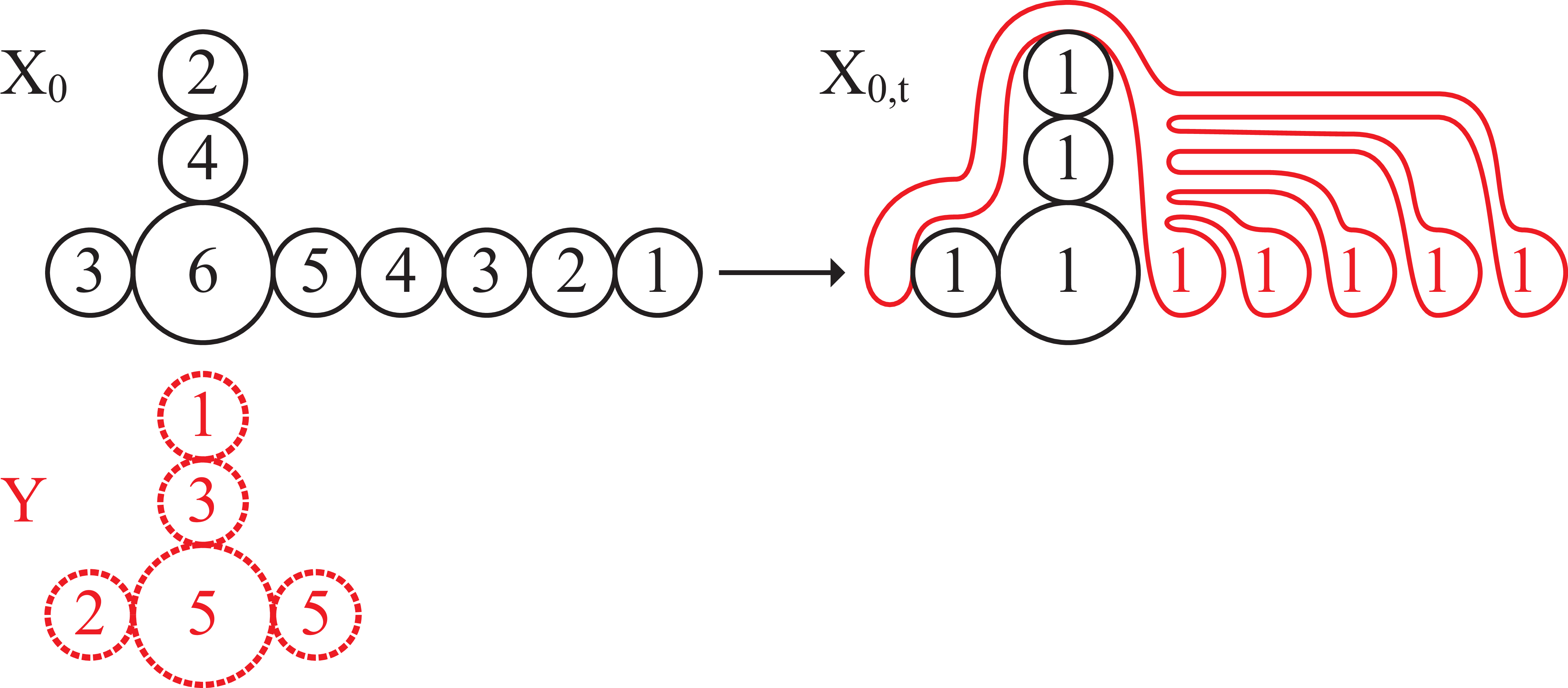}\end{center}
\newpage

$\boldsymbol{[II^*.5]} \ II^* \barkarrow I_3^*$
\begin{center}\includegraphics[width=10cm,clip]{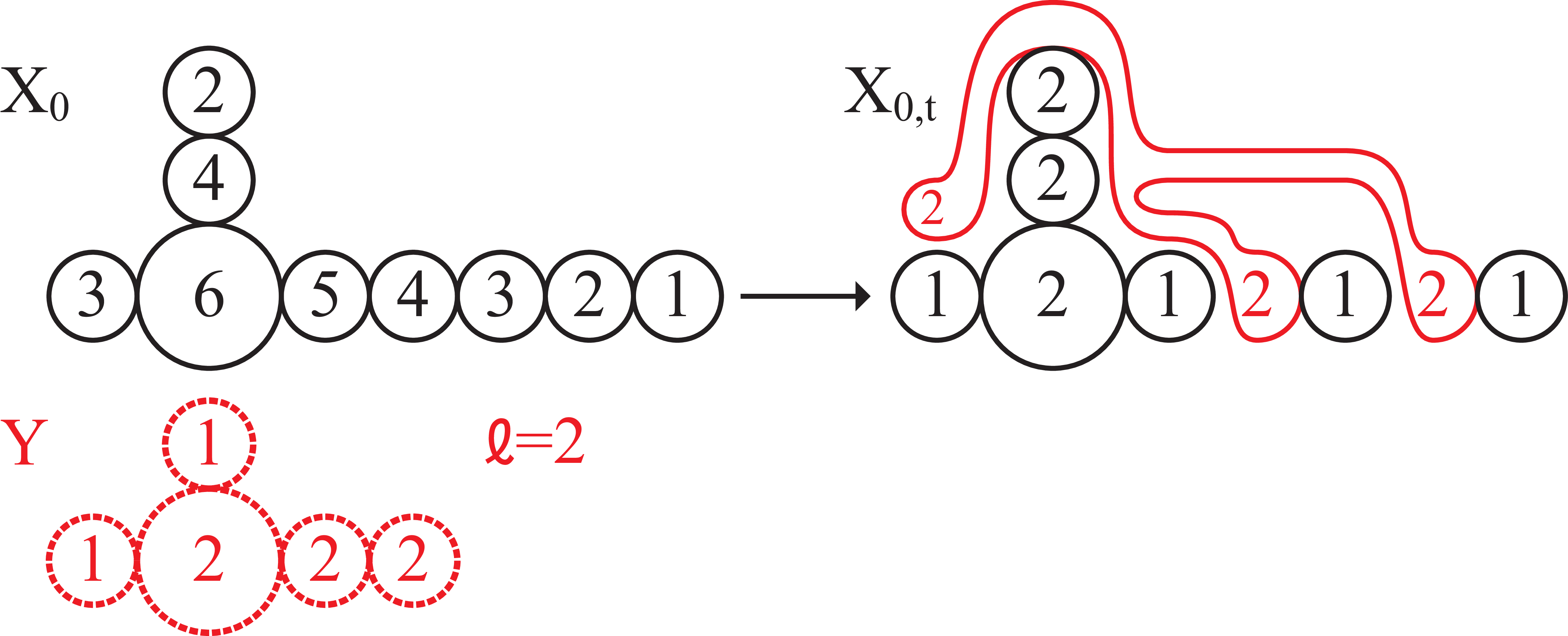}\end{center}

$\boldsymbol{[II^*.6]} \ II^* \barkarrow I_3^*$
\begin{center}\includegraphics[width=10cm,clip]{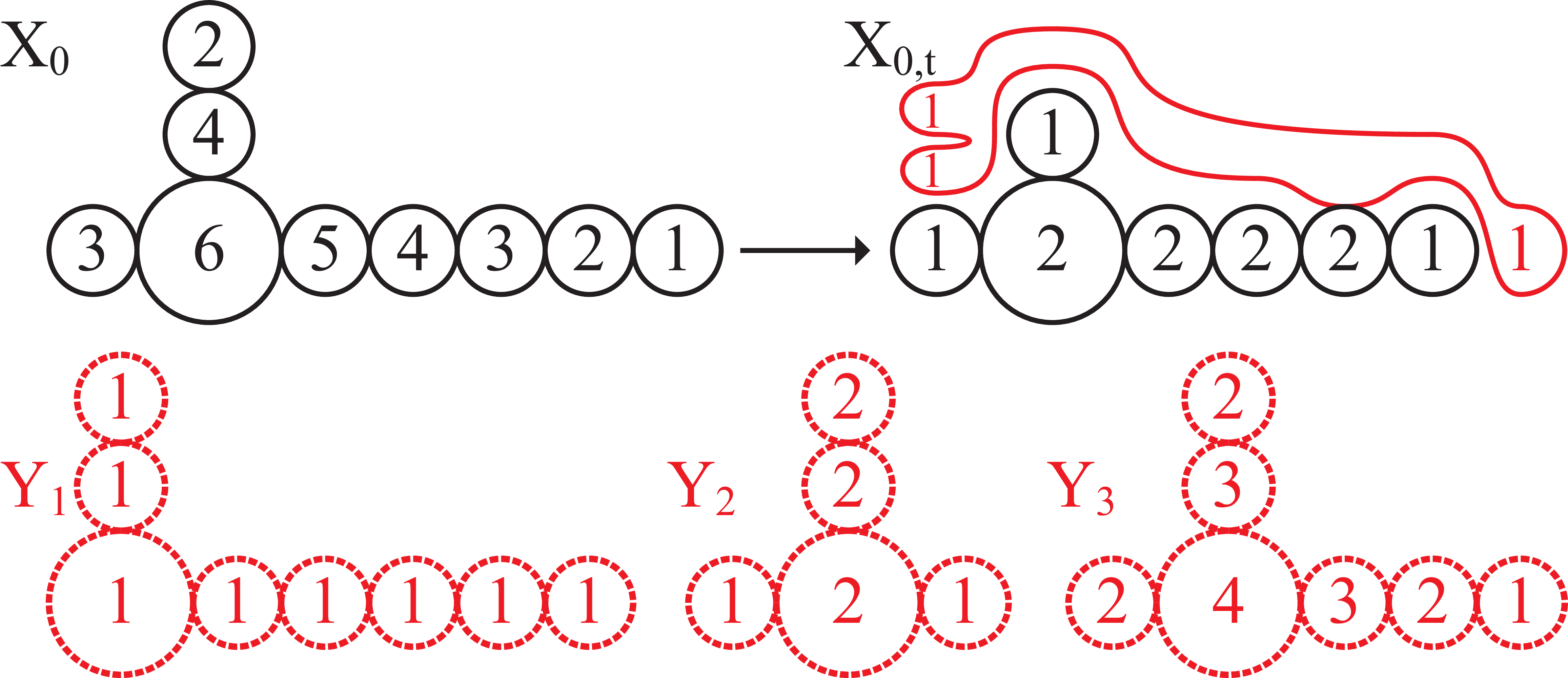}\end{center}

$\boldsymbol{[II^*.7]} \ II^* \barkarrow I_8$
\begin{center}\includegraphics[width=10cm,clip]{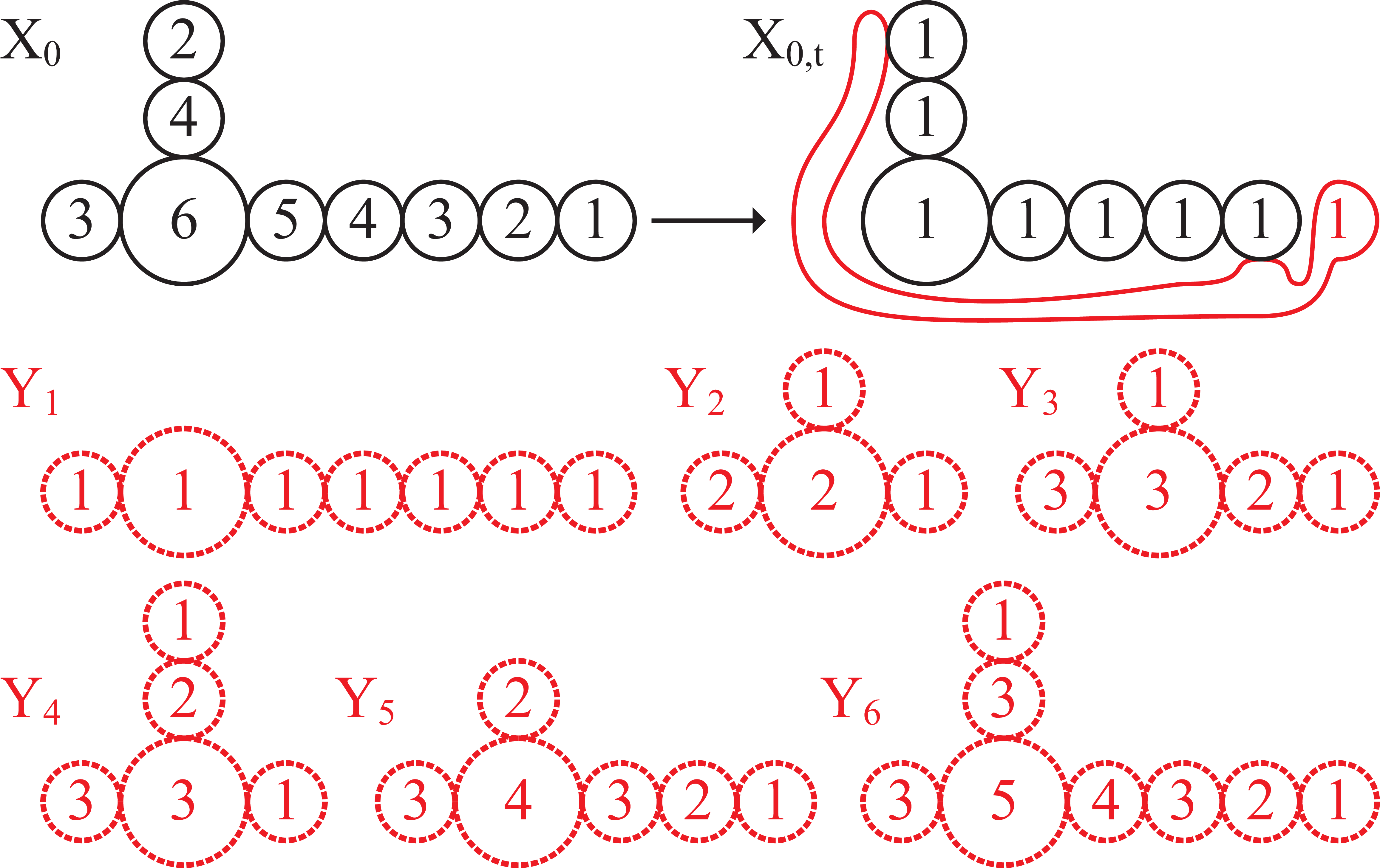}\end{center}
\newpage

$\boldsymbol{[II^*.8]} \ II^* \barkarrow III^*$
\begin{center}\includegraphics[width=10cm,clip]{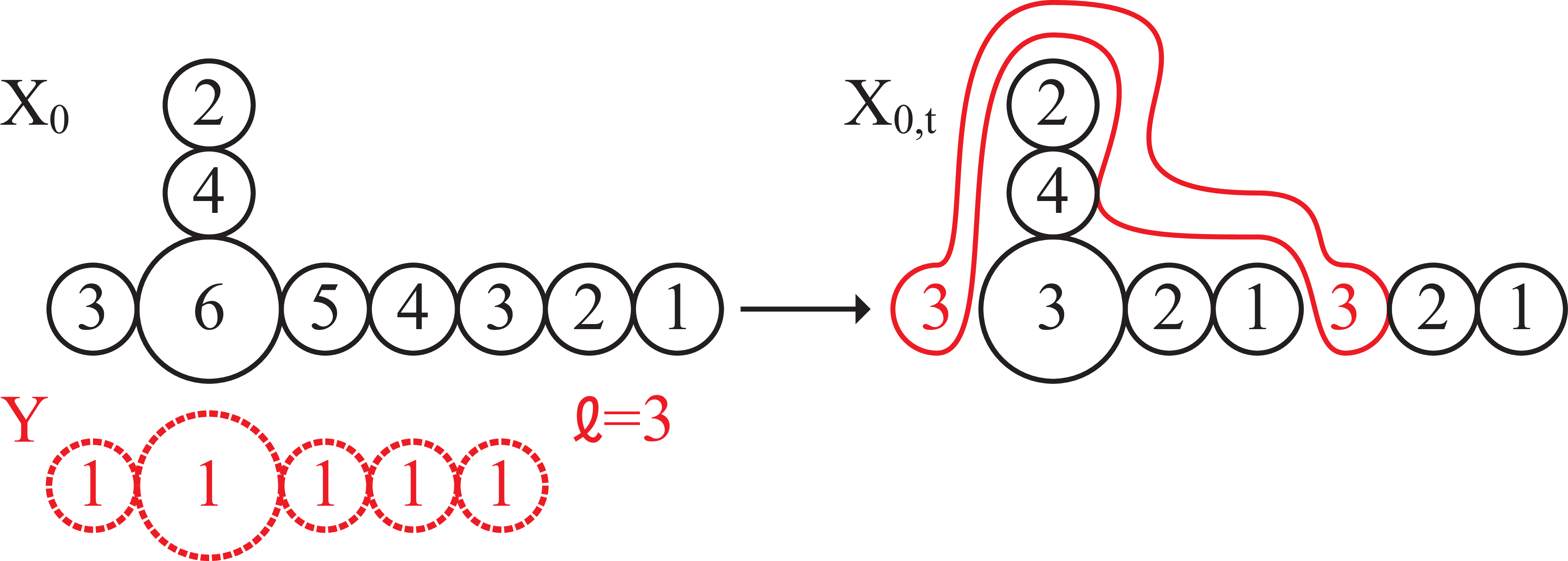}\end{center}

$\boldsymbol{[II^*.9]} \ II^* \barkarrow III^*$
\begin{center}\includegraphics[width=10cm,clip]{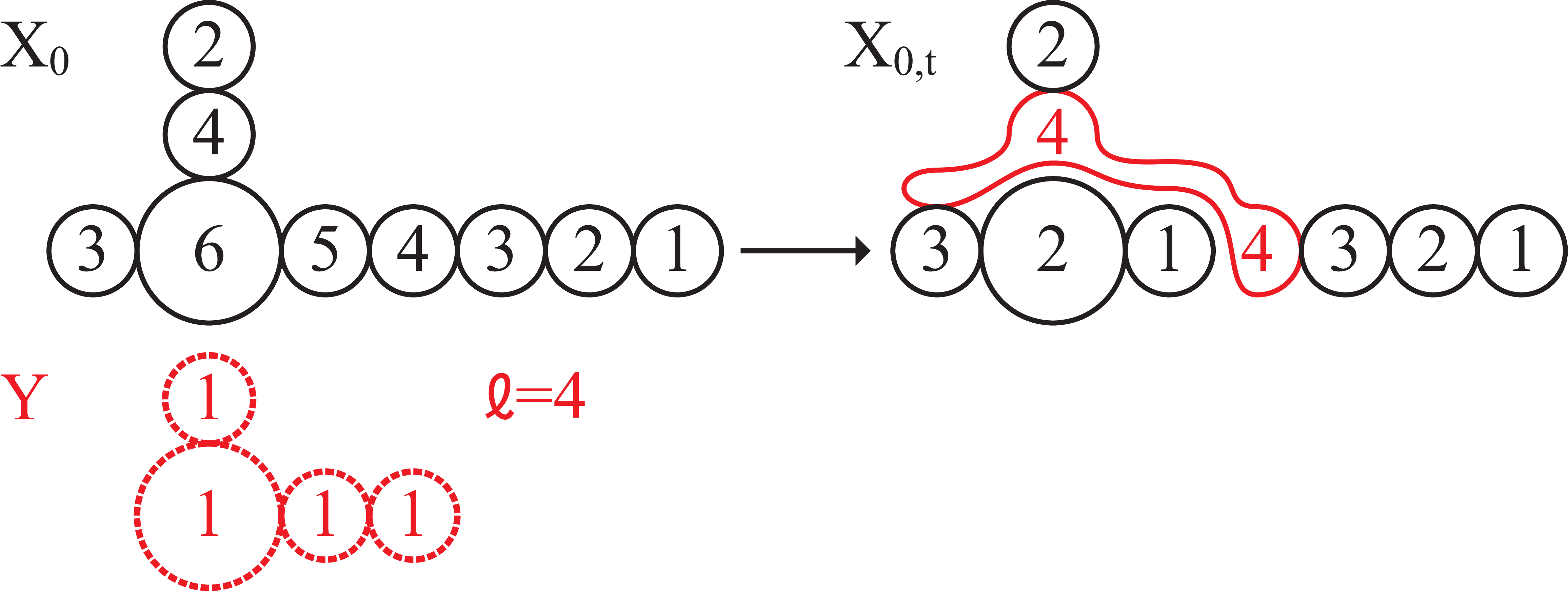}\end{center}

\vspace{2.0em}

$\boldsymbol{[III.1]} \ III \barkarrow I_2$
\begin{center}\includegraphics[width=5.5cm,clip]{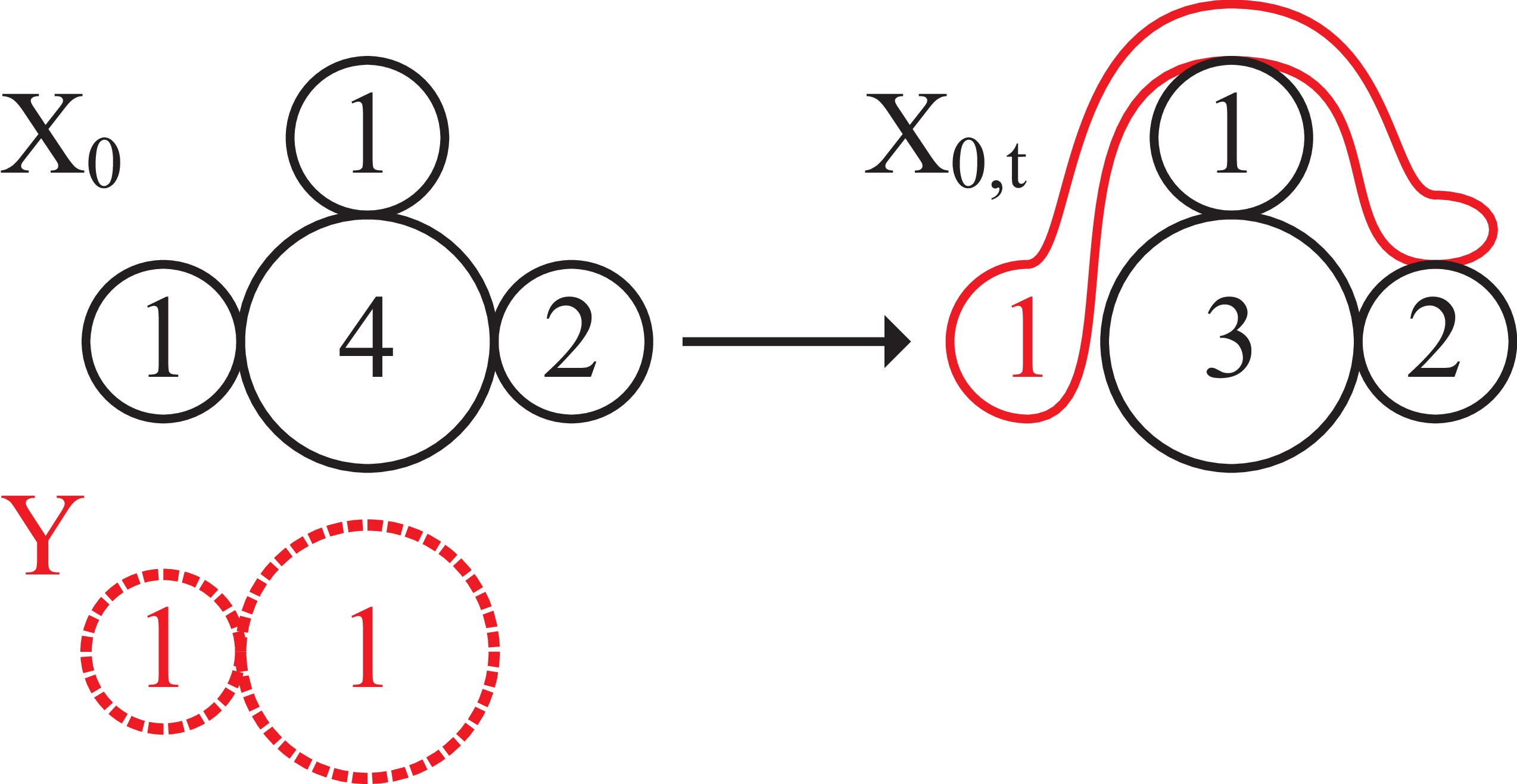}\end{center}

$\boldsymbol{[III.2]} \ III \barkarrow I_1$
\begin{center}\includegraphics[width=5.5cm,clip]{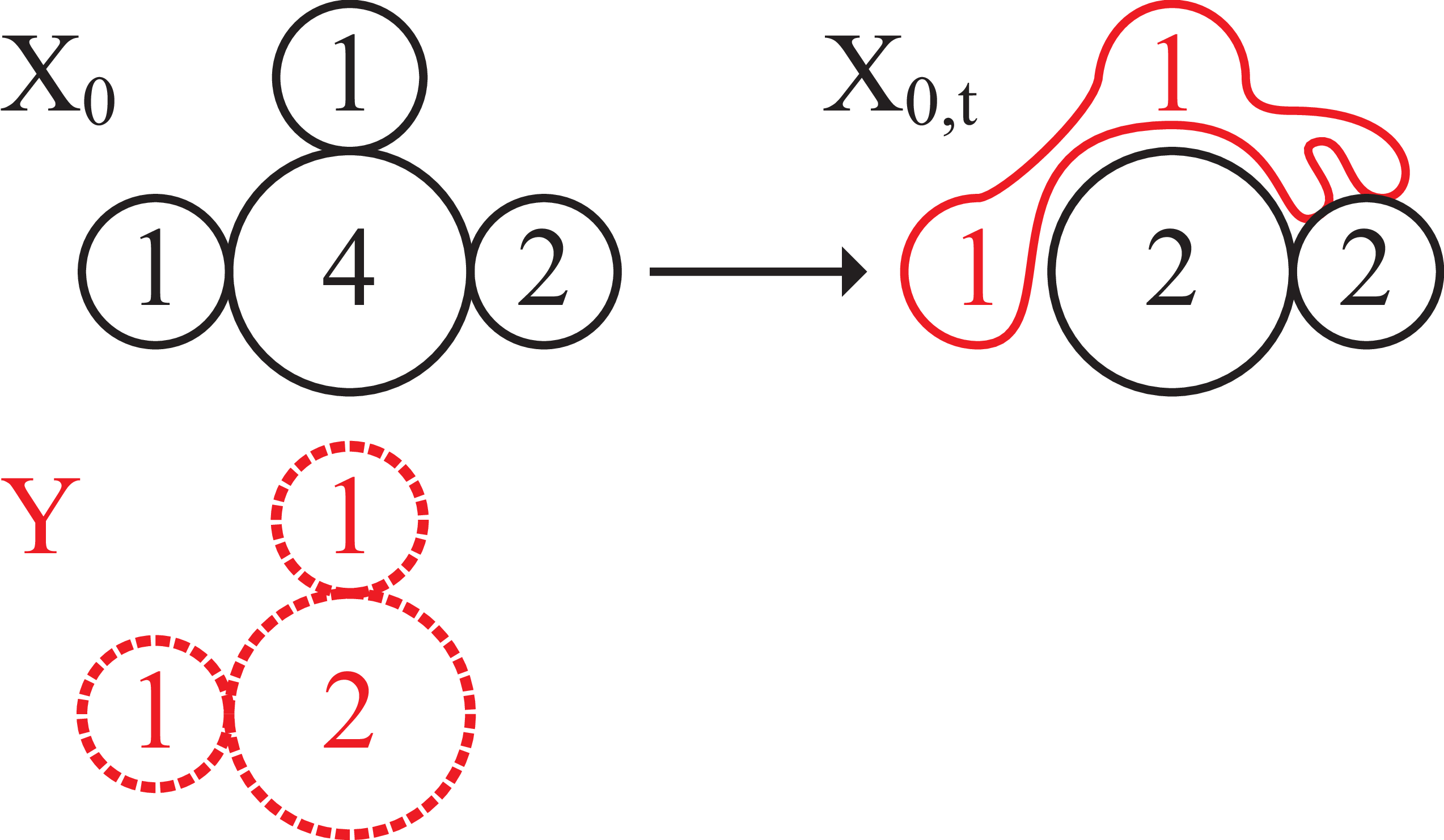}\end{center}
\newpage

$\boldsymbol{[III.3]} \ III \barkarrow I_2$
\begin{center}\includegraphics[width=5.5cm,clip]{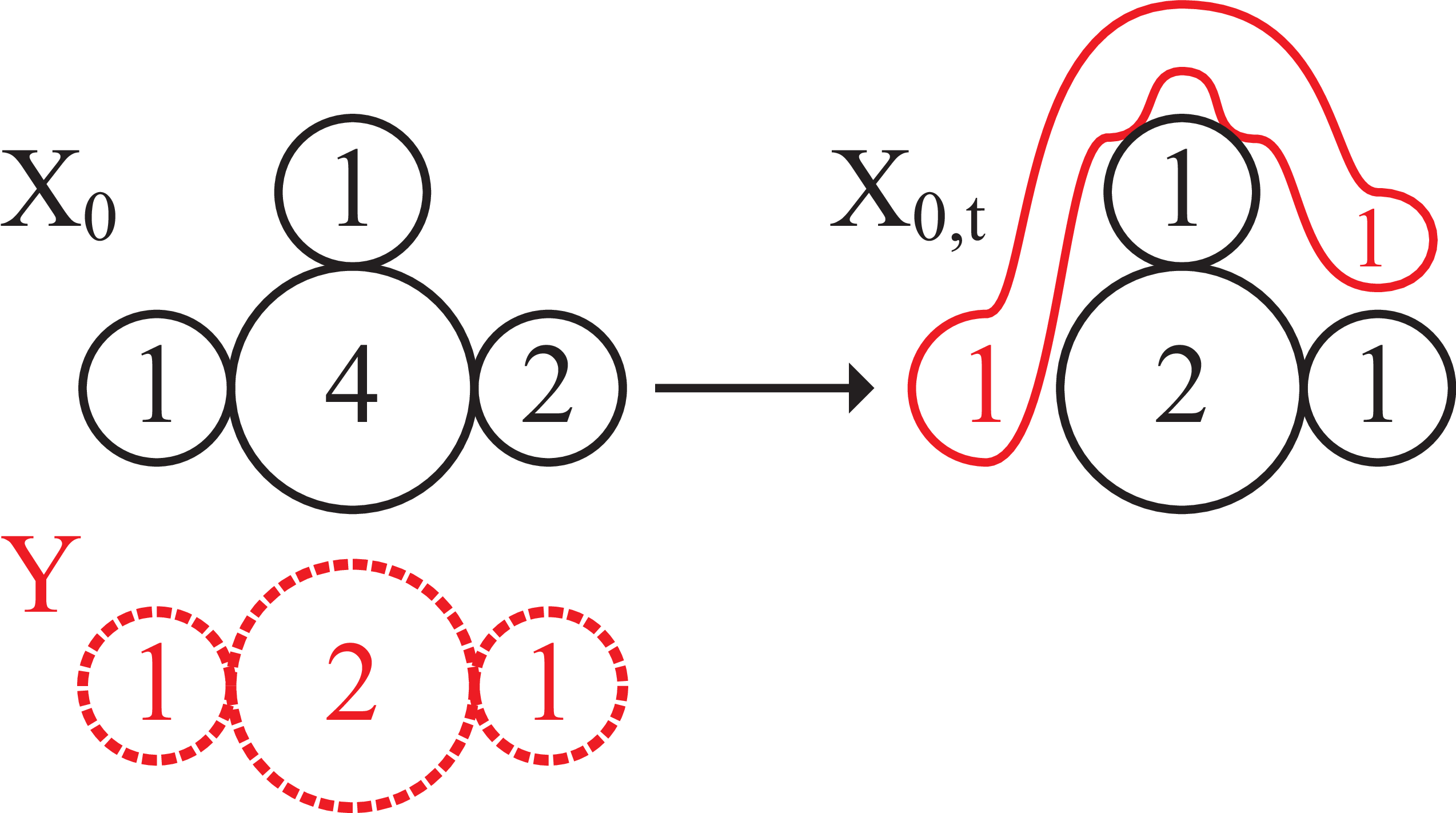}\end{center}

\vspace{2.0em}

$\boldsymbol{[III^*.1]} \ III^* \barkarrow IV^*$
\begin{center}\includegraphics[width=10cm,clip]{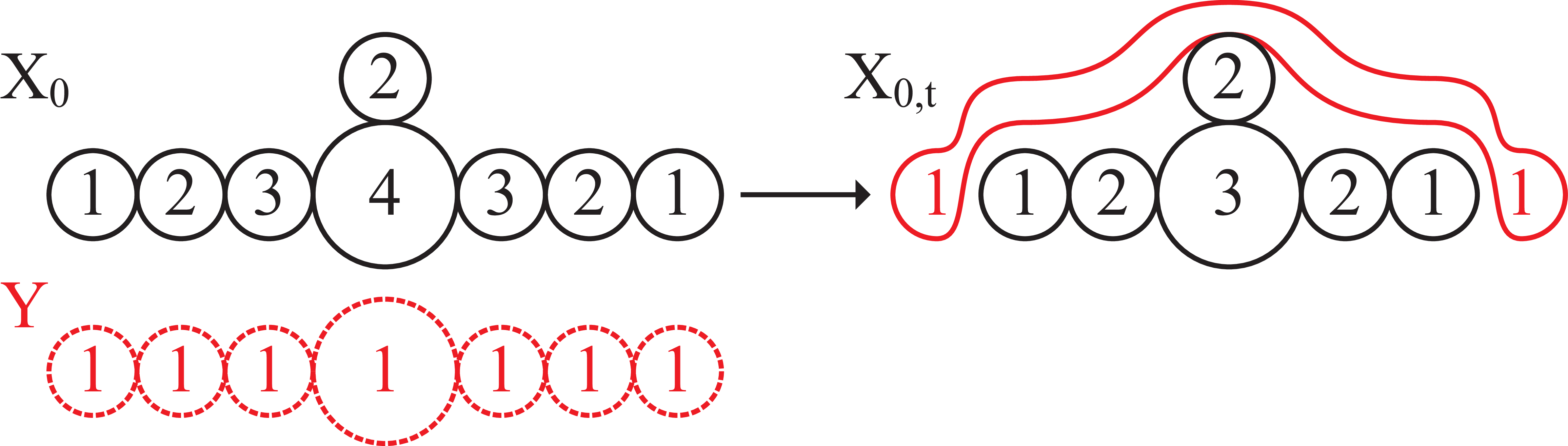}\end{center}

$\boldsymbol{[III^*.2]} \ III^* \barkarrow I_1^*$
\begin{center}\includegraphics[width=10cm,clip]{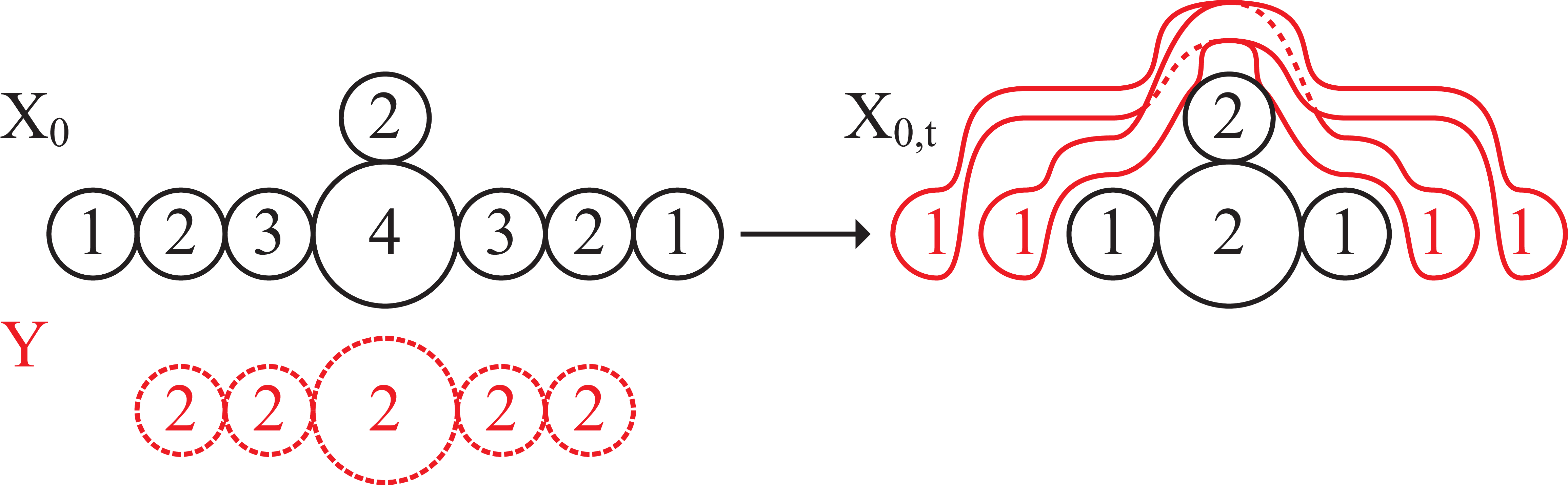}\end{center}

$\boldsymbol{[III^*.3]} \ III^* \barkarrow I_2^*$
\begin{center}\includegraphics[width=10cm,clip]{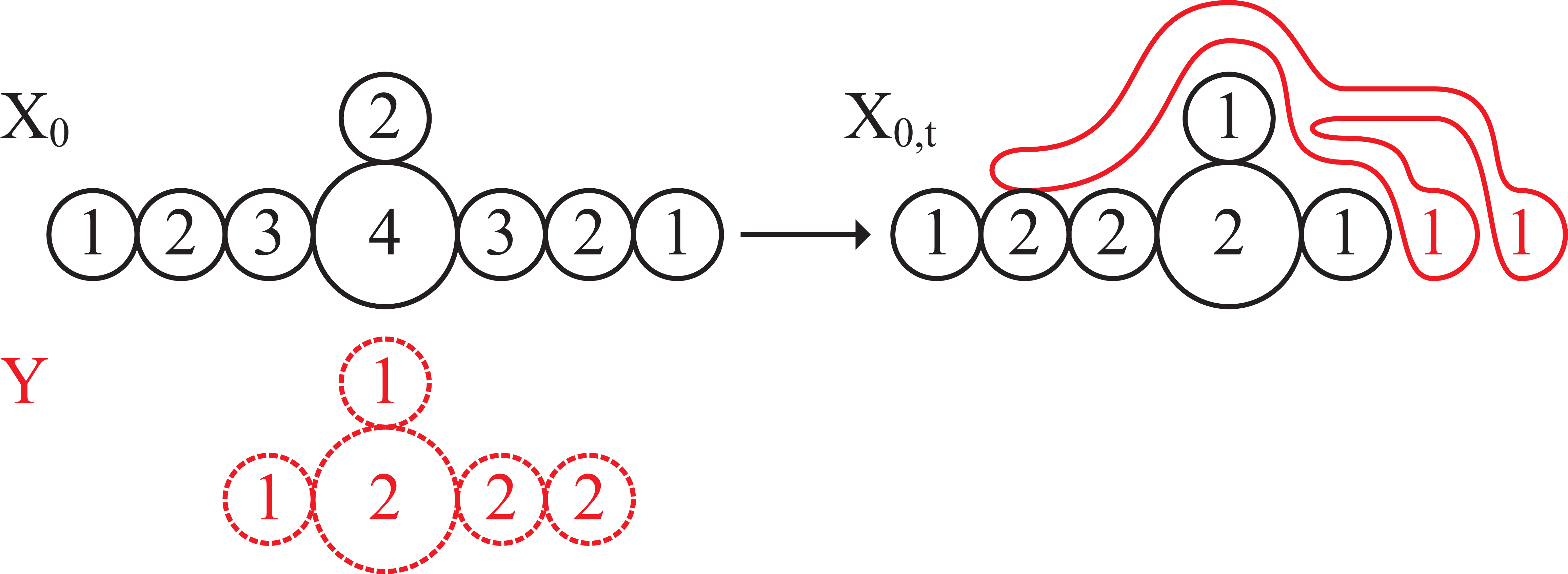}\end{center}
\newpage

$\boldsymbol{[III^*.4]} \ III^* \barkarrow I_0^*$
\begin{center}\includegraphics[width=10cm,clip]{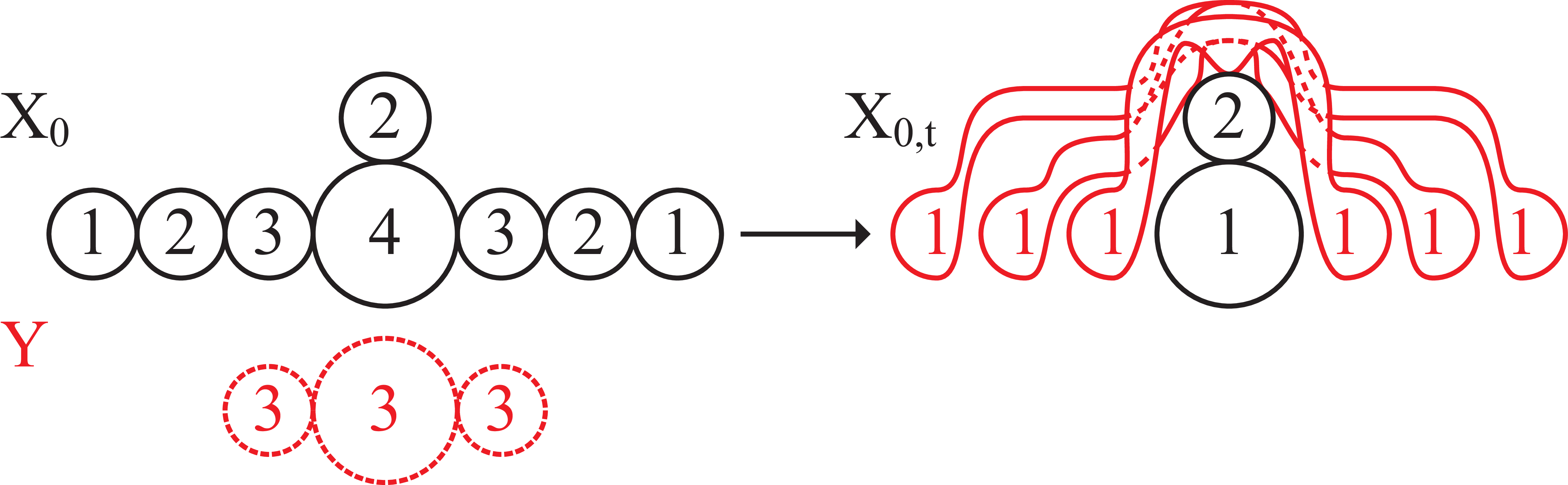}\end{center}

$\boldsymbol{[III^*.5]} \ III^* \barkarrow I_6$
\begin{center}\includegraphics[width=10cm,clip]{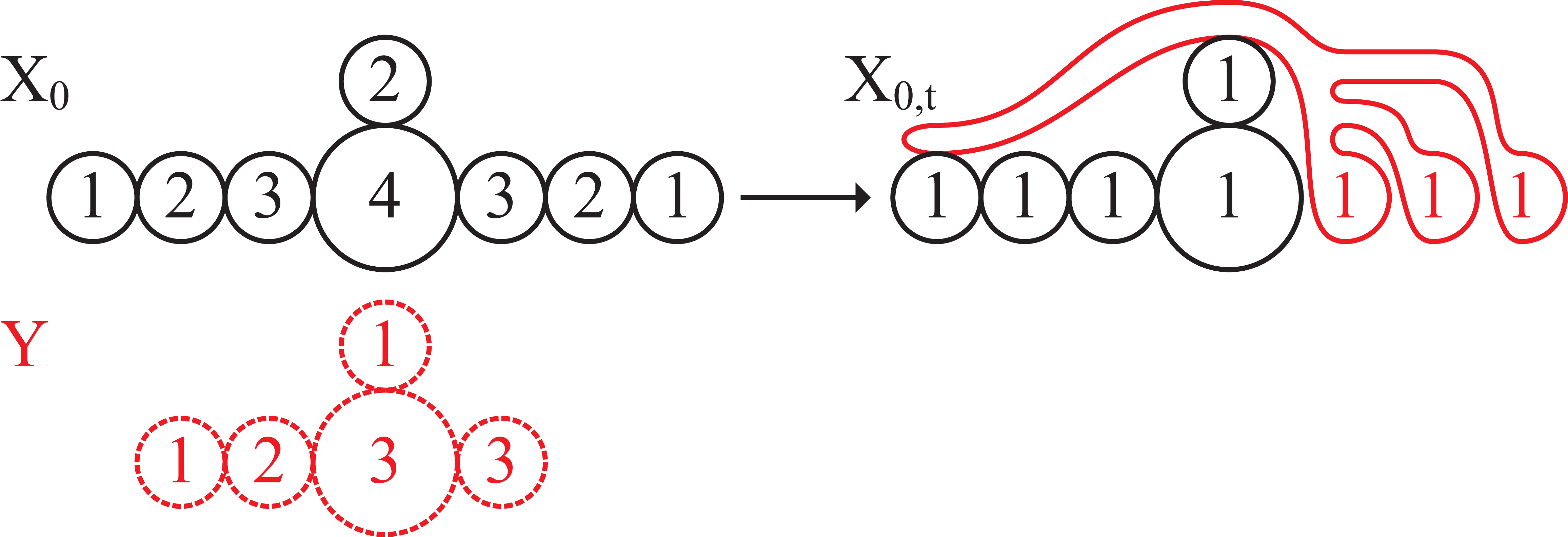}\end{center}

$\boldsymbol{[III^*.6]} \ III^* \barkarrow I_2^*$
\begin{center}\includegraphics[width=10cm,clip]{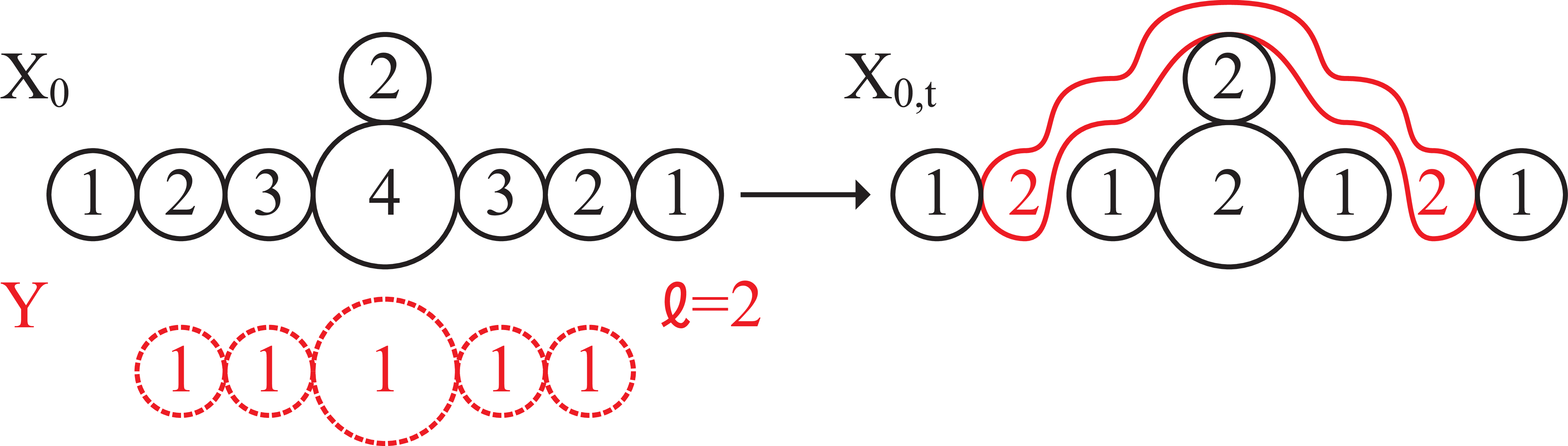}\end{center}

$\boldsymbol{[III^*.7]} \ III^* \barkarrow I_7$
\begin{center}\includegraphics[width=10cm,clip]{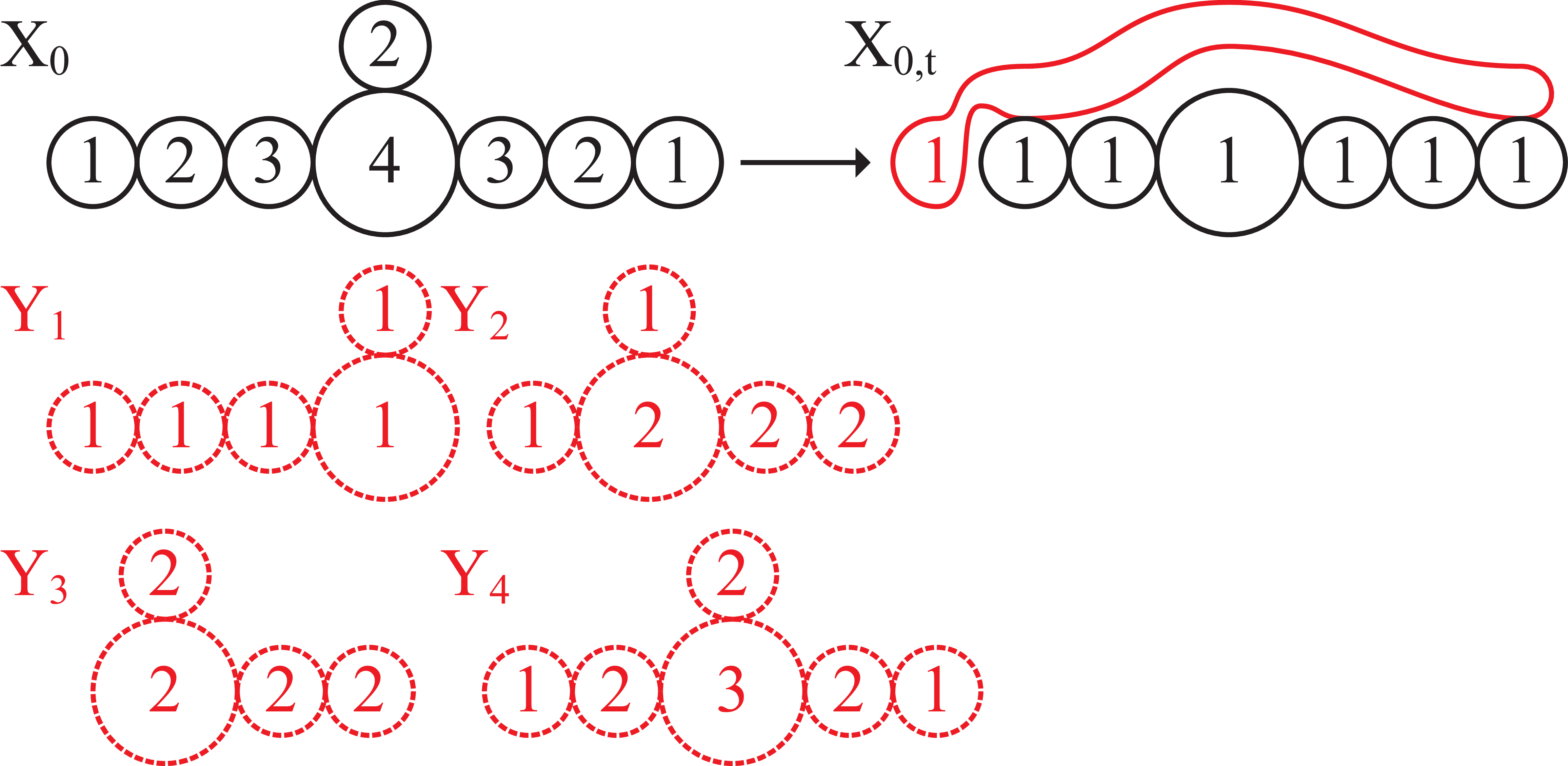}\end{center}
\newpage

$\boldsymbol{[III^*.8]} \ III^* \barkarrow I_6$
\begin{center}\includegraphics[width=10cm,clip]{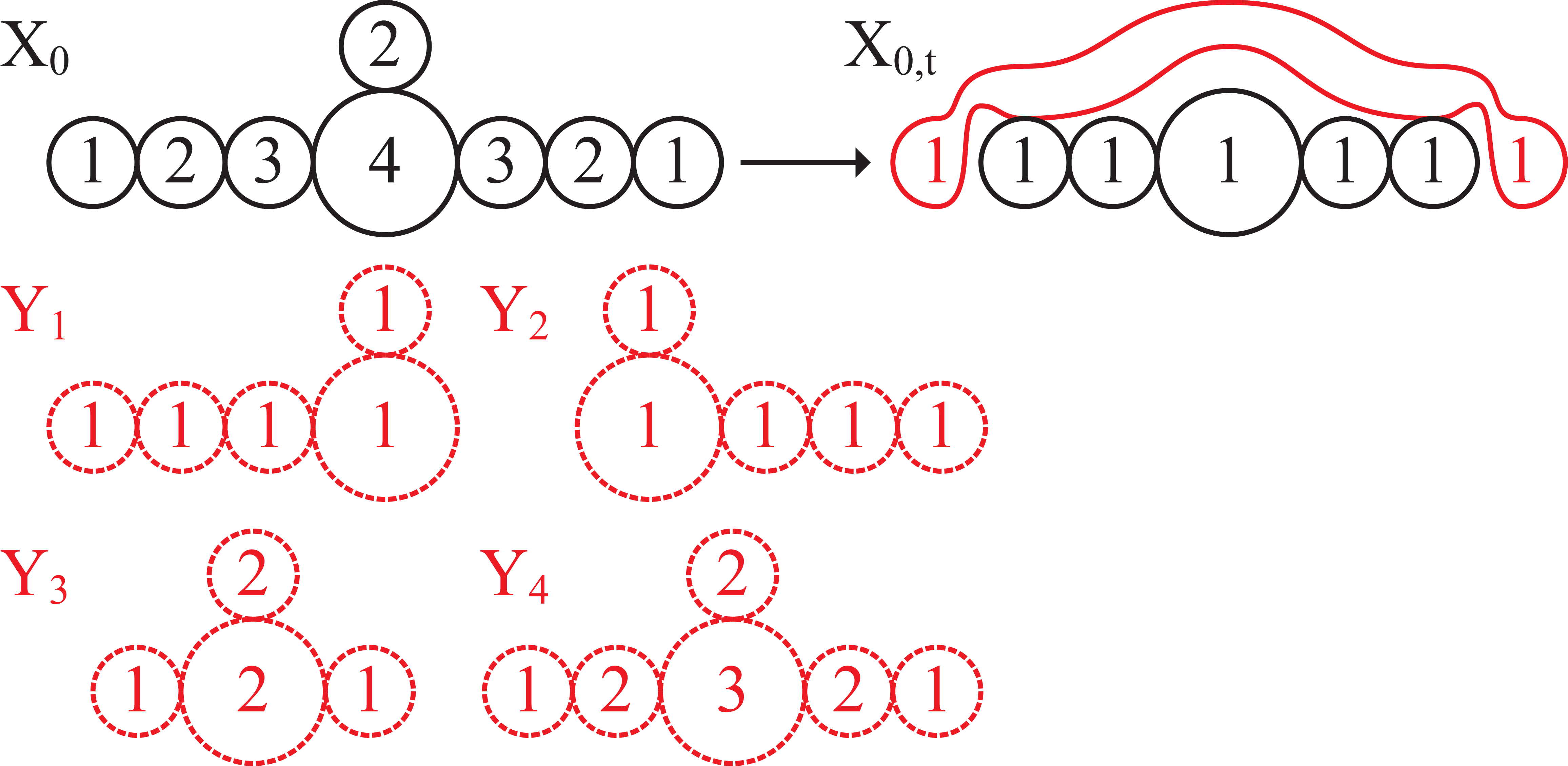}\end{center}

$\boldsymbol{[III^*.9]} \ III^* \barkarrow IV^*$
\begin{center}\includegraphics[width=10cm,clip]{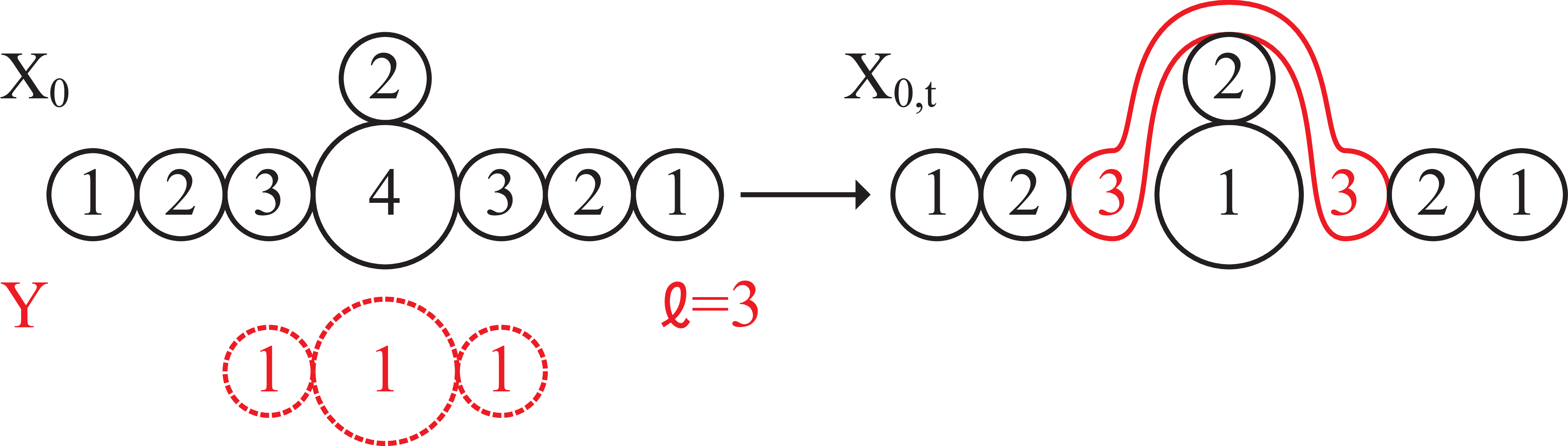}\end{center}

\vspace{2.0em}

$\boldsymbol{[IV.1]} \ IV \barkarrow I_3$
\begin{center}\includegraphics[width=5.5cm,clip]{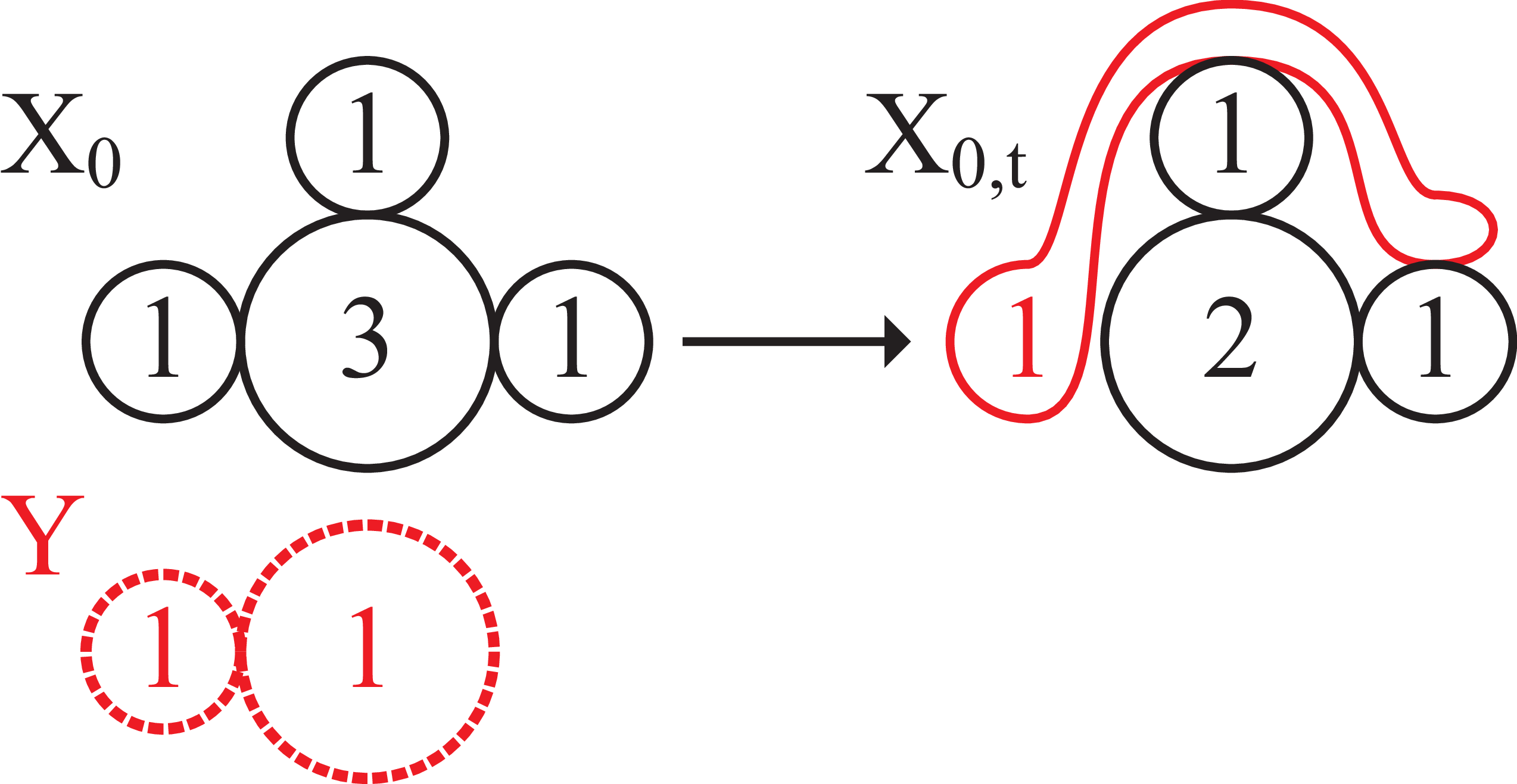}\end{center}

$\boldsymbol{[IV.2]} \ IV \barkarrow I_2$
\begin{center}\includegraphics[width=5.5cm,clip]{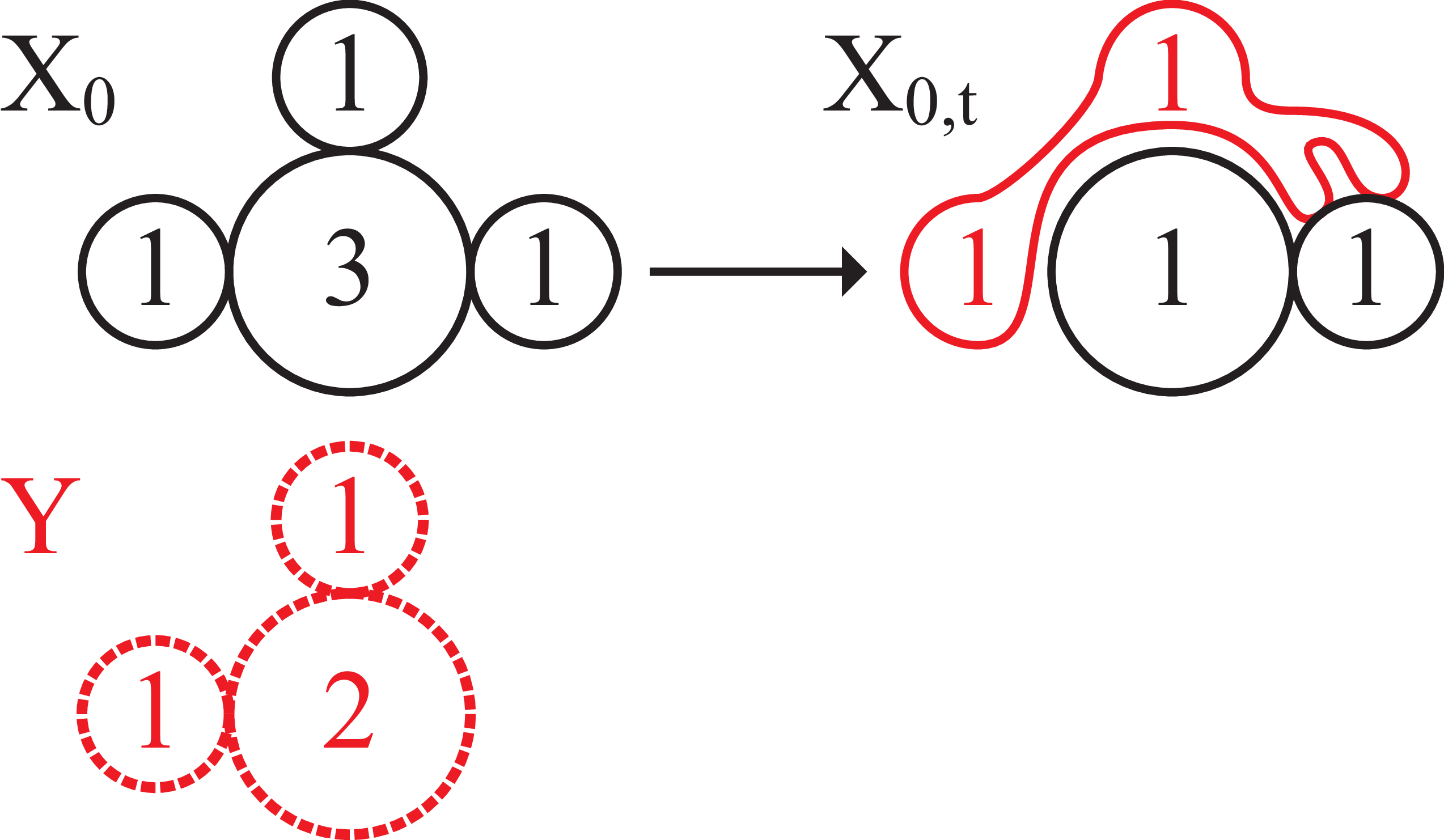}\end{center}
\newpage

$\boldsymbol{[IV.3]} \ IV \barkarrow I_2$
\begin{center}\includegraphics[width=5.5cm,clip]{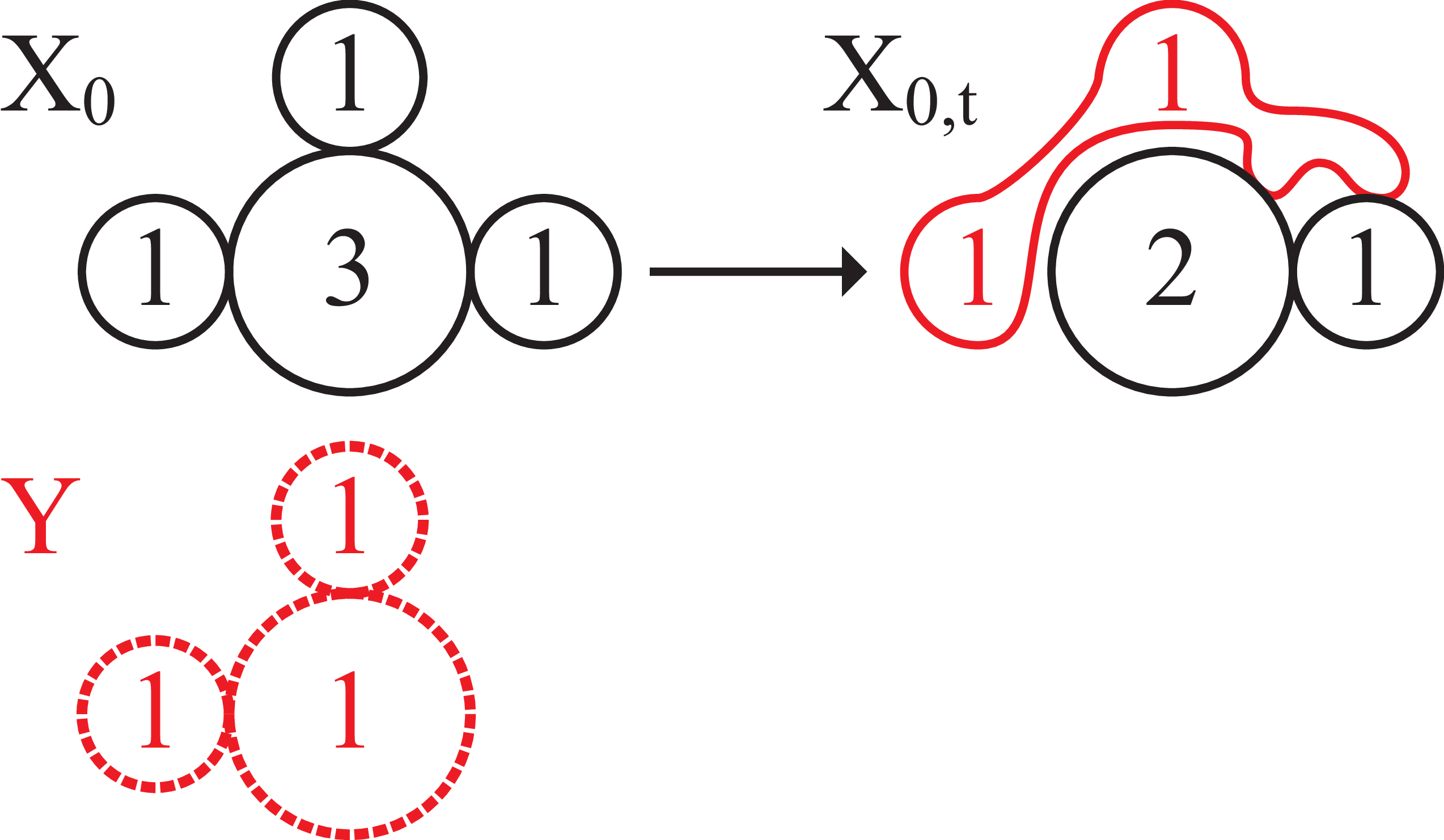}\end{center}

$\boldsymbol{[IV.4]} \ IV \barkarrow II$
\begin{center}\includegraphics[width=5.5cm,clip]{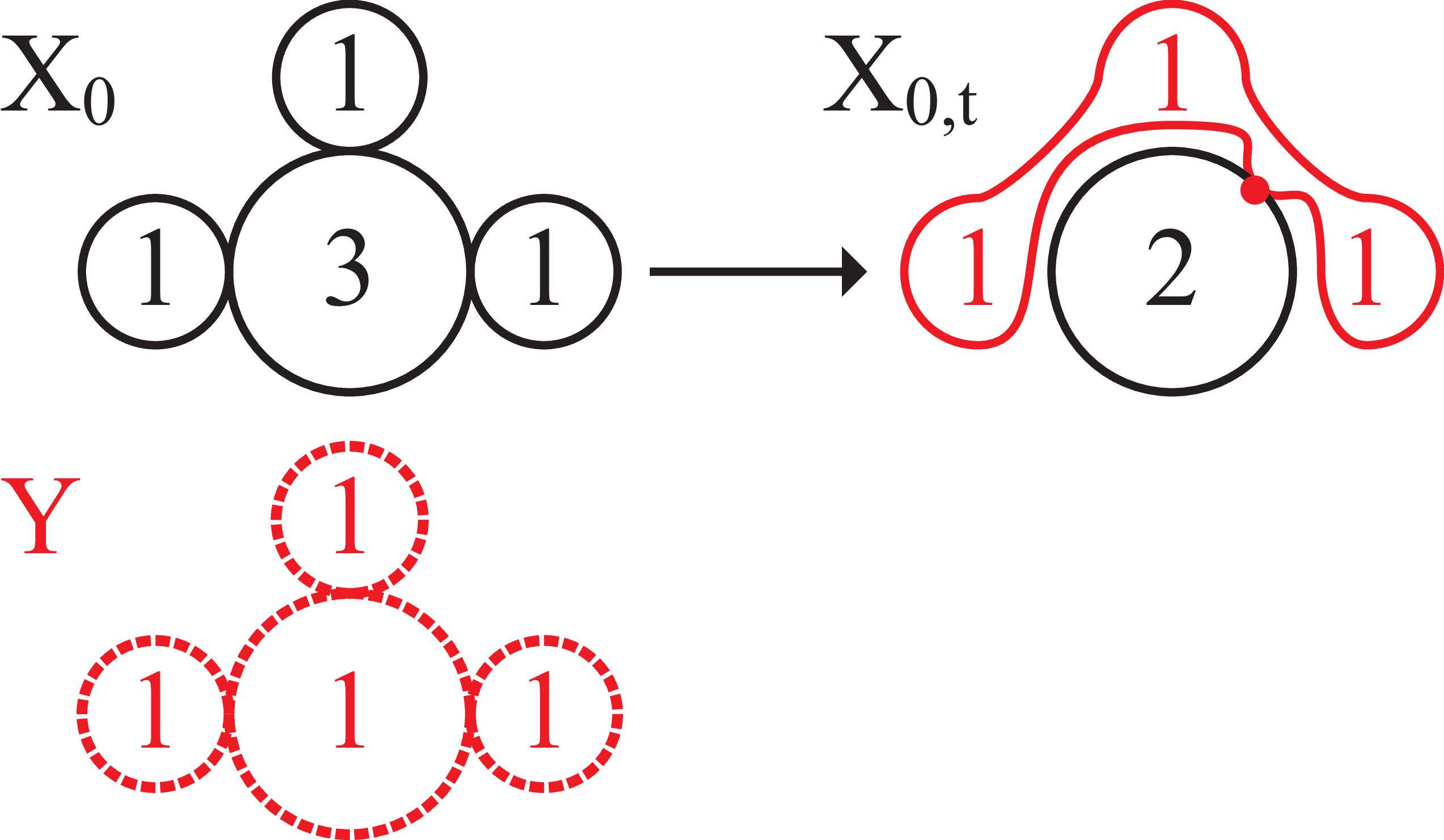}\end{center}

$\boldsymbol{[IV^*.1]} \ IV^* \barkarrow I_1^*$
\begin{center}\includegraphics[width=8cm,clip]{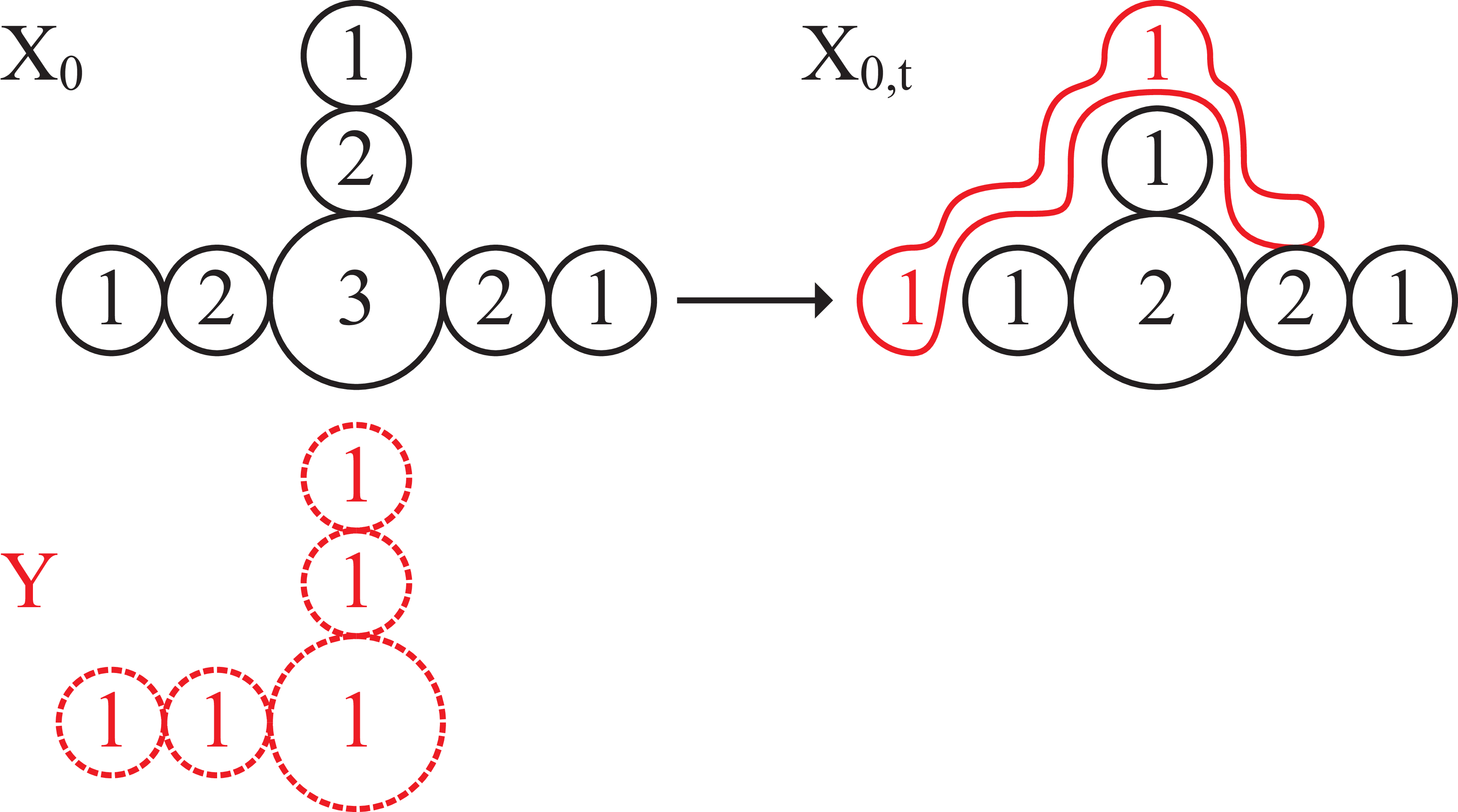}\end{center}

$\boldsymbol{[IV^*.2]} \ IV^* \barkarrow I_0^*$
\begin{center}\includegraphics[width=8cm,clip]{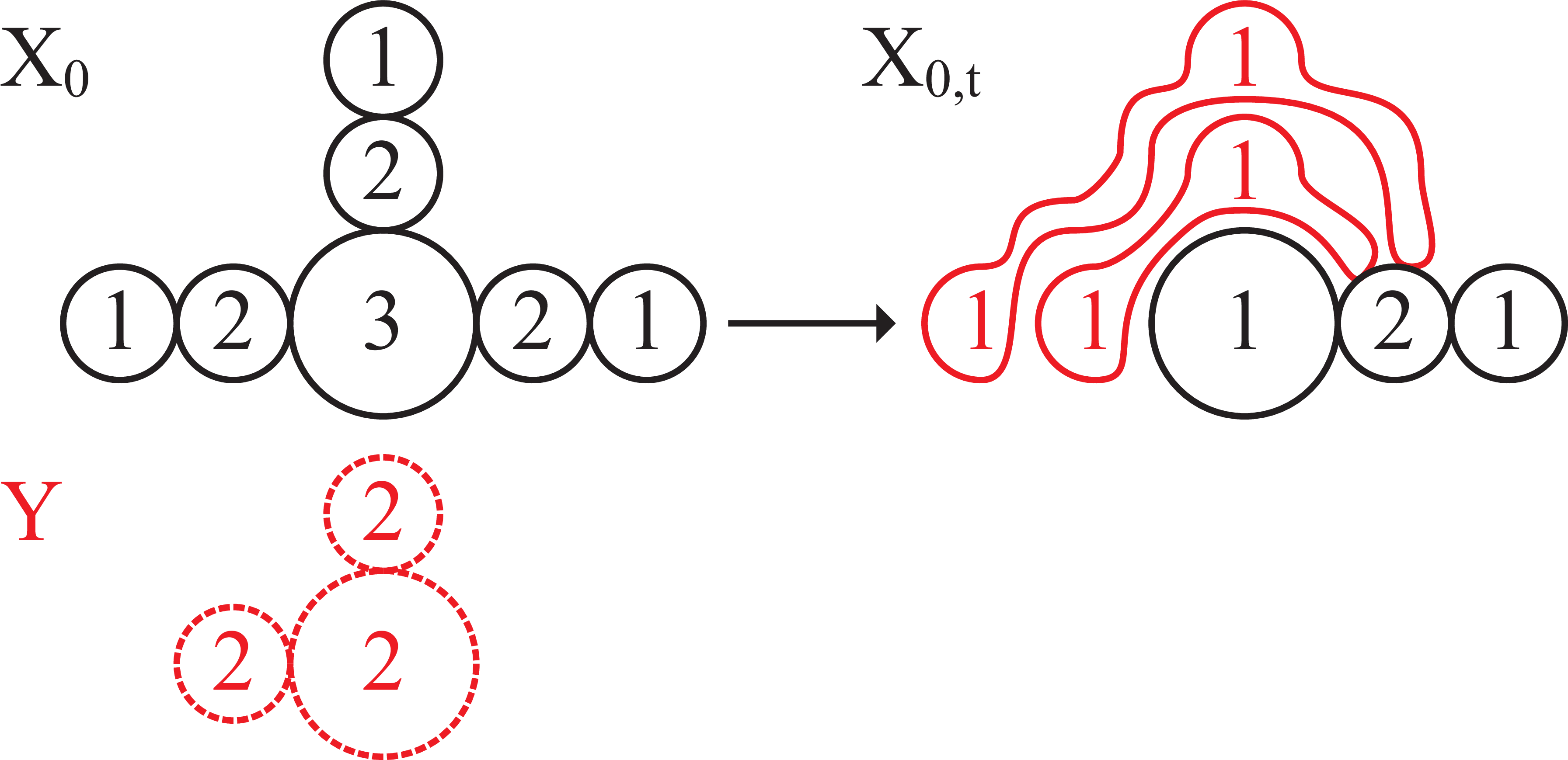}\end{center}
\newpage

$\boldsymbol{[IV^*.3]} \ IV^* \barkarrow I_6$
\begin{center}\includegraphics[width=8cm,clip]{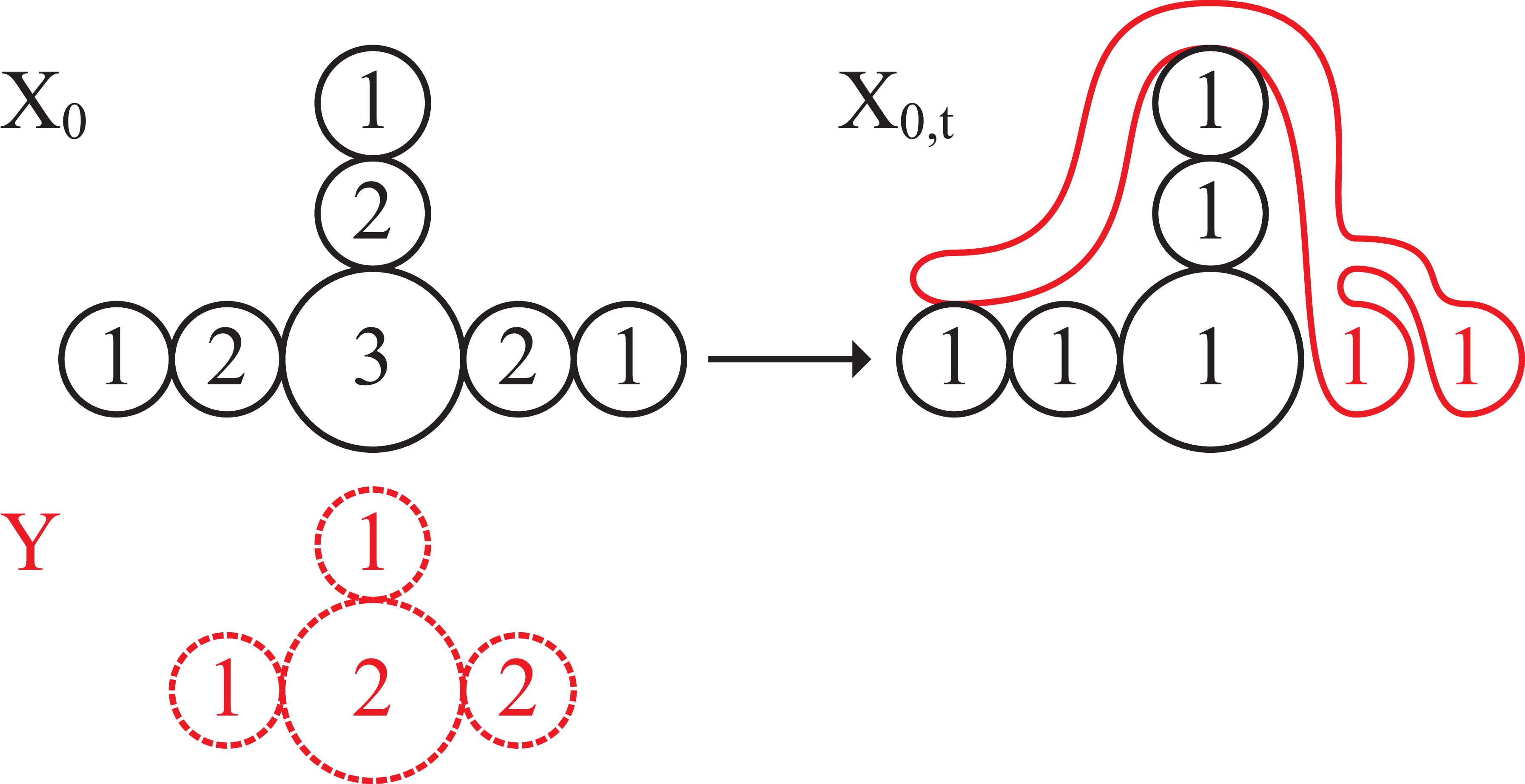}\end{center}

$\boldsymbol{[IV^*.4]} \ IV^* \barkarrow I_1^*$
\begin{center}\includegraphics[width=8cm,clip]{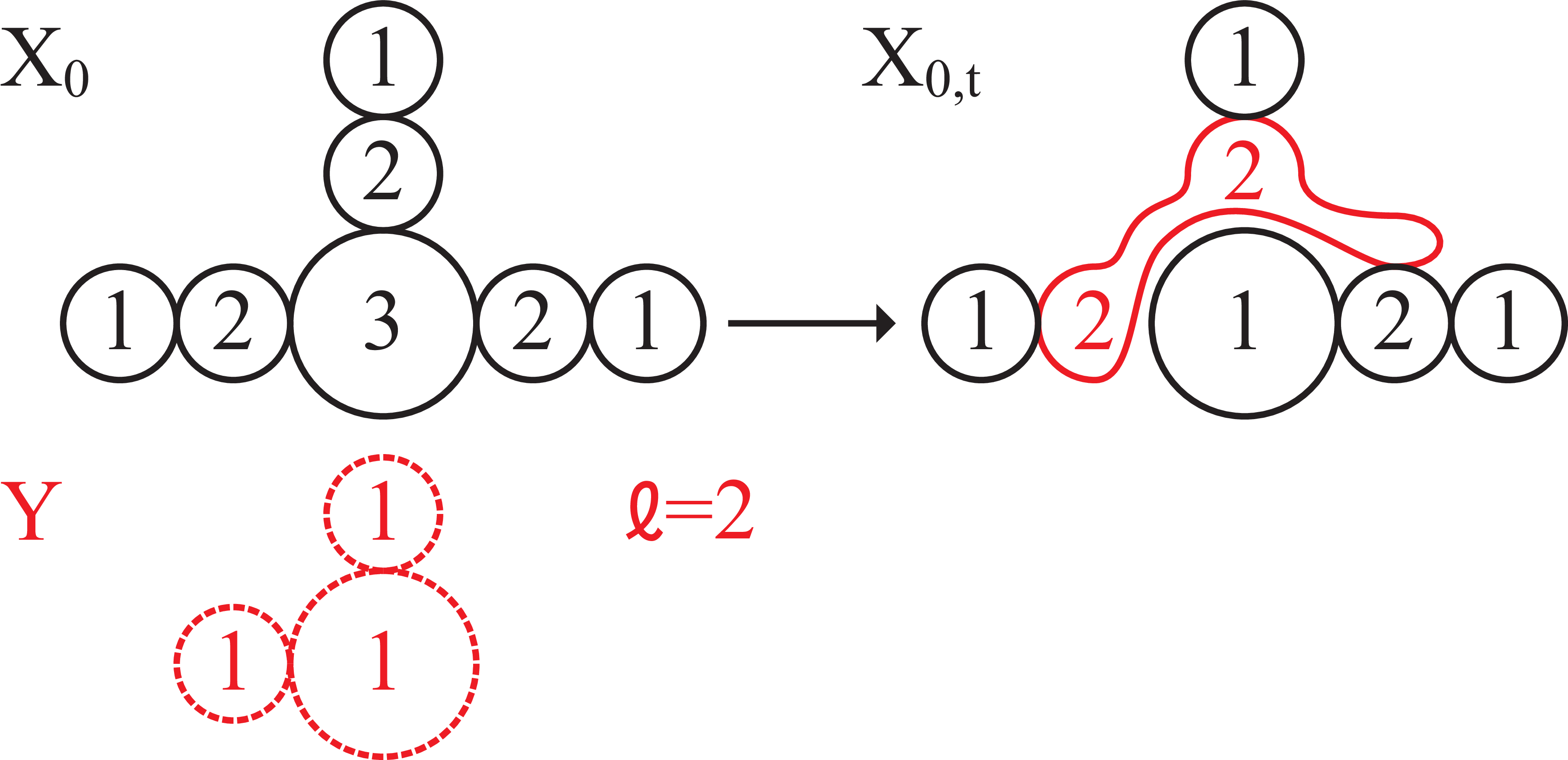}\end{center}

$\boldsymbol{[I_0^*.1]} \ I_0^* \barkarrow I_4$
\begin{center}\includegraphics[width=5cm,clip]{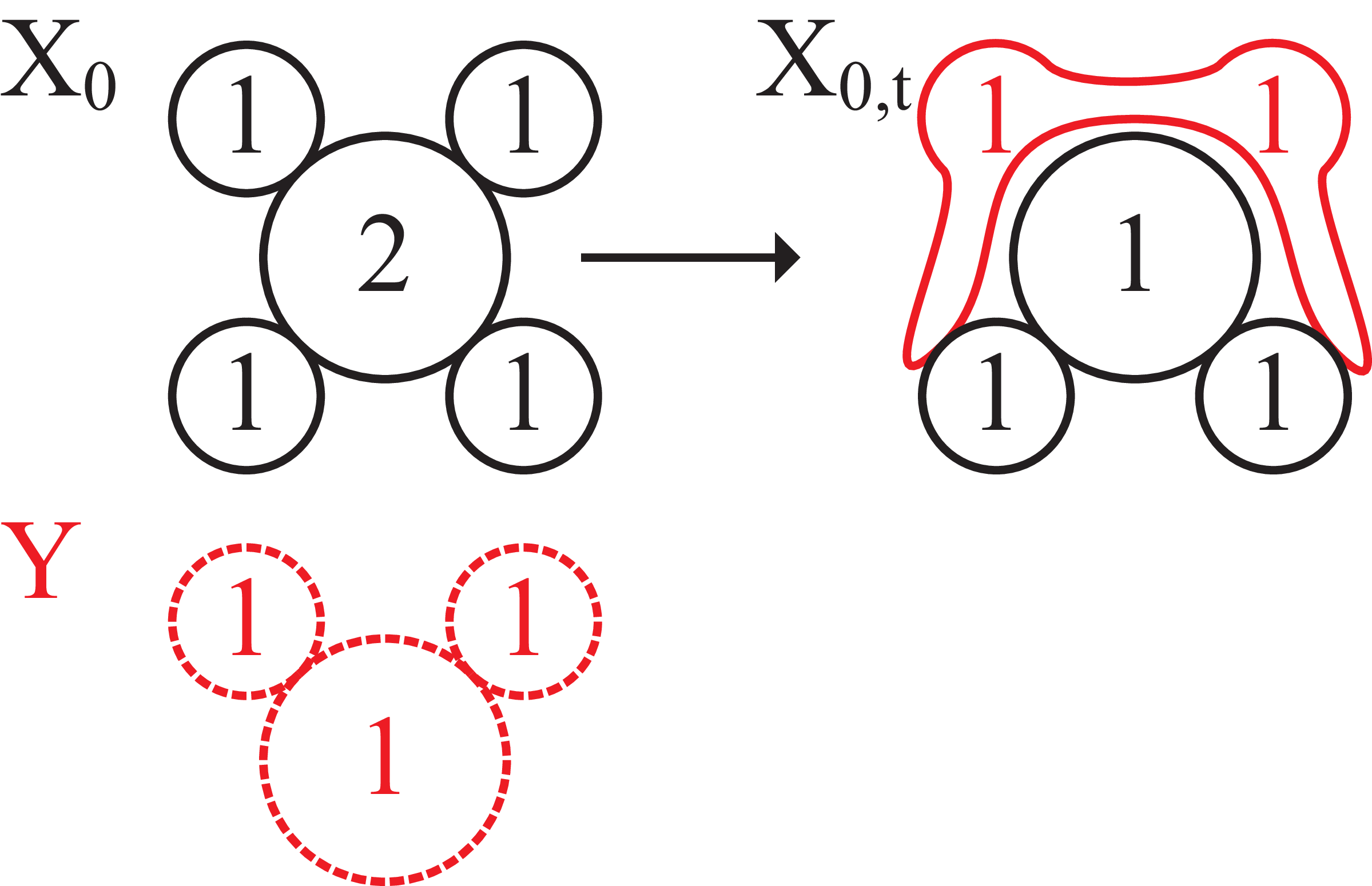}\end{center}

$\boldsymbol{[I_0^*.2]} \ I_0^* \barkarrow I_3$
\begin{center}\includegraphics[width=5cm,clip]{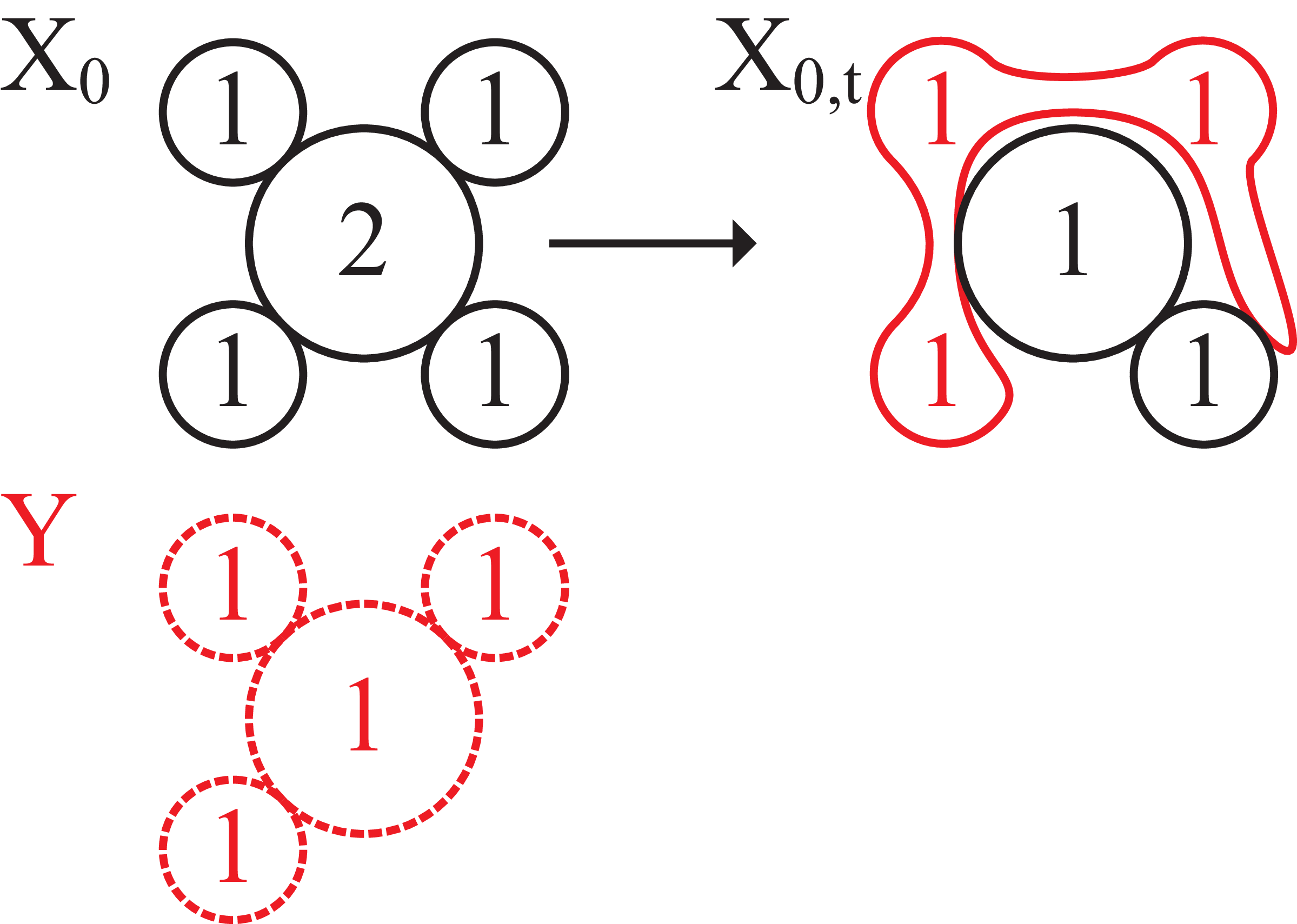}\end{center}
\newpage

$\boldsymbol{[I_n^*.1]} \ I_n^* \barkarrow I_{n-1}^*$ \\
\begin{center}\includegraphics[width=10cm,clip]{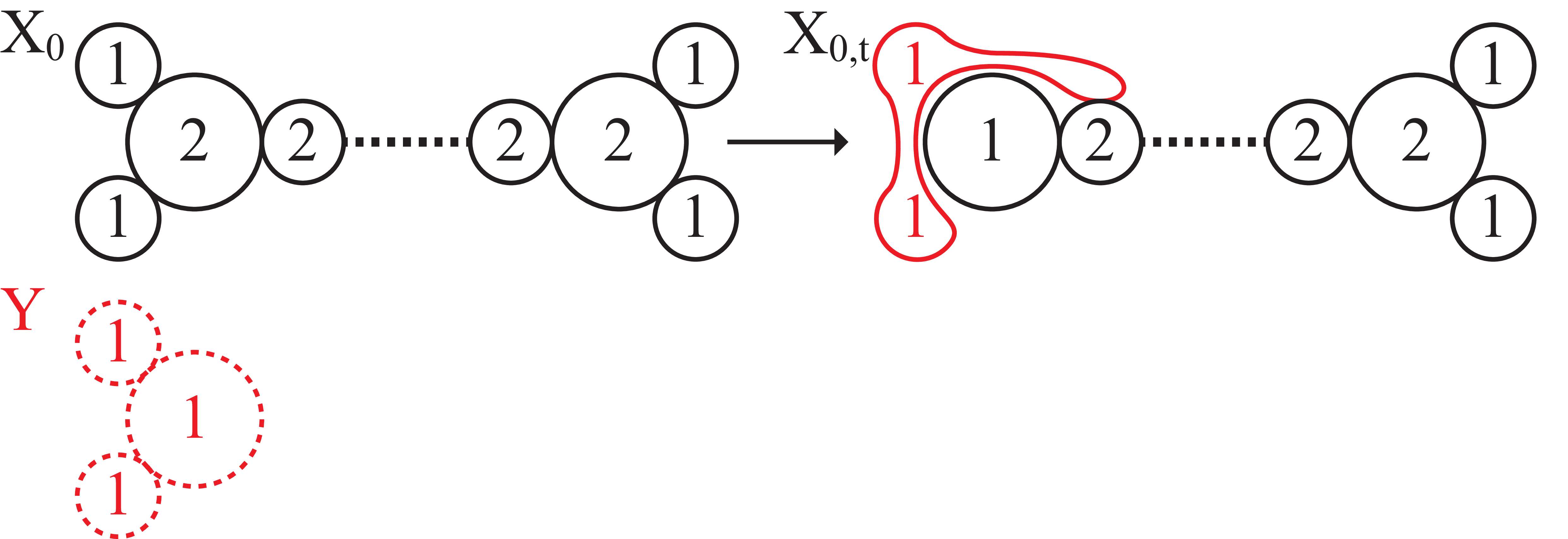}\end{center}

$\boldsymbol{[I_n^*.2]} \ I_n^* \barkarrow I_{n+4}$ 
\begin{center}\includegraphics[width=10cm,clip]{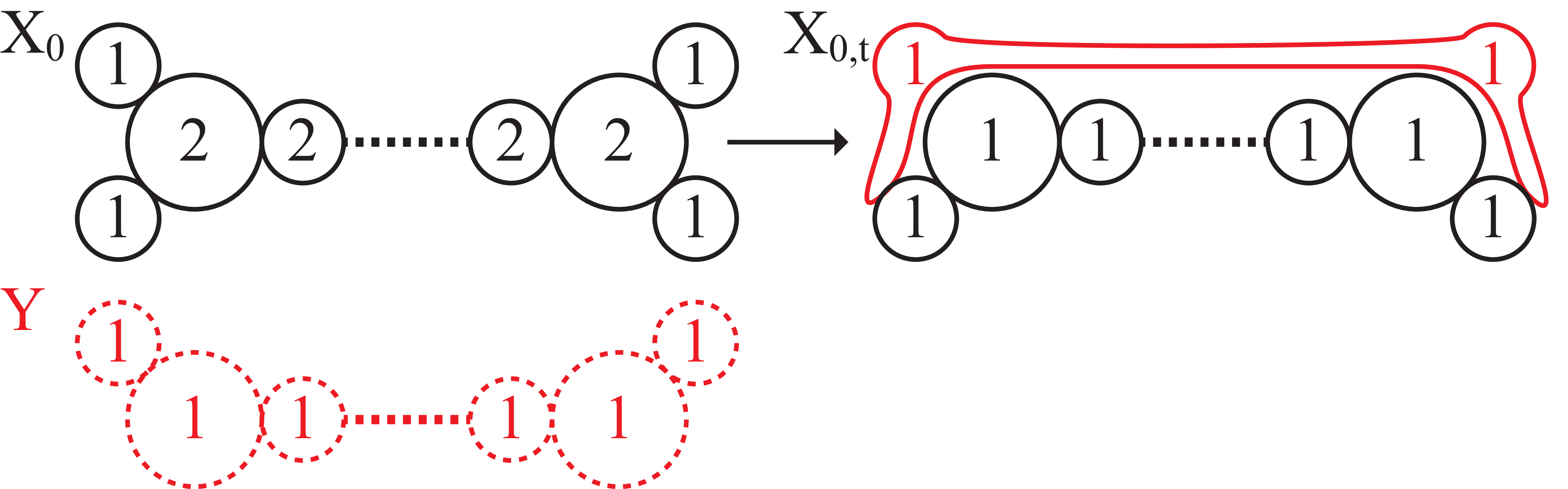}\end{center}

\begin{rem}
\begin{description}
\item[\rm (a)]
Takamura \cite{Ta3} introduced 
not only a barking family associated with \emph{one} simple crust 
(which we reviewed in Section \ref{sec:Ta}) 
but also a barking family associated with \emph{several} crusts. 
The latter is called a \emph{compound barking family}. 
Note that 
the barking families 
$\boldsymbol{[II^*.6]}$, 
$\boldsymbol{[II^*.7]}$, 
$\boldsymbol{[III^*.7]}$, 
$\boldsymbol{[III^*.8]}$ 
are compound barking families. 

\item[\rm (b)]
The singular fiber 
$I_n^* \, (n \geq 1)$ 
is \emph{constellar} (constellation-shaped), 
that is, 
it is obtained by bonding stellar singular fibers. 
So 
$\boldsymbol{[I_n^*.1]}$ 
and 
$\boldsymbol{[I_n^*.2]}$ 
are barking families 
for constellar case 
rather than 
for stellar case. 
See \cite{Ta3} for details. 

\item[\rm (c)]
This list contains no barking families for a degeneration with the singular fiber
$m I_n$. 
In fact, 
for 
$m I_n \, (m \geq 2)$, 
we use another method to construct a splitting family, 
which splits 
$m I_n$ 
into 
$m I_{n-1}$ 
and 
$I_1$. 
See \cite{Ta1} for details. 
\end{description}
\end{rem}


{\fontshape{sc}\selectfont
\noindent
Takayuki OKUDA 

\noindent
Graduate School of Mathematics, Kyushu University, \\
Motooka 744, Nishi-ku, Fukuoka 819-0395, Japan
}

\noindent
\textit{E-mail address}: {\fontfamily{cmtt}\selectfont {\text t-okuda@math.kyushu-u.ac.jp}}


\end{document}